\documentclass{article}
\usepackage[utf8]{inputenc}
\usepackage[russian,english]{babel}
\usepackage{amsmath,amsfonts,amssymb,mathrsfs,amscd,comment,latexsym}

\usepackage[matrix,arrow,curve]{xy}

\usepackage[usenames]{color}
\usepackage{colortbl}

\newtheorem{theorem}{Theorem}

\newtheorem{lemma}{Lemma}[section]
\newtheorem{sublemma}{Sublemma}[lemma]
\newtheorem{corollary}{Corollary}

\newtheorem{definition}{Definition}

\newtheorem{remark}{Remark}

\newcommand{\steparhead}[1]{%
 \vskip 5pt\par\noindent\emph{\textcolor{blue}{(#1)}} \par\vskip 3pt}

\newcommand*{\hm}[1]{#1\nobreak\discretionary{}%
{\hbox{$\mathsurround=0pt #1$}}{}}

\newcommand{\aff}{{\mathbb{A}}}

\newcommand{\bFN}{{\mathscr{F}_{Nis}}}
\newcommand{\bO}{\mathscr{O}}
\newcommand{\bF}{\mathscr{F}}
\newcommand{\BiMod}[3]{{\phantom{.}_{#1}} {{#2}_{#3}}}
\DeclareMathOperator{\coker}{coker}

\newcommand{\bZ}{\mathbb{Z}}

\usepackage[matrix,arrow,curve]{xy}

\newcommand{\pri }[1]{ { {#1}^\prime } }
\newcommand{\ppri }[1]{ { {#1}^{\prime\prime} } }
\newcommand{\pppri }[1]{{{#1}^{\prime\prime\prime}}}

\newcommand{\affl}{\mathbb{A}^1}
\newcommand{\prl}{\mathbb{P}^1}
\newcommand{\pro}{{\mathbb{P}}}

\newcommand{\bL}{\mathscr{L}}

\newcommand{\ovWCor}{\overline{WCor}}

\newcommand{\Zwtr}{\mathbb{Z}_{wtr}}

\newcommand{\ShNW}{ShNisWtr}

\newcommand{\bHom}{{\mathcal Hom}}

\begin{document}
\title{
  Triangulated category of 
    effective $Witt$-motives 
  $DWM^-_{eff}(k)$
  \thanks{Research is supported 
    by Chebyshev Laboratory 
      (Mathematics and Mechanics Faculty, 
       St.Petersburg State University), 
    by the Russian Science Foundation grant No.14-21-00035, 
    by "Native towns", 
      a social investment program of PJSC "Gazprom Neft", 
and by the RFBR grant 14-001-31095}
  }
\date{January 20, 2016}
\author{A.~E. Druzhinin\footnote{Chebyshev Laboratory, St. Petersburg State University, 14th Line, 29b, Saint Petersburg, 199178 Russia, andrei.druzh@gmail.com}} 

%
%


\maketitle


\begin{abstract}  
The category of effective $Witt$-motives $DWM^-(k)$ 
  for perfect field $k$, $char k\neq 2$,
  with functor $WM\colon Sm_k\to DWM^-(k)$
    defining motives of smooth affine varieties  
is constructed.
In the construction
Voevodsky-Suslin method
is applyed to a category
of $Witt$-correspondence
between affine smooth varieties 
$WCor_k$
  that morphisms are defined by 
    class in $Witt$-group
    of quadratic space $(P,q_P)$ with 
    $P$ being $k[X\times Y]$-module finitely generated projective over $k[X]$
    and 
    $q_P\colon P\to Hom_{k[X]}(P,k[X])$ being $k[X\times Y]$-liner isomorphism.
And
  the natural isomorphism
  $$Hom_{DWM^-_{eff}(k)}(WM(X),\bF[i]) \simeq H^i_{Nis}(X,\bF) $$  
  for any smooth affine $X$ and homotopy invariant Nisnevich sheave $\bF$ with $Witt$-transfers 
  (that is presheave on the category $WCor_k$ such that its restriction on the category $Sm_k$ is a sheave) 
  is proved. 

\end{abstract}


\section{Introduction.}

This work is devoted to the problem of 
the construction of 
  the triangulated category of $Witt$-motives $DWM^-(k)$
by thes
  Voevodsky-Suslin method that was originally used for 
  the construction of the category of motives $DM^-(k)$. 

Let's explain what is meant under
  Voevodsky-Suslin method.
In \cite{VSF_CTMht} V. Voevodsky constructed 
the triangulated category of motives $DM^-(k)$
  for perfect field $k$.
The construction starts with 
  the category of correspondence 
    between smooth affine varieties 
  $Cor_k$, 
and proceeds by proving
  some properties of 
  the presheaves of abelian groups 
    on this category.
These properties express certain compatibility between
  the structure of the category $Cor_k$,
  the topological structure on the category $Sm_k$ and 
  the structure of the category with an interval specified by 
    the affine line $\affl$. 
Due to this compatibility 
the construction results in obtaining 
  the category of motives 
in sufficiently explicit form as 
  stabilization in the $\mathbb G_m^{\wedge 1}$ direction from 
  a full subcategory of a derived category.    
 
So we wish to construct the category of $Witt$-motives starting 
from some category of correspondences between smooth affine varieties. 
It was suggested by I.A.~Panin 
to use so called category of $Witt$-correspondence $WCor_k$.
It is additive category 
that objects are smooth affine varieties 
and $Hom$-groups are the $Witt$-groups of some exact categories. 
Namely, 
for two smooth affine varieties $X$ and $Y$, 
the quadratic spaces $(P,q_P)$ 
defined by 
the $k[X\times Y]$-module $P$ 
  finitely generated projective over $k[X]$
and 
$k[X\times Y]$-linear symmetric isomorphism 
$q_P\colon P\simeq Hom_{k[X]}(P,k[X])$ 
are regarded.

\vspace{5pt}

In current text
following described above approach
we define 
the category of effective $Witt$-motives $DWM^-_{eff}(k)$
  as 
    full subcategory of the derived category $D^-(\ShNW_k)$
    of the category of Nisnevich sheaves with $Witt$-transfers,
    consisting of 
    the motivic complexes,
    i.e. complexes $A^\bullet$ 
      whose Nisnevich sheaves of cohomologies 
        $\underline{h}_{Nis}^i(A)$
      are homotopy invariant
  and
  define
  the functor $WM\colon Sm_k\to DWM^-_{eff}(k)$
    of $Witt$-motives of smooth affine varieties
  $$WM(X) = 
  \{U\mapsto 
    (\dots\to WCor_{Nis}(U\times\Delta^i,X)
     \dots\to WCor_{Nis}(U\times\Delta^1,X)
     \to WCor_{Nis}(U,X))\}  
  $$ where $WCor_{Nis}(-,Y)$ denotes
       Nisnevich sheafication of the presheave $WCor(-,Y)$   
          
We prove
{\par\noindent Theorem A.\em
  The category $DWM^-_{eff}(k)$ is equivalent to 
  the localization of the derived category $D^-(\ShNW_k)$ 
    by the morphisms 
    corresponding to the projections $X\times \affl\to X$
}\par    
and 
{\par\noindent Theorem B.\em  
  There is natural isomorphism
  \begin{equation}\label{bpWM}
  Hom_{DWM^-_{eff}(k)}(WM(X),\bF[i]) \simeq H^i_{Nis}(X,\bF)
  \end{equation}
  for any
   smooth affine varieties $X$ and 
   homotopy invariant sheaves with $Witt$-transfers $\bF$. 
}

For this purpose we 
  strengthens the result of \cite{AD_WtrSh}
    about preservation of homotopy invariance
          under Nisnevich sheafication
          for presheaves with $Witt$-transfers 
  and prove 
{\par\noindent Theorem C.\em    
  Nisnevich sheafication $\bF_{Nis}$ of homotopy invariant presheave
    with $Witt$-transfers $\bF$
  is strictly homotopy invariant  
}\par
and prove 
{\par\noindent Theorem D.\em
  For any presheave with $Witt$-transfers $\bF$
  the Nisnevich sheafication $\bF_{Nis}$ and 
  presheaves of cohomologies $H^i_{Nis}(\bF_{Nis})$
  are equipped with $Witt$-transfers in canonical way     
}
that implies that 
  the category $\ShNW$ 
  is abelian 
and allows to consider 
  the category $D^-(\ShNW)$. 
  
Let's note that
to give
the definition of 
  the category $DWM^-_{eff}(k)$ and functor $WM$
it is enough to prove the results on the 
  preservation of $Witt$-transfers and homotopy invariance only
  for Nisnevich sheafication of the presheave.  
But to prove described properties of this objects
it requires 
strictly version of this results, 
i.e. to prove similar properties of 
 presheaves of Nisnevich cohomologies.
(In fact to prove the isomorphism \eqref{bpWM} we use 
 more stronger statement 
 (theorem \ref{ReprChWtr})
 then theorem D.)     

Now we describe the contents of the text.
In sections \ref{sec_relAffEx} and \ref{sec_ndimEtEx}
  two excision isomorphisms are proved.
Namely 
  Zariski excision isomorphism
  on the affine line over local base
  $$\frac{\bF(\affl_U-0_U)}{\bF(\affl_U)}
    \simeq
    \frac{\bF(V-0_U)}{\bF(V)} 
  $$
  in the section \ref{sec_relAffEx}  
  and 
  etale excision isomorphism
  $$\frac{\bF(X^\prime-Z)}{\bF(X^\prime)}
    \simeq
    \frac{\bF(X-Z)}{\bF(X)} 
  $$ 
  in the section \ref{sec_ndimEtEx}.

In the section \ref{sec_DHI}
this excision isomorphisms are applied to prove 
that homotopy invariant Nisnevich sheaf with $Witt$-transfers
is strictly homotopy invariant.
In combination with the result proved in \cite{AD_WtrSh}
this implies that 
  Nisnevich sheafication 
  of homotopy invariant presheave with $Witt$-transfers
  is strictly homotopy invariant.  
  
The section \ref{secWtrNTop}
is devoted to analyse of the behaviour of $Witt$-transfer 
in relation to the Nishevich topology. 
In particular we prove theorem D.
  
In the section \ref{sec_DWMeff} 
the construction of the 
  category of effective $Witt$-motives with
  the functor defines $Witt$-motives of 
    smooth affine varieties
  is given.
And the basic property \ref{bpWM} 
  of $Witt$-motives for smooth affine schemes
is proved.

The author expresses gratitude for I.A.~Panin for 
formulation of the problem and consultations at its decision.

\section{The $Witt$-correspondence category}

In this section  
 the definition and basic properties of 
 category $Witt$-correspondense
   introduced in \cite{AD_WtrSh} 
are given.

At first let's take next technically useful definition.

\begin{definition}[$Proj(p)$]$\phantom{d}$
\label{P_S_U}
For any morphism of affine schemes $p\colon S\to U$
we define exact category with duality 
$Proj(p)$
 as follows.
The category $Proj(p)$ is
  full subcategory 
    of the category $k[S]$-mod 
  consisting of modules 
    that are finitely generated projective over $k[U]$.  
To define the duality 
we consider the functor 
$$\begin{aligned}
D_p\colon 
  k[S]-mod\to k[S]-mod\\
  M \mapsto D_p(M)=Hom_{k[U]}(M,k[U])
\end{aligned}$$
(where 
  the structure of $k[S]$-module on $D_p(M)$
  is defined by the rule: 
    $(f\cdot \rho)(m)=\rho(f\cdot m)$
       for $\rho \in D_p(M)$ and $f \in k[S]$).
Since
  the functor $Hom_{k[U]}(-,k[U])$
  is the duality 
  on the exact category $Proj(k[U])$ 
      of finitely generated projective $k[U]$-modules 
and since 
  the forgetful functor $Proj(p)\to Proj(k[U])$ is exact
,  
$D_p$ induce on $Proj(p)$ 
the structure of exact category with duality.        

Also we will denote by $Proj(X,Y)$
the category $Proj(pr_X)$ 
  where $pr_X\colon X\times Y\to X$ is canonical projection.

\end{definition}
\begin{remark}
\par 
The scalar restriction  along 
  a morphism $\pri S{\xrightarrow{f}} S$ 
induce a functor
$f_*\colon Proj(p\hm\circ f)\to Proj(p)$
respecting duality
for any morphism $S\stackrel{p}{\rightarrow} U$. 
\par 
The base change (or reverse image functor) 
  along the morphism $u\colon U^\prime\to U$
induce the functor of categories with duality
  $u^*\colon Proj(p)\to Proj(p^\prime)$
for any morphism $p\colon S\to U$ and
$p^\prime = p\times_U U^\prime 
  \colon S\times_U U^\prime\to U$.
\end{remark}

Using this definition the additive category $WCor_k$ 
can be defined as follows. 

\begin{definition}[$WCor_k$]$\phantom{d}$
\begin{itemize}
\item[$\Diamond$] 
  $Ob\,WCor_k\,\hm=\,Ob\,Sm_k$;
\item[$\Diamond$] 
  $WCor_k(X,Y) \hm= W(Proj(X,Y), D_{pr})$ where 
    $W(Proj(X,Y), D_{X,Y})$ is 
      the Witt group of exact category with duality 
        (with respect to definition from \cite{BW}).
        
  Typical example of morphism from the variety $X$ to $Y$ 
    is quadratic space $(\BiMod{k[Y]}{P}{k[X]}, \phi)$ where 
      $\BiMod{k[Y]}{P}{k[X]}$ is $k[Y\times X]$-module 
        finitely generated projective as $k[X]$-module and
      $$\phi: P \to Hom_{k[X]}(P, k[X])$$ is
        $k[X\times Y]$-linear isomorphism. 
\item 
  The composition of the morphisms 
  $\Phi\in WCor_k(X,Y)$ and $\Psi\in WCor_k(Y,Z)$ 
    is defined by tensor product of corresponding quadratic spaces.
\item 
{\sloppy
  The identity morphism is defined by diagonal. 
  It means that
    $id_X$ for any variety $X$ is defined by
    the bi-module $\BiMod{k[X]}{k[X]}{k[X]}$ and 
    canonical isomorphism $k[X] \hm\simeq Hom_{k[X]}(k[X],k[X])$
      regarded as quadratic form on it.                  

}
\end{itemize}
\end{definition}

\begin{definition}[Presheaves and sheaves with $Witt$-transfers]\label{DeShT}
  Abelian groups presheaf with $Witt$-transfers is
  a functor $F\colon WCor_k\to Ab$
    satisfying additivity condition on disjoint unions :
    $\bF(X_1\coprod X_2)=\bF(X_1)\oplus\bF(X_2)$ for any $X_1$ and $X_2$.
  The Sheaf with $Witt$-transfers is 
    a presheaf with $Witt$-transfers
    that becomes a sheaf after restriction to $Sm_k$.  
\end{definition}

\newcommand{\kX}{{k[X]}}
\newcommand{\kY}{{k[Y]}}
\newcommand{\kZ}{{k[Z]}}
\newcommand{\yPx}{\BiMod{\kY}{P}{\kX}}

\newcommand{\zQy}{\BiMod{\kZ}{Q}{\kY}}
\newcommand{\Dual}[2]{Hom_{#1}(#2,#1)}

\newcommand{\kW}{{k[W]}}
\newcommand{\xPw}{\BiMod{\kX}{P}{\kW}}
\newcommand{\yQx}{\BiMod{\kY}{Q}{\kX}}

\newcommand{\zRy}{\BiMod{\kZ}{R}{\kY}}

\begin{remark}
\label{i_Sm_WCor}
There is a functor $\;i: Sm_k\to WCor_k$
  that maps 
    the regular map $f\colon X\to Y$ to 
    the morphism defined by 
      the module $\BiMod{k[Y]}{k[X]}{k[X]}$ and 
      canonical isomorphism 
        $\BiMod{k[Y]}{k[X]}{\phantom{.}}\hm\simeq 
          Hom_{k[X]}(\BiMod{k[Y]}{k[X]}{\phantom{.}},(k[X]).$
Moreover 
  the composition of a morphism with the image of a regular map at the right or left 
is equal to 
  direct or reverse image 
    of corresponding quadratic space 
  along the $f$:
$$\begin{aligned}
f\circ w=red_{f\times W}(w),& w\in WCor_k(W,X),\\
w\circ f=ind_{f}(w),& w\in WCor_k(Y,Z).\end{aligned}$$

\newcommand{\yPw}{\BiMod{\kY}{P}{\kW}}
\newcommand{\wPx}{\BiMod{\kW}{P}{\kX}}

\end{remark}
\begin{remark}
For any smooth variety $X$
$$WCor(X,pt) = W(X)$$
and there is an diagonal embedding
$$W(X)\to WCor(X,X).$$
So 
the functor of $Witt$-groups is 
a presheave with $Witt$-transfers
and
any presheave with $Witt$-transfers is 
a presheaf of modules over the Witt ring. 

Also let's note that 
for any invertible function $\lambda\in k[X]^*$
we will denote by 
  $(\lambda)_X \in QSpace(Proj(X,X))$ 
  quadratic spaces
    of rank one 
    on a diagonal
    with quadratic form defined by $\lambda$,
    i.e.
    $(\BiMod{k[X]}{k[X]}{k[X]},\lambda)$.
And 
  for by regular map $f\colon X\to Y$
let's denote by 
  $(\lambda)_f$ or $(\lambda)_{\Gamma_f}$ 
  quadratic space in $Proj(X,Y)$
   of rank one 
   that support is graphic $\Gamma_f$ of map $f$      
   and that quadratic form defined by $\lambda$   
   i.e.
     $(\BiMod{k[Y]}{k[\Gamma_f]}{k[X]},\lambda)$.
\end{remark}

\begin{remark}[Essential smooth schemes] 
\label{rem_essmWCor}
Sometimes we will suppose 
  presheaves with $Witt$-transfers 
  and
  bi-functor $WCor_k(-,-)$  
 to be well defined
  on essential smooth schemes over $k$.

The category $WCor_k$ can be defined
on essential smooth schemes, 
i.e. it is possible to define 
$Witt$-correspondence 
  between such schemes
agreed with the definition of 
  $Witt$-correspondence between smooth affine schemes,
by considering of essential smooth schemes as pro-objects.
However,
  in fact 
  we will use only 
  germs of presheaves on essential smooth schemes  
  and 
  $Witt$-correspondence between smooth schemes over 
    fields of functions $k(X)$ of smooth varieties
    and
    local rings $\mathcal O_{X,x}$ at smooth point $x$.
So we will describe in details only this two constructions.

1)
The germ of the presheave $\bF$ 
on any essential smooth variety $U$ is 
  $$\bF(\varprojlim\limits_{i\to\infty} V_i) = 
  \varinjlim\limits_{i\to\infty} \bF(V_i).$$
And any morphism $g\colon U\to X$, $X\in Sm_k$
induces homomorphism $\bF(g)\colon \bF(X)\to\bF(U)$
agreed with composition with any morphism $f\colon X\to X^\prime$
  $\bF(f\colon g) = \bF(g)\colon \bF(f)$.

2) 
For any smooth affine variety $S$ and a point $s\in S$
there are functors
$$\begin{aligned}
  PreWtr_k&\to PreWtr_{k(S)}\\
  \bF\; &\longmapsto \bF_{k(S)}\\
  PreWtr_k&\to PreWtr_{\mathcal O_{S,s}}\\
  \bF\; &\longmapsto \bF_{k[S]_s} \end{aligned},$$
such that 
 $$\bF_{k(S}(X) \simeq \bF(X)$$
 ($\bF(X)$ at the right side is
  the germ of $\bF$ on $X$ considered as
    essential smooth scheme over $k$),
and such that
  for any morphism $f\colon X\to Y$ of smooth schemes 
    over open subscheme $S^\prime\subset S$ 
  section group homomorphism 
    $$\bF_{k(S)}(f_{k(S)})\colon \bF_{k(S)}(X)\to \bF_{k(S)}(Y) $$
  coincides with germ group homomorphism induced by $f$,
and similar for $\bF_{k[S]_s}$.   

To define these functors it is enough 
to construct for any 
  $Witt$-correspondence $\Phi\in WCor_{k(S)}(X,Y)$
  ($\Phi\in WCor_{k[S]_s}(X,Y)$) 
  between smooth affine $k(S)$-schemes ($k[S]_s$ schemes),
some 
  $Witt$-correspondence $\Phi^\prime\in WCor_{k}(X^\prime,Y^\prime)$
  between smooth affine $S^\prime$-schemes, 
  for some open subscheme $S^\prime\subset S$,
that goes to $\Phi$ under base chance 
    along the morphism $k[S^\prime]\to k(S)$ ($k[S^\prime]\to k[S]_s$).
To construct such $\Phi^\prime$ 
we can 
choose algebraic parametrisation 
  of the data defining $X$, $Y$ and $\Phi$%
, i.e. parametrise it 
  by finite set of algebraic parameters from $k(S)$ ($k[S]_s$), 
  satisfying some algebraic equations.
And then set $S^\prime\subset S$ to be
   open subscheme
   on that all parameters are well defined.
     
\end{remark}

It is useful 
  for proofs of excision isomorphisms   
  to extend the category $WCor_k$ to pairs 
    of smooth scheme and its open subscheme.

\begin{definition}[Category $WCor^{\cdot\hookrightarrow\cdot}$]
The objects of the category  
  $WCor_k^{\cdot\hookrightarrow\cdot}$ 
  are 
  the pairs $(X_1,X_2)$ 
    of smooth affine scheme $X_1$ and 
    its open subscheme $X_2$.  
Any morphism
    $\Phi\in WCor_k^{\cdot\hookrightarrow\cdot}
      ((X_1,X_2) \to (Y_1,Y_2))$  
  is a pair of morphisms 
    $\Phi_i\in WCor(X_i,Y_i),\;i=1,2$
  such that
    $ \Phi_1\circ i_X = i_Y \circ\Phi_2 $. 

It is equivalent to say that 
  $WCor_k^{\cdot\hookrightarrow\cdot}$
  is full subcategory 
  of the category of arrows of the category $WCor_k$
  consists of the open embeddings.
\end{definition}
\begin{definition}[Category $WCor^{pair}$]
Additive category $WCor_k^{pair}$ 
is a factor-category 
  of the additive category 
  $WCor_k^{\cdot\hookrightarrow\cdot}$
  by the ideal generated by identity morphisms
  of the objects $(X,X)$ for all varieties $X$.
\end{definition}
\begin{remark}\label{rem_WCorpair}
The last definition is equivalent to say that 
$WCor_k^{pair}$
is full subcategory 
  of homotopy category $\mathcal K(WCor_k)$
  of additive category $WCor_k$
consists of the complexes 
  concentrated in 2 adjoin degrees 
  with differential homomorphisms being 
    open embedding $X_2\to X_1$.

More precisely it means that
$Hom$-groups are defined as follows
  \begin{multline*}
  WCor_k((X_1,X_2),(Y_1,Y_2)) \stackrel{ref}{=}\\
  H(
    WCor_k(X_1,Y_2)
      \xrightarrow{i_Y \circ - , - \circ i_X}
    WCor_k(X_1,Y_1) \oplus WCor_k(X_2,Y_2)
      \xrightarrow{- \circ i_X , i_Y \circ -}
    WCor_k(X_2,Y_1)
  ),  
   \end{multline*} 
  where $H$ denotes 
     homology group in middle term of 
     the complex of the length 3.

{\sloppy 
So morphism $\Phi\colon (X_1,X_2) \to (Y_1,Y_2)$ 
  is defined by pair of morphisms 
    $\Phi_i\in WCor(X_i,Y_i),\;i=1,2$
  such that 
    left diagram is commutative
  and 
  the pair $(\Phi_1,\Phi_2)$ defines zero morphism 
    if and only if there exists 
      $\Xi\in WCor_k(X_1,Y_2)$
    such that right diagram is commutative.    
$$\xymatrix{
X_2\ar@{^(->}[r]^{i_X}\ar[d]^{\Phi_2}&X_1\ar[d]^{\Phi_1}\\
Y_2\ar@{^(->}[r]^{i_Y}&Y_1\\
}\phantom{GGGGG}
\xymatrix{
X_2\ar@{^(->}[r]^{i_X}\ar[d]_{\Xi\circ i_X}&X_1\ar[d]^{i_Y\circ \Xi}\ar[ld]^{\Xi}\\
Y_2\ar@{^(->}[r]^{i_Y}&Y_1,\\
}$$

}In terms of quadratic spaces
  to construct the morphism $\Phi\colon (X_1,X_2)\to(Y_1,Y_2)$
  in $WCor_k^{pair}$
it's enough
  to find the quadratic space  
\begin{equation}\label{dPM}
(P,q_P)\in Proj(pr_{Y_1\times X_1\to X_1}) 
:\quad
k[Y_2]\otimes_{k[Y_1]}
  P
    \otimes_{k[X_1]}k[X_2] \simeq 
      P
        \otimes_{k[X_1]}k[X_2].
\end{equation}
Because
  such $(P,q_P)$ defines the morphism in $WCor_k(X_1,Y_1)$.
And 
  due to isomorphism \eqref{dPM} 
  the module $P\otimes_{k[X_1]}k[X_2]$ 
    has the canonical structure 
    of the module over $k[Y_2]\otimes k[Y_1]$.
Hence 
  the quadratic space $i_X^*(P,q_P) $ 
  in fact lies in  $Proj(pr_{Y_2\times X_2\to X_2})$
  and defines the morphism in $WCor_k(X_2,Y_2)$.   

It is noteworthy that
  not every morphism has such representation 
  but 
  all used in the constructions morphisms have.   
\end{remark}
     
\begin{remark}\label{rem_PreshPair}
For any presheaf with $Witt$-transfers $\bF$
one can define 
  the presheaf $\bF^{pair}$ on the category of pairs 
  such that $$\bF^{pair}(X_1,X_2) = \frac{\bF(X_2)}{\bF(X_1)}.$$ 
\end{remark}

The rest part of the section is devoted to 
homotopy invariant presheaves with $Witt$-transfers.

Let's define the category $\ovWCor_k$.
\begin{definition}[Category $\ovWCor_k$]\label{H0W}
$\ovWCor_k$ is 
a factor-category of additive category $WCor_k$
such that 
$$\overline {WCor_k}(X,Y)\hm=
  \coker(
    WCor_k(\affl\hm\times X,Y)
    \xrightarrow{(-\circ i_0)-(-\circ i_1)}
    WCor_k(X,Y)
   ),
$$ 
  where 
  $i_0,i_1\colon{X}\hookrightarrow\affl\hm\times{X}$ 
  are zero and unit sections.
  
And similar we define factor-category $\overline{WCor_k^{pair}}$
of the category $WCor_k^{pair}$  
\end{definition}

\begin{remark}[Homotopy invariance]
\item[1)] 
By definition homotopy invariant presheaf with $Witt$-transfers is                
  a presheaf with $Witt$-transfers that is homotopy invariant. 
So such presheaves are exactly 
  the presheaves on the category $\ovWCor_k$.  
And this equivalence is agreed with 
  continuation of presheaves to essential smooth schemes
  in sense of the remark \ref{rem_essmWCor}.
\item[2)]  
Let's take used later  descriptions
  of the pairs of homotopy morphisms. 

Two morphisms $\Phi_1,\Phi_2\colon X\to Y$ in $WCor_k$
  represented by the spaces $(P_1,q_1)$ and $(P_2,q_1)$
becomes equal in $\ovWCor_k$
if and only if 
exists the quadratic space $(H,q)$, 
such that
  the equalities
  $$[j_0^*(H,q_H)]=[(P_0,q_0)] ,\quad
    [j_1^*(H,q_H)]=[(P_1,q_1)]$$  
  holds in $Proj(pr_{X\times Y\to Y})$. 
                 
For coincidence 
  of two morphisms of pairs in $\overline{WCor_k^{pair}}$  
  represented by spaces in sense of the remark \ref{rem_WCorpair}
it is enough the existence of  
the quadratic spaces 
  $(H,q)\in Proj(pr_{Y,\affl\times X})$ and
  $(G_0,\pri q_0),(G_1,\pri q_1)\in Proj(pr_{Y_1,X_1}))$,
such that
  \begin{gather*}
  k[Y_2]\otimes_{k[Y_1]}H\otimes_{k[X_1]}k[X_2] = 
     H\otimes_{k[X_1]}k[X_2]\\
  [{j_i}^*(H,q)] = [(P_i,q_i)\oplus(G_i,\pri q_i)],\\
  k[Y_2]\otimes_{k[Y_1]}G_i= G_i,\quad i=0,1.
  \end{gather*}
\end{remark}


\section{Excision on $\affl_U$.}\label{sec_relAffEx}

This section is devoted to 
  some particular case of
  excision isomorphism 
  in Zariski topology 
  on the relative affine line $\affl_U$ 
    over the local essential-smooth scheme $U$. 
Namely 
  it's the case of excision 
  outside the zero section
  of $\affl_U$.

\begin{theorem}\label{locrelAfZarEx}
Let 
  $\bF$ is be homotopy invariant sheave with $Witt$-transfers
  and
  $U=Spec\, \mathcal O_{X,x}$ be spectre of the local ring 
    of a smooth variety $X$ at any point $x$.
Then 
  for any Zariski open subvariety $V\subset\affl_U$
    containing zero section $0_U\subset \affl_U$
  restriction homomorphism   
    $$i^*\colon \frac{\bF(\affl_U-0_U)}{\bF(\affl_U)}\hm\to \frac{\bF(V-0_U)}{\bF(V)}$$
  is an isomorphism 
    ($i$ denotes embedding of $V$ into  $U$).  
\end{theorem}

\begin{remark} 
\item{1)}
In fact 
In terms of the remark \ref{rem_PreshPair}
theorem statement means that 
  $i^* \colon 
    \bF^{pair}(\affl_U-0_U,\affl_U) \to \bF^{pair}(V-0_U,V)$
is isomorphism
and follows from the following lemma.
\item{2)} 
Obviously theorem is equivalent to the statement
  that
    restriction homomorphism 
      $$
      i^*\colon 
        \frac{\bF(V-0_U)}{\bF(V)} 
        \hm\to 
        \frac{\bF(V^\prime-0_U)}{\bF(V^\prime)}
      $$
    is an isomorphism 
    for any two embedded open subvarieties containing zero section:
      $$0_U\subset V^\prime \subset V \subset\affl_U.$$
\end{remark}

\begin{lemma}\label{arW_ExrelU} 
Let $i$ denotes embedding of $V$ into $\affl_U$.
Then its class
  $[i] \in \overline {WCor_k}((V,V-0_U) , (\affl_U,\affl_U-0_U))$
  is isomorphism (as a morphism in $\overline {WCor_k}$).
\end{lemma}

\begin{remark}\label{relZarExDiag}
$\phantom{f}$\begin{list}{}{\leftmargin 0ex \itemsep   3pt}

\item[\em{a)}] 
The right inverse morphism to $[i]$  in $\overline{WCor_K}$ is 
a morphism 
$$
\Phi \hm\in 
    WCor_K( (\affl_U,\affl_U\hm-0_U) , (V,V\hm-0_U) )
:\quad
[i\circ\Phi] \hm= [id] 
  \hm\in 
  \overline {WCor_K}( (\affl_U\hm-0_U) , (\affl_U\hm-0_U) )
$$  
And it is equivalent to existence of 
\begin{equation}\label{HE_i_relU}
\Theta \hm\in 
  WCor_K(
     (\affl_U\hm\times\affl,(\affl_U\hm-0_U)\hm\times\affl) ,
     (\affl_U,\affl_U\hm-0_U)
        )
:\quad
\Theta \circ j_0 = i \circ \Phi ,\; 
\Theta \circ j_1 = id,
  \end{equation} 
  (where 
   $j_0,j_1 \colon 
     (\affl_U,\affl_U\hm-0_U) 
     \to 
     (\affl_U,\affl_U\hm-0_U)\times\affl$ 
   denotes 
     embeddings of zero and unit sections  
   respectively).
   
\vspace{0pt}

\item[\em{b)}] 
The left inverse to $[i]$ in $\overline{WCor_K}$  
is an morphism  
$$\Psi \hm\in 
    WCor_K( (\affl_U , \affl_U\hm-U_0)) , 
            (V , V\hm-U_0) )
:\quad
  [\Psi \circ i] = [id]
    \in 
    \overline{WCor_K}( (V,V\hm-0_U),
                       (V,V\hm-0_U) ).
$$
It means the existence of 
\begin{equation}\label{HE_s_relU}
\Xi \hm\in 
    WCor_K( (V\hm\times\affl, (V\hm-0_U) \hm\times \affl) ,
            (V,V-0_U) ) 
:\quad
  \Xi \hm\circ j_0  \hm = \Psi\hm\circ i ,\;
  \Xi \hm\circ j_1 \hm = id
\end{equation}
  ($j_0,j_1 \colon (V,V\hm-0_U)\to (V,V\hm-0_U)\times\affl$ 
      denotes 
      zero and unit sections respectively). 
\end{list}
\end{remark}

\emph{Proof of the lemma \ref{arW_ExrelU}}
\par
We will construct required morphisms in the category  
  $WCor_U$ 
  using it's embedding into $WCor_k$ 
  extended in sense of remark \ref{rem_essmWCor}.

\par {\em a)}

\steparhead{Description of quadratic spaces} 
To construct right inverse to 
  $i\in \overline{WCor} ((V,V-0_U),(\affl_U,\affl_U\hm-0_U))$
means 
to find quadratic spaces $(P,q_P)$ and $(H,q_H)$ 
  corresponding to the morphisms 
  $\Phi$ and $\Theta$ from remark \ref{relZarExDiag}.a).
In the terms of the quadratic spaces 
the properties of $\Phi$ and $\Theta$ 
means following: 

\begin{list}{}{\leftmargin 1ex \itemsep   3pt}

\item[\em 1)] 
  $P\in k[V\hm\times_U \affl_U]-mod$, 
  $P$ is finitely generated projective over $k[\affl_U]$ 
  and
  $q_P\colon P\hm\simeq Hom(P,k[\affl_U]),$
  is $k[V\hm\times_U \affl_U]$-linear symmetric isomorphism;
\item[\em 2)] 
  $H\in k[\affl_U\hm\times_U \affl_U\hm\times\affl]-mod$, 
  $H$ is finitely generated over $k[\affl_U\hm\times\affl]$ 
  and
  $q_H\colon H\hm\simeq Hom(H,k[\affl_U \hm\times \affl]),$
  is $k[\affl_U\hm\times_U \affl_U\hm\times\affl]$-linear 
     symmetric isomorphism;
\item[\em 3)] 
  canonical homomorphisms
  \begin{equation}\begin{array}{c}\label{PM_i_relU}
  P\otimes_{k[\affl_U]}k[\affl_U-0_U]\hm\to 
    k[V-0_U]\otimes_{k[V]}P\otimes_{k[\affl_U]}k[\affl_U-0_U],\\
  H \otimes_{k[\affl_U]}k[\affl_U-0_U]\hm\to 
    k[\affl_U-0_U]\otimes_{k[\affl_U]}H
    \otimes_{k[\affl_U]}k[\affl_U-0_U]							
  \end{array}\end{equation} 
  (one of structure  of $k[\affl_U]$-module 
   is regarded here as left and another as right)
  are an isomorphisms 
  (it means that spaces $(P,q_P)$ and $(H,q_H)$ 
   defines the morphisms of pairs); 
\item[\em 4)] 
  in the Witt group of the category 
  $Proj(\aff^2_U\to\affl_U)$ 
  following equalities holds:
  \begin{equation}\begin{aligned}\label{WHE_i_relU}
  \big[(H,q_H)
    \otimes_{k[{\affl_U}\hm\times\affl]}
    k[\affl_U\hm\times 0]\big] \hm=
      \big[ \phantom{.}_{k[\affl_U]}{k[V]}\hm
      \otimes_{k[V]}
      (P,q_P)\big],
   \\ 
   \big[(H,q_H)
     \otimes_{k[{\affl_U}\hm\times\affl]}
     k[\affl_U\hm\times 1]\big] \hm=
       \big[E_{k[\affl_U]}\big],
  \end{aligned}\end{equation}
  (it equivalent to the equalities \eqref{HE_i_relU}). 
\end{list}

\steparhead{Definition of modules $P$ and $H$}
Let's now give one lemma (sublemma 3.1.1 from \cite{AD_ShNisWtr})
used in the proof of excision isomorphisms. 
\begin{lemma}\label{surcon}
Let 
  $X$ be projective scheme over netherian ring,
  $Z$ be closed subscheme,  
  $\mathscr{F}$ be a coherent sheave and 
  $\bL$ be very ample bundle on $X$. 
Then 
  for all $n$ larger some $k$ 
  the restriction 
  $
  \Gamma(
    \mathscr{F}\otimes \bL^{\otimes n}
    )
    \to 
  \Gamma(
    (\mathscr{F}\otimes\bL^{\otimes n})\big|_{Z}
    )
  $ is surjective.
\end{lemma}
And let's involve one used definition-denotetion.
\newcommand{\bI}{\mathcal I}
\begin{definition}[sheaf of iseals $\bI(s)$ and subscheme $Z(\bI)$]
\label{def_ZeroSubs-SectIdeal}
\newcommand{\shI}{\mathcal I}
Any non-zero global section
  $s\in \Gamma(\bL)$ 
    of some invertable sheaf $\bL$ 
    on any irreducible scheme $X$
  defines 
    a sheaf of ideals in $\bO(X)$ 
    isomorphic to the sheaf $\bL^{-1}$
    (%
     The section $s$ defines homomorphism $\bO(X)\to \bL)$,
     And multiplying it on the sheaf $\bL^{-1}$
     we get homomorphism $\bL^{-1}\to\bO(X)$.   
     Or equivalently it is the sheaf of functions
     $\{f\in \bO(X),\,div_0\,f-div_0\,s < 0$.%
     ) 
  We will denote this sheaf of ideals by $\bI(s)$.
Also for any sheaf of ideals $\bI\subset \bO(X)$ 
  we will denote by $Z(\bI)$
  corresponding closed subscheme of $X$. 
And finally we denote 
  the closed subscheme $Z(\bI(s))$ by $Z(s)$
  for any section of line bundle $s$. 
\end{definition}

Let's identify
  $\affl_U\hm\times_U \affl_U$ and $V\hm\times_U \affl_U$ 
  with subsets of $\prl_{\affl_U}$ 
  $$
  \affl_U\hm\times_U \affl_U = \affl_{\affl_U} 
    \subset \prl_{\affl_U},\,
  V\hm\times_U \affl_U = V_{\affl_U} 
    \subset \prl_{\affl_U},\,
  $$    
Let $T$ and $D$ be its complement
and let
  $\Delta\subset \prl_{\affl\times U}$ be 
    the graph of embedding $\affl_U\hookrightarrow \prl_U$.
  $$
  T=\prl_{\affl_U}\setminus \affl_{\affl_U},\, 
  D=\prl_{\affl_U}\setminus V_{\affl_U}
  \Delta\subset \prl_{\affl\times U} = 
     \Gamma(\affl_U\hookrightarrow \prl_U).
  $$ 
Then let' choose a sections 
  $$
  \nu,\delta\in\Gamma(\prl_{\affl_U},\bL(T))\colon
  div_0\,\nu = 0_{\affl_U} ,\; 
  div_0\,\delta = \Delta  
  .$$ 
We can assume in addition that
  $$\nu\big|_T = \delta\big|_T.$$
Indeed
since $T$ is the infinity section of $\prl_{\affl\times U}$ 
the fraction 
  $u=\frac{\nu\big|_T}{\delta\big|_T}$ 
can regarded as the function on $\affl\times U$
and 
since 
intersections of 
  zero divisors of $\nu$ or $\delta$ 
  with $T$ 
  are empty 
$u$ is invertible function.
So 
  if we multiply $\delta$ 
  by the inverse image of $u$ along the projection 
  $$\prl_{\affl\times U}\to \affl\times U,$$
we don't change zero divisor of $\delta$
and 
make the required equality holds.

Subscheme $V\subset  \prl_{\affl\times U}$ contains 
  zero section $0_{\affl\times U}$.
Hence 
  the intersection 
    of $D$ with $0_{\affl\times U}$ 
  is empty
and 
  by the sublemma \ref{surcon} 
  for sufficiently large $n$ 
  there is a section 
  $$s_0\in 
    \Gamma(\prl_{\affl_U},\bL(nT)):\; 
    s_0\big|_D=\nu^n,\, 
    s_0\big|_{0_{\affl_U}}=\delta^n
  .
  $$
\label{correct}
Let 
  $s = 
   s_0\cdot (1-t)+ \delta^n\cdot 
   \in 
   \Gamma(\prl_{\affl_U} \times\affl ,\bL(nT))
  .$
Then 
  $$
  s\big|_{0_{\affl}\times\affl} = \delta^n,\,
  s\big|_{T_{\affl}\times\affl} = \nu^n ( = \delta^n)
  .
  $$
Since
  $s$ is invertible on $T_{\affl}$ and
  $s_0$ is invertible on $D$, 
they defines closed subschemes 
$$
S_0 = Z(s_0)/ \subset V_{\affl_U} ,\quad 
S   = Z(s)/ \subset {\affl}_{\affl_U}\times\affl,
$$
(see definition \ref{def_ZeroSubs-SectIdeal} for $Z(s)$.) 
Thus we can put 
$$ 
  P \hm= 
  k[S_0]_{k[V\times_U \affl_U]} 
  ,\quad 
  H \hm= 
  k[S]_{k[\affl_U\times_U \affl_U \times\affl]}.
$$

\steparhead{Checking of point 3,i.e. equalities \eqref{PM_i_relU}} 
The condition \eqref{PM_i_relU} holds 
because 
  $s_0\big|_{0\times\affl\times U} = \delta^n$ and
  $s\big|_{0\times\affl\times U\times\affl} = \delta^n$
and hence 
  $$S_0\cap 0\times\affl\times U=0\times 0\times U$$ and 
  $$S\cap (0\times\affl\times U\times\affl) = 
   0\times 0\times U\times\affl.$$

\steparhead{Definition of quadratic forms $q_P$ and $Q_H$}
To prove that
    $P$ and $H$ 
  are finitely generated projective over 
    $k[\affl\times\affl_{U\times\affl}]$ and  
    $k[V\times \affl_{U}]$               respectively
and
to construct 
  quadratic forms $q_P$ and $q_H$ 
it is useful 
  to consider 
    the morphism of projective 
    smooth schemes over $\affl\times{U\times\affl}$  
    $$
    \overline F=([s:\mu^n]\colon 
      \prl_{\affl_U\times \affl}\to 
      \prl_{\affl_U\times \affl}.
    $$
The morphism $\overline F$ 
is finite 
because it is projective and quasi-finite 
and it is flat  
because it is finite morphism of  
  essential smooth schemes
  of the same dimension.

Next we consider 
  base change of $\overline F$ 
  along the embedding 
  $\affl_{\affl_U\times \affl}\hookrightarrow
    \prl_{\affl_U\times \affl}$.
Since $div\,\mu = T$
we get 
  the morphism of affine schemes 
  $$F = \frac{s}{\mu^n} \colon \
    \affl_{\affl_U\times \affl} \to 
           \affl_{\affl_U\times \affl}$$ 
  that is also finite and flat.
  
Let's denote by 
  $B_A$ 
  the algebra corresponding to the morphism $F$, 
  it means that 
    both of $B$ and $A$ as $k$-algebras 
      are isomorphic to 
      the function algebra of relative affine line
      $k[\affl_{\affl\times U\times \affl}]$
    and 
    homomorphism $A\to B$ 
      is induced by the function $\frac{s}{\mu^n}$.
Let's fix the trivialisation 
  of the canonical class of $\affl\times\affl\times U\times\affl$. 
Then using the proposition 2.1. from \cite{OP}
we get 
  $k[B]$-linear  isomorphism $$q_B\colon k[B]\simeq Hom_A(k[B],k[A]).$$ 

Now let's consider 
  commutative diagram with Cartesian squaes 
$$\xymatrix{
  V\times_U\affl_U\ar[r] & \affl_U\times_U\affl_U\times\affl & &
  \\
  S_0\ar[u]^{e_0}\ar[r]\ar[d] & 
    S\ar[u]^{e}\ar[r]\ar[d] & 
      \affl\times\affl_{U\times \affl}
        \ar[r]\ar[d]_F                 & 
      \prl\times\affl_{U\times \affl}
        \ar[d]_{\overline F}           \\
  \affl_U\ar[r]^{j_0}               & 
  \affl_{U}\times\affl\ar[r]^z & 
    \affl\times\affl_{U}\times\affl\ar[r] & 
    \prl\times\affl_{ U}\times \affl
}$$
where
  $$e_0\colon  ,\quad e$$
  denotes embedding of closed subschemes $S_0$ and $S$
  and $d_0$ and $z$ defines zero section along  
    affine lines at left and right sides respectively 
  $$j_0 = \affl_U\times 0 ,\quad z = 0\times \affl_U .$$
Then let's apply to 
quadratic form $q_B$
base changes along $z$ and $j_0$
and
restrictions of scalars along $e$ and $e_0$
and define
quadratic forms 
$$\begin{aligned}
& q_{S} = z^*(q_B) 
  \colon  K[S]\simeq Hom_{}(K[S],K[\affl_U\times\affl])    &,\; 
& q_{S_0} = {j_0}^*(q_{S}) 
  \colon  k[S_0]\simeq Hom_{K[\affl_U]}(K[S_0],K[\affl_U]) \\
& q_{P} = {e_0}_*(q_{S_0}) &,\;
& q_H = {e}_*(q_{S})       . 
\end{aligned}$$ 

\steparhead{Checking of point 4,i.e. equalities \eqref{WHE_i_relU}}
The first equality of \eqref{WHE_i_relU} holds 
because functoriality of  
  restriction of scalars 
  in respect to base changes  
  $$
  i_*(q_P) = (i\circ e_0)_*(q_{S_0}) = 
  (i\circ e_0)_*({j_0}^*(q_S)) = 
  {j_0}^*(e)_*((q_S)) = {j_0}_0(q_H).
  $$

The second equality doesn't necessary true yet.
But it becomes true after 
 multiplication booth of the forms $q_H$ and $q_P$ 
 by some invertible function $\lambda\in k[\affl\times U]$.

In fact
let's denote   
$(H,q_H) \otimes_k k[\affl\times\affl\times U\times 1]$
by
$(H_1,q_{H_1})$
.
Then 
  by definition of $H = k[Z(s)]$
  so $H_1 = k[Z(s\big|_{\prl\times\affl\times U\times 1})]$
and thus  
$$
H_1 = H \otimes k[\affl\times\affl\times U\times 1] \simeq 
 k[\affl\times\affl\times U]/(\delta^n) 
$$ 
  ($\delta$ is regarded here as 
       function on $\affl\times\affl\times U$ 
     by trivialising of $\bL(T)$ 
       on this subscheme of 
       $\prl\times\affl\times U$).
Then since 
  $q_H$ is 
    $k[\affl\times\affl\times U \times \affl]$-linear 
    symmetric 
    isomorphism
  and  
  $q_{H_1}$
  is $k[\affl\times\affl\times U]$-linear 
    symmetric isomorphism, 
By
  the sublemma \ref{NilpLagr}
for any 
$k[\affl\times\affl\times U]$-linear quadratic form  $q$ 
   on $k[\affl\times\affl\times U]/(\delta^{n})$
ideals $(\delta^i)\hm\subset k[\affl\times\affl\times U]$
  are sublagrangian subspaces 
  for $i \leqslant n/2$.

\begin{sublemma}[%
sublemma 3.3.1 from \cite{AD_ZarShWtr}%
]
\label{NilpLagr}
Let 
  $B$ be $A$-algebra, and 
  $q\colon B\hm\simeq Hom_A(B,A)$ be 
  $B$-linear nondegenerate quadratic form 
  on $B$ over $A$. 
Then 
  for any ideal $I\subset B$ 
  its orthogonal $I^{\bot}\subset B$ 
    in respect to quadratic form $q$ 
  coincides with annihilator $Ann(I)\subset B$.
\end{sublemma}
     
So
  since $n$  is odd
  by sublagrangian reduction 
  the quadratic space 
    $(k[\affl\times\affl\times U]/(\delta^{n}),q_1)$
  is equal in Witt-group 
  to some one-ranged quadratic space 
  $(k[\affl\times\affl\times U]/(\delta^{n}),\lambda)$
  where $\lambda$ is invertible function of diagonal, 
  i.e. $l\in k[\affl_U]^*$. 

Then 
  if we multiply the forms $q_P$ and $q_H$ by $l^{-1}$, 
we 
  will not violate  
  first equality of \eqref{WHE_i_relU} and
  equalities \eqref{PM_i_relU}
and
  make  
  the class of
    $(H_1,q_{H_1})$
    in $WCor(\affl_U,\affl_U)$
  equal  
    to the class of $(1)_{k[\Delta]}$
  that is identity morphism of $\affl_U$.

\vskip 10pt
\steparhead{Point b)}
\par {\em b)} \label{affex_point-b}
\steparhead{Description of quadratic spaces}
To construct the left inverse to  
  $i\in \overline{WCor} ((V,V-0_U),(\affl_U,\affl_U\hm-0_U))$ 
means 
to find quadratic spaces $P$ and $H$ 
  corresponding to
   $\Psi$ and $\Xi$ 
   from remark \ref{relZarExDiag}.b). 
I.e. to to find following
\begin{list}{}{\leftmargin 1ex \itemsep   3pt}
\item[\em 1)] 
  $P\in k[V\hm\times_U \affl_U]-mod$ 
  finitely generated projective over $k[\affl_U]$ 
  and
  $k[V\hm\times_U \affl_U]$-linear symmetric isomorphism 
  $q_P\colon P\hm\simeq Hom(P,k[\affl_U]),$
\item[\em 2)] 
  $H\in k[\stackrel{1}V \hm\times_U \stackrel{2}V\hm\times\affl]-mod$ 
  finitely generated over $k[\stackrel{1}V \hm\times\affl]$  
  and
  $k[\stackrel{1}V \hm\times_U 
     \stackrel{2}V\hm\times\affl]$-linear symmetric isomorphism  
  $q_H\colon 
    H\hm\simeq 
    Hom_{k[\stackrel{2}V \hm\times \affl]}
      (H,k[\stackrel{2}V \hm\times \affl]),$
\item[\em 3)] canonical homomorphisms
\begin{equation}\begin{array}{c}\label{PM_s_relU}
P\otimes_{k[\affl_U]}k[\affl_U-0_U] \hm\to 
  k[V-0_U]\otimes_{k[V]}
    P\otimes_{k[\affl_U]}k[\affl_U-0_U]
,\\
H \otimes_{k[\stackrel{2}V]}k[V-0_U] \hm\to 
  k[V-0_U]\otimes_{k[\stackrel{1}V]} H 
  \otimes_{k[\stackrel{2}V]}k[V-0_U]
\end{array}\end{equation}
(one of structure  of $k[V]$-module is regarded as left and another as right)
are an isomorphism. 
(It means that spaces $(P,q_P)$ and $(H,q_H)$ defines the morphisms of pairs.) 
\item[\em 4)] 
in the Witt group of the category $Proj(V\times_U V\to V)$
  holds the equalities
  \begin{equation}\begin{aligned}\label{WHE_s_relU}
  \big[
    (H,q_H)
    \otimes_{k[\stackrel{2}V\times\affl]}
             k[\stackrel{2}V\times 0]
  \big]=
    \big[(P,q_P)\otimes_{k[\affl_U]}{k[V]}\big],
  \\ 
  \big[
    (H,q_H)
    \otimes_{k[\stackrel{2}V\hm\times\affl]}
             k[\stackrel{2}V\hm\times 1]\big]=
     \big[(1)_{V}\big],
  \end{aligned}\end{equation}
(This equalities are equivalent to the equalities \eqref{HE_s_relU}.)
\end{list}

\steparhead{Definition of modules P and H}
Since intersection of $\Delta$ with $D_V$ into $\prl_V $ is empty, 
  $\delta$ is invertible on $D_V$.  
Let's denote by $\delta^{-1} \in \Gamma(\bL(-T)_{D_V})$ it's inverse.
Next by sublemma \ref{surcon} 
  for sufficiently large  $n$
  there are exist the sections 
   $$\begin{aligned}
  \pri s\in \Gamma(\bL(n\cdot T),\prl_{\affl_U})&
     \pri s\big|_{D_{\affl_U}}=\nu^n,&&
         pri s\big|_{0_{\affl_U}}=\mu^{n-1}\cdot\delta
  &\\
   g\in \Gamma(\bL((n-1)\cdot T),\prl\times V)&
     g\big|_{D_V}   =\nu^n\cdot\delta^{-1},&
       g\big|_{\Delta}   =\mu^{n-1},&
         g\big|_{0_U\times V}=\mu^{n-1}
   &,\end{aligned}$$
because 
  intersection 
    of $D_{\affl_U}$ 
    with $0_{\affl_U}$ 
    into $\prl_{\affl_U}$ 
  is empty
  and
  intersection 
    of $D_V$ 
    with $\Delta$ 
    into $\prl_V$ 
  is empty too.  

Now we can define 
  the sections 
  $$\begin{aligned}
  &s_0 & \in \;&
      \Gamma(\bL(n\cdot T),\prl\times{V})
    \colon\quad  &
    s_0 = {{\prl}_{V\hookrightarrow {\affl_U}}}^*(\pri s) &\\
  &s_1 & \in \;& 
      \Gamma(\bL(n\cdot T),\prl_{V})
    \colon\quad  &
    s_1 = g\cdot\delta \in \Gamma(\bL(n\cdot T),\prl_{V}) &\\
  &s  & \in \;& 
      \Gamma(\bL(n\cdot T),\prl_{V\times\affl})
    \colon\quad &
    s=s_0\cdot (1-t)+s_1\cdot t
  &.\end{aligned}$$ 
    (where in firs line
     $\prl_{V\hookrightarrow \affl_U}$     
     denotes the embedding $\prl_V\hookrightarrow \prl_{\affl_U}$
     ) 
Then by definition of $s^\prime$, $g$, $s_0$, $s_1$ and $s$
  $$
  s_0\big|_{D_V} =
     s_1\big|_{D_V} =
       s\big|_{D_{V\times\affl}} =
          \nu^n,\; 
  s_0\big|_{0_V} =
     s_1\big|_{0_V} =
       s\big|_{0_{V\times\affl}} = 
          \mu^{n-1}\cdot\delta ,\;
  s_1\big|_{\Delta}=0 ,\; 
  div\,s_1 = div\,g \coprod  \Delta
  .
  $$ 
Then since 
  $s$ is invertible on $D_{V\times\affl}$ and 
  $s_0$ is invertible on $D_{\affl_U}$, 
 $$
 S_0 = Z(s_0) \subset V\times_U\affl_U,\;
 S = Z(s) \subset V\times_U V  \times\affl
 .$$
So we can put 
 $$
 P = k[S_0]_{k[V\times_U\affl_ U]},\;
 H = k[S]_{k[V\times_U V\times\affl]}
 .$$
 
\steparhead{Checking of the point 3, i.e. equalities \eqref{PM_s_relU}}
\label{parhead-checking-PMs_relU}
We should also 
  to  define quadratic forms 
but 
we can 
  check the equalities \eqref{PM_s_relU} just now 
because 
  it deals only with modules structures of $P$ and $H$.
To do it first of all note that
\begin{equation*}\begin{aligned}
P\otimes_{k[\affl_U]}k[\affl_U-0_U] 
&\simeq 
k[S_0 \times_{\affl_U} (\affl_U-0_U)]
,\\
k[V-0_U]\otimes_{k[V]}P\otimes_{k[\affl_U]}k[\affl_U-0_U]
&\simeq 	
k[(V - 0_U) \times_V S_0 \times_{\affl_U} (\affl_U-0_U)]
.
\end{aligned}\end{equation*} 
So first equality of \eqref{PM_s_relU}
  is equivalent to the equality
$$
S_0 \times_{\affl_U} (\affl_U-0_U) 
\simeq 
(V - 0_U) \times_V S_0 \times_{\affl_U} (\affl_U-0_U)
$$
that means that  
$$
0_U \times_V  S_0 \times_{\affl_U} (\affl_U-0_U)
 = \emptyset. 
$$
But
$$s_0\big|_{0\times{\affl_U}}=\delta$$
so
$$(0\times\affl_U) \cap S_0=0\times 0\times U.$$
Similarly 
second equality is equivalent to
$$0_U \times_V S \times_{V} (V - 0_U) = \emptyset,$$ 
but 
$$S\cap (0\times\affl_U\times\affl) = 0\times 0\times U\times\affl$$
because 
$s\big|_{0\times{\affl_U\times\affl}}=\delta$.

\steparhead{Definition of quadratic spaces}

Let's now check that 
  $P$ and $H$ are finitely generated projective 
  over $k[\affl_U]$ and $k[V\times_U\affl]$ respectively 
and
define the  quadratic forms 
  $q_P\colon P\simeq Hom(P,k[\affl_U])$ and 
  $q_H\colon H\simeq Hom(H,k[V\times\affl])$. 
It will be done by using the maps
  $$
  \overline{\pri F} = 
    [\pri s:\nu^n]\colon \prl_{\affl_U}\to\prl_{\affl_U},  \quad
  \overline{F_0} = 
    [s_0:\nu^n]\colon \prl_{V}\to\prl_{V},             \quad
  \overline{F} = 
    [s:\nu^n]\colon \prl_{V\times\affl}\to\prl_{V\times\affl}
  .$$
This maps 
  are finite surjective flat morphisms of essential smooth schemes
  and 
  are agreed by base changes, i.e
    $\overline{F_0}$ 
    coincides with  
    base changes of 
      $\overline{\pri F}$ and 
      $\overline{F}$ 
    along the embeddings 
      $i   \colon 
        V \hookrightarrow \affl_U$ and 
      $j_0 \colon 
        V \stackrel{id_V\times 0}{\hookrightarrow}
          {V\times\affl}$ 
    respectively.

\noindent

\newcommand{\Fp}{{F^\prime}}
\newcommand{\Ap}{{A^\prime}}
\newcommand{\Bp}{{B^\prime}}
\newcommand{\qp}{{q^\prime}}
\newcommand{\dpri}{{d^\prime}}

Next to 
  construct quadratic forms $q_P$ and $q_H$
  agreed in sense of \eqref{WHE_s_relU}
we apply to the maps 
    $\overline{F^\prime}$, $\overline{F_0}$ and $\overline{F}$
  the same construction 
    as in the point a) 
  simultaneously.

It means following. 
We consider 
  regular maps 
    $F^\prime$, $F_0$ and $F$ that are
    base changes 
    of $\overline{F^\prime}$, $\overline{F_0}$ and $\overline{F}$
     along the embeddings
     of affine lines into projective lines 
  $$
  {\pri F} = 
    \frac{\pri s}{\nu^n} \colon \affl_{\affl_U}\to\affl_{\affl_U},
      \quad
  {F_0} = 
    \frac{s_0}{\nu^n}\colon \affl_{V}\to\affl_{V},  
      \quad
  {F} = 
    \frac{s}{\nu^n}\colon \affl_{V\times\affl}\to\affl_{V\times\affl}
  .$$
Thus we get commutative diagram with Cartesian squares
$$\xymatrix{
&
 S\ar@_{_(->}[dl]\ar[dd]_p&
  &
   S_0\ar@_{_(->}[dl]\ar[dd]\ar@_{_(->}[ll]\ar@_{_(->}[rr]&
   &S'\ar@_{_(->}[dl]\ar[dd]\\
\affl_{V\times\affl}\ar[dd]_<<<<<{F}&
 &
  \affl_{V}\ar[dd]_<<<<<{F_0}\ar@_{_(->}[ll]\ar@_{_(->}[rr]&
   &
    \affl_{\affl_U}\ar[dd]_<<<<<{F'}&
     \\
&
 V\times\affl\ar@_{_(->}[ld]^{z}&
  &
   V\ar@_{_(->}[ld]^{z_0}
    \ar@_{_(->}[ll]^{ id_{V}\times 0}
    \ar[rr]_<<<<<<<{ i}&
     &
      \affl_U\ar@_{_(->}[ld]^{\pri z}\\
   \affl_{V\times\affl}&
 &
  \affl_{V}
    \ar@_{_(->}[ll]^{
     id_{\affl}\times id_{V}\times 0}
    \ar[rr]_{id_{\affl}\times i}&
   &
    \affl_{\affl_U}&.\\
}$$  
where  
  $$
  z = 0_{V\times\affl} \colon
     V\times\affl \hookrightarrow \affl_{V\times\affl} ,\quad
  z_0 = 0_{V} \colon
    V \hookrightarrow \affl_{V} ,\quad
  \pri z = 0_{\affl_U} \colon 
    \affl_U \hookrightarrow \affl_{\affl_U}
  .$$
  are embeddings by zero section. 
Let's also denote embeddings of subschmes $S$, $S_0$ and $\pri S$  by
  $$
  e\colon S\hookrightarrow \affl_{V\times \affl} ,\quad
  e_0\colon S_0\hookrightarrow \affl_{V} ,\quad
  \pri e\colon \pri S\hookrightarrow \affl_{\affl_U}
  .$$
Then
we construct symmetric isomorphisms
    \begin{gather*}\label{fs_aq} 
    \qp\colon B^\prime \simeq Hom_{A^\prime}(B^\prime,A^\prime),\quad
    q_0\colon B_0 \simeq Hom_{A_0}(B_0,A_0), \quad
    q\colon B \simeq Hom_{A}(B,A):\\
    j_0^*(q^\prime) = q_0 = i^*(q), 
    \end{gather*}
    where
      $A^\prime$ denotes $k[\affl_{\affl_U}]$, 
      $B'_{A'}$ denotes the algebra correspondent to $F'$
      and similarly for 
      $A_0$, ${B_0}_{A_0}$, $A$ and $B_A$
  using isomorphism of
      line bundle of relative canonical class 
      and dual module to function ring
    for finite flat morphisms
  and 
  agreed trivialisation of canonical classes.
  
More detailed,
  by the proposition 2.1 of \cite{OP}
  there are agreed isomorphisms
    \begin{multline*}
    \dpri\colon\omega(\Fp)    \simeq 
      Hom(k[\affl_{\affl_U}],k[\affl_{\affl_U}]), \quad
    d_0\colon\omega(F_0)  \simeq 
      Hom(k[\affl_{V}],k[\affl_{V}]),\\
    d\colon\omega(F)      \simeq 
      Hom(k[\affl_{V\times\affl}],k[\affl_{V\times\affl}])
    .\end{multline*}
  Thus 
    to find $q$, $q_0$ and $\qp$ from \ref{fs_aq} it is enough  
    to define agreed trivialisations of
      $\omega(\pri F)$, $\omega(F_0)$ and $\omega(F)$.      
  
  Let's note that 
    for any smooth schemes $X$, $Y$ and
    for any relative morphism $f\colon Y_T\to X_T$ 
      over smooth scheme $T$
    $$     \omega(f) 
    \simeq \omega(Y_T) \cdot f^*(\omega(X_T))^{-1} 
    \simeq c^*(\omega(Y)) \cdot f^*(c^*(\omega(X)))^{-1}
    .
    $$ 
    where $c$ denoted 
    the projections $X_T\to X$ and $Y_T\to Y$ along $T$.
  So trivialisation of canonical classes of $Y$ and $X$
    defines canonical trivialisation of $\omega({f})$
  for all $T$ and $f\colon Y_T\to X_T$.
  In our case  
    $X=\affl$, $Y = \affl$,
    $T_1 = \affl_U$, $T_2 = \affl_V$ and $T_3 = V\times\affl$  
  and fixing trivialisation of canonical class of $\affl$
  we get agreed trivialisations of 
    $\omega(F')$, $\omega(F_0)$ and $\omega(F)$.

get $q_P$ and $q_H$  
  we apply to $\qp$ and $q$ 
  base changes along
    $\pri z$, and $z$ 
  and 
  restriction of scalars along
    $\pri e$, and $e$:
  $$    
  q_P = 
    {\pri e}_*({\pri z}^*(\qp)),\;
  q_H = 
    {e}_*({z}^*(q))
  .$$  

\steparhead{Checking of point 4,i.e. equalities \eqref{WHE_s_relU}}

The first equality of \eqref{WHE_s_relU} holds
due to the existence of $q_0$
because
  $$i^*(q_P) = {\pri e}_*({\pri z}^*(\pri q)) =
    {e_0}_*({z_0}^*(q_0)) 
      {e}_*({z}^*(q)) = j_0^*(q_H).$$
      
To complete the construction it is enough
  to make the second equality of \eqref{WHE_s_relU} true
  (Because it isn't necessary true yet). 
Since
  $S = Z(s)$ is support of quadratic space $(H,q_H)$
  and
  $s\big|_{\prl\times V\times\affl} = s_1$,
$S_1 = Z(s_1)$ is support of 
  it's base change $(H_1,q_{H_1})$
  over unit section 
    $V\times 1\subset V\times\affl$
Moreover
  $H_1 = k[S_1]_{k[\affl\times V}$ 
  (because $H = k[S]_{k[\affl\times V\times\affl]}$). 
Next since 
  $s_1 = \delta\cdot g$ and 
  $g\big|_\Delta$ is invertible, 
$$S_1 = \Delta \coprod Z(g)$$. 
So 
  quadratic space $(H,q_H)$  
  splits over the unit section 
    $V\times 1\subset V\times\affl$
  into sum 
  $$
  (H_1,q_{H_1} = q_H\otimes_{k[V\times\affl} k[V\times 1] = 
    (E,q_E)\oplus (G,q_G),\;
    E\simeq k[\Delta]_{\affl\times V\times \affl},\,
    G\simeq k[Z(g)]. $$    
Since the first summand $E$ 
  is free module of rank 1 over $k[U]$,
the quadratic form $q_E$
  is defined by some 
  invertible function $\lambda$ on $V$. 
Then let's 
  multiply
  $q_P$ and $q_H$ by the inverse function $\lambda^{-1}$ 
    (using 
     the inverse images along the 
     projections of 
       $V\times_U \affl_U$ and 
       $V\times_U V\times \affl$ 
     on the first multiplicator),
Or equivalently, let's
  compose corresponding morphisms in $WCor$
  at the left side with endomorphism of $V$ defined by $\lambda^{-1}$:
   $$
   \Phi\leadsto (\lambda^{-1})_{\Delta}\circ \Phi,\quad 
   \Theta\leadsto (\lambda^{-1})_{\Delta}\circ \Theta
   .$$
Then     
  the quadratic form $q_E$ 
  becomes unit
and it doesn't violate other conditions 
  (i.e. equalities \eqref{PM_s_relU}
   and 
   first equality from \eqref{WHE_s_relU}
   )%
.      

\section{Etale excision}\label{sec_ndimEtEx}
 \label{toveri_etale}

In this section 
  etale excision in an arbitrary dimension $n$
  is proved.
It's generalisation of
   etale excision on curves proved in \cite{AD_WtrSh}.

\begin{theorem}\label{ndimEtEx} 
Let 
  $\bF$ be homotopy invariant presheave with $Witt$-transfers
and 
  $\pi\colon  X^\prime\hm\to X$ be etale morphism of smooth varieties over the field $K$
    that is field of fractions of some variety over the base field $k$.
Let 
  $Z\subset X$ be closed subscheme of codimension 1, 
    such that  
    $\pi$ induces isomorphism between 
      $Z$ and its preimage  $Z^\prime=\pi^{-1}(Z)$
and let
  $z$ be a closed point of $Z$ and $z^\prime$ be it's preimage.

Then $\pi$ induces the isomorphism  
   $$\pi^*\colon \frac{\bF(U-Z)}{\bF(U)}  
              \stackrel{\sim}{\to}  
                     \frac{\bF(U^\prime-Z^{\prime})}{\bF(U^\prime)},$$
  where    
     $U=Spec(\mathcal O_{X,z})$, 
     $U^\prime=Spec(\mathcal O_{X^\prime,z^\prime})$.
\end{theorem}

\begin{remark}\label{remarkEtResn}
In terms of the remark \ref{rem_PreshPair}  
theorem \ref{ndimEtEx}  
means that 
$$i^*\colon \bF^{pair}(U,U\hm-Z)  \to
         \bF^{pair}(U^\prime,U^\prime\hm-Z^\prime)$$
is an isomorphism.
So it follows from following lemma. 
\end{remark}

\begin{lemma}\label{lemW_etexn} 
Let
  $\pi\colon X\to X^\prime$ be etale morphism of 
    smooth varieties with trivial canonical class, 
  $z$ and $z^\prime$ be closed points of $X$ and $X^\prime$
     such that
      $\pi(z')\hm = z$ and 
      residue fields $k(z)$ and $k(z^\prime)$ are isomorphic.    
  Let
  $Z$ and $Z^\prime$ be closed subschemes of 
    $X$ and $X^\prime$ containing $z$ and $z^\prime$
    such that 
      $Z^\prime=\pi^{-1}(Z)$. 
  And let  
    $U \hm= 
      \varprojlim\limits_{z\in V\subset X} V$, 
    $U^\prime \hm= 
      \varprojlim\limits_{z^\prime\in V^\prime\subset X^\prime} V'$. 
Then
\begin{list}{}{\leftmargin 0.7ex \itemsep   5pt}
\item[a)]  
there exists a morphism 
  $$
  \Phi \in WCor_K((U,U\hm-Z),(X^\prime,X^\prime\hm-Z^\prime))
   :\quad
  [\pi\circ\Phi]\hm=[i]\in\overline{WCor_K}((U,U\hm-Z),(X,X\hm-Z))
  $$
;
\item[b)] 
there exists a morphism 
  $$
  \Psi \in WCor_K((U,U\hm-Z),(X^\prime,X^\prime\hm-Z^\prime))
  :  \quad
  [\Psi\circ\pi]\hm=[i^\prime]\in 
    \overline{WCor_K}((U^\prime,U^\prime\hm-Z^\prime),
                                  (X^\prime,X^\prime\hm-Z^\prime))
  $$
.
\end{list}\end{lemma}

\begin{remark}\label{rem_EtExDiad} 
In terms of category $\overline{WCor_K}$ 
statement of previous lemma means that
there exists
  $$\begin{aligned}
  & \Phi \in WCor((U,U-Z),(\pri X,\pri X- \pri Z)) ,\, 
    & \Psi \in WCor((U,U-Z),(\pri X,\pri X- \pri Z))
  \\
  & \Omega\in WCor_K(U,X\hm-Z) ,\, 
    & \Omega^\prime \in WCor_K(U^\prime,X^\prime\hm-Z^\prime)
  \end{aligned}$$ 
such that 
in the group $\overline{WCor_K}((U,U\hm-Z),(X,X\hm-Z^\prime))$ 
following equalities holds: 
  \begin{gather*}
  [\pi\circ\Phi]\hm=[i]+[\Omega],\\
  [\Psi\circ\pi]\hm=[i^\prime]+[\Omega^\prime].
  \end{gather*}

And in terms of the category $WCor_K$ it means 
the existence of following commutative diagrams
$$\xymatrix{
  &
    &
      (U\hm-Z)\times\affl
        \ar@{^(->}[r]
        \ar[dd]^{\Theta^\prime}
        &
          U\times \affl\ar[dd]_\Theta
            &
              &
  \\
  (U\hm-Z)
    \ar@{^(->}[r]
    \ar[drr]^{\quad\pi\circ\Phi}
    \ar[rru]_>>>>>>>>>>{\quad j_0}
    &
      U\ar[drr]\ar[urr]&
        &
          &
            (U\hm-Z)\ar@{^(->}[r]\ar[dll]\ar[ull]
            &
                U
                \ar[dll]_>>>>>>>>>>>>>>{i+\Omega\quad}
                \ar[ull]^>>>>>>>>>>{j_1\quad}
   \\
      &
        &
          X\hm-Z\ar@{^(->}[r]&
            X&
              & 
  }$$
$$\xymatrix{
   &
     &
       (U^\prime\hm-Z^\prime)
       \times
       \affl\ar@{^(->}[r]
         \ar[dd]^{\Xi^\prime}&
           U\times \affl\ar[dd]_\Xi &
             &
  \\
  (U^\prime\hm-Z^\prime)
    \ar@{^(->}[r]
    \ar[drr]^{\quad\Psi\circ\pi}
    \ar[rru]_{\quad j^\prime_0}
    &
      U^\prime\ar[drr]\ar[urr]&
        &
          &
           (U^\prime\hm-Z^\prime)
           \ar@{^(->}[r]\ar[dll]\ar[ull]
           &
             U
             \ar[dll]_>>>>>>>>>>>>>>>{i^\prime+\Omega^\prime\quad}
             \ar[ull]^>>>>>>>>>>{j^\prime_1\quad}
  \\
    &
      &
        X^\prime\hm-Z\ar@{^(->}[r]&X^\prime&
          &,
}$$
such that 
  $j_0$, $j^\prime_0$, $j_1$ and $j^\prime_1$ 
  are zero and unit sections respectively.

\end{remark}

\emph{Proof of the lemma \ref{lemW_etexn}.}

\newcommand{\bX}{\mathcal{X}}
\newcommand{\bXp}{\pri{\mathcal{X}}}
\newcommand{\bXpp}{\ppri{\mathcal{X}}}
\newcommand{\ovbX}{ \overline{ \mathcal X } }
\newcommand{\ovbXp}{ \overline{ \pri{\mathcal X} } }
\newcommand{\ovbXpp}{ \overline{ \ppri{\mathcal X} } }
\newcommand{\bZp}{ \pri{\mathcal Z} }
\newcommand{\bZpp}{ \ppri{\mathcal Z} }    
    
\newcommand{\mcX}{\mathcal{X}}
\newcommand{\mcXp}{\pri{\mathcal{X}}}
\newcommand{\mcXpp}{\ppri{\mathcal{X}}}
\newcommand{\mcZ}{\mathcal{Z}}
\newcommand{\mcZp}{\pri{\mathcal{Z}}}
\newcommand{\mcZpp}{\ppri{\mathcal{Z}}}

Any morphism 
  from $U$ to $X$ (or from $\pri U$ to $X$)  
  is defined by
  quadratic space in $Proj(X,U)$ (or in $Proj(X,\pri U)$ 
  that are categories of projectives modules
  along the morphisms $X\times U\to U$ (or $X\times \pri U\to\pri U$)
  that has the same relative dimension as $X$, i.e. $d$. 
In the proof of etale excision for curves
like as in case of Zariski excision on relative affine line
we deals with morphisms of relative dimension 1. 
To make current situation similar and to 
reduce it to the case of morphisms of relative dimension 1
we use Quillen's trick.

\steparhead{Construction of relative curves} 
I.e. following Quillen's trick 
let's 
  fix some finite  surjective morphism 
    $p\colon X\to \aff^d_K$
    ($d = dim\,X$)
and then 
  construct 
  commutative diagram 
  with Cartesian parallelograms
\newcommand{\priv}{ {\pri v} }
\begin{equation}\label{CurC}\xymatrix{
  \bXp
    \ar[d]_{\varpi}  \ar[rrd]^{\pri v}  \ar[dr]|{\pri p_U}
    &&&\\
  \bX
    \ar[rrd]  \ar[d]_{f^X_U}  \ar[r]|{p_U} 
    &
      \aff^{d}_U
        \ar[ld]  \ar[rrd]^v
        & 
          \pri X
            \ar[d]|<<<{\pi}  \ar[dr]^{\pri p}
            &  \\
  U
   \ar[rrd]_{f^U}  \ar[rr]^i 
   && X
       \ar[r]_p  \ar[d]^{f^X}
       &  \aff^d_K
            \ar[ld]^{pr_\aff}
            \\
  && \aff^{d-1}_K &
  },\end{equation}
  by choosing some linear projection $pr_\aff$.
(The property that all parallelograms are Cartesian
 means that 
    $\mathcal{X}$ is fibred product
      of $pr\circ p$ and $pr\circ p\circ i$, and 
    $\mathcal{X}^\prime$ is fibred product 
      of $pr\circ p\circ\pi$ and $pr\circ p\circ\pi$.)
In addition we may choose $pr_\aff$ in such way that
the projection 
  of closed subscheme $Z\subset X$ to $\aff^{d-1}$ 
is finite and unramified at $z$. 
(These assumptions isn't necessary for following construction 
  but allows to simplify it .)

\steparhead{Description of quadratic spaces}
Then since
  $\mcX$ and $\pri\mcX$ have agreed 
  structures of varieties over $X\times U $ and $\pri X\times U$
  to construct required morphisms 
  $\Phi$, $\Theta$ and $\Omega$ 
    from the remark \ref{rem_EtExDiad} 
  it's enough 
    to find following quadratic spaces: 
\begin{list}{}{\leftmargin 0pt \itemsep 3pt}
\item[\em 1)] 
  Quadratic space $P$ in $Proj(pr_U)$, 
    where $pr_U\colon \mathcal{X}^\prime\to U$
          is canonical projection. 
  I.e. 
  $P\in K[\mathcal{X}^\prime]-mod$
  finitely generated projective over $K[U]$ 
  and
  $K[\mathcal{X}^\prime]$-linear isomorphism 
    $q_P\colon P\simeq Hom_{K[U]}(P,K[U])$  
\item[\em 2)] 
  Quadratic space $H$ in $Proj(pr_{\aff\times U})$, 
    where $pr_{U\times\aff}\colon \mathcal{X}\times \aff\to U\times\aff$ 
          is canonical projection.
  I.e. 
  $H\in K[X\times \aff\times U]-mod$ 
  finitely generated projective over $K[U\times\aff]$ 
  and 
  $K[\mathcal{X}]$-linear isomorphism 
  $q_H\colon H\simeq Hom_{K[U\times\aff]}(H,K[U\times\aff]),$
\item[\em 3)] 
canonical homomorphisms:
\begin{equation}\label{EtndimEx_PM_i}\begin{aligned}
  P\otimes_{K[U]}K[U\hm-Z] \to 
    K[X^\prime\hm-Z^\prime]\otimes_{K[X^\prime]}
      P  \otimes_{K[U]}K[U\hm-Z]
    ,\\
    H\otimes_{K[U]}K[U\hm-Z] \to 
      K[X\hm-Z]\otimes_{K[X]}
        H  \otimes_{K[U]}K[U\hm-Z]
\end{aligned}\end{equation}
are isomorphisms; 
\item[\em 4)] 
  Following isomorphisms of quadratic spaces holds:
  \begin{equation}\begin{aligned}\label{EtndimEx_WHE_ei}
  {j_0}^*(H,q_H)\simeq {\varpi}_*(P,q_P),\\
  {j_1}^*(H,q_H)\simeq (K[\Delta],1)\oplus (G,q_G),
  \end{aligned}\end{equation}
  where 
  $(K[\Delta],1)$ denotes space with unit form on $K[\Delta]$ 
   i.e. the form gotten from unit 
     by isomorphism $K[\Delta]\simeq K[U]$, 
  and
  $G$ is $K[\mathcal{X}]$-module such that 
    $G\simeq K[\mcX\hm-\mcZ]\otimes_{K[\mcX]}G$).
\end{list}

\newcommand{\ovmcX}{\overline{\mcX}}
\newcommand{\ovmcXp}{\overline{\mcXp}}
\newcommand{\ovmcXpp}{\overline{\mcXpp}}
\newcommand{\ovpU}{\overline{p_U}}
\newcommand{\ovppU}{\overline{p^\prime_U}}
\newcommand{\ovvarpi}{\overline{\varpi}}
\newcommand{\ovmcZ}{\overline{\mcZ}}
\newcommand{\ovmcZp}{\overline{\mcZp}}
\newcommand{\ovmcZpp}{\overline{\mcZpp}}

Like as in the proof of
  etale ecxision for curves {see theorem 1 from \cite{AD_ShNisWtr}} and 
  Zariski excision isomorphism 
  on relative affine line 
  made in the previous section 
the construction of 
  quadratic spaces 
  uses some global sections of line bundles \eqref{sectblDiscr}
  on relative projective curves
  that are some compactifications of $\bX$ and $\bXp$.

\steparhead{Compactification}

Now we 
  construct these projective curves 
  and
  involve some additional used later 
    definitions and denotations.

Firstly let's 
  fix an embedding 
  $e_{\affl}\colon \affl\hookrightarrow \prl$,
and 
  consider normalisations of 
    compositions 
      of $p_U$ and $\varphi$ 
      with ${e_{\aff}}_U$
      that is base change of $e_{\aff}$. 

\newcommand{\ovmpU}{\overline{p_U}} 
\newcommand{\ovarpi}{\overline\varpi}  
I.e. let 
  $\ovmpU\colon\ovbX\to \prl_U$ 
    be normalisation of the composition 
    ${e_{\affl}}_U\circ p_U\colon X\to\prl_U$
  and
  $\varpi \colon \ovbXp\to \ovbX$
    is normalisation of composition of 
      $\varpi$ with 
      embedding of $\bX$ into $\ovbX$ 
\begin{equation}\label{compCurU}\xymatrix{
   &
      \ovbXp \ar[r]^{\ovarpi}& 
         \ovbX \ar[r]^{\ovmpU}& 
            \prl_U \\
   &
      \bXp \ar[r]^{\varpi} \ar@{^(->}[u] \ar[dl]_{\pri v} \ar[drr]& 
        \bX \ar[r]^{p_U}   \ar@{^(->}[u] \ar[ld]_{v}      \ar[dr]&  
          \affl_U          \ar@{^(->}[u]_{{e_{\affl}}_U}  \ar[d]\\
   \pri X \ar[r]^{\pi}&
     X&
       &
         U
}. \end{equation}
Let's note that schemes 
  $\ovbX$ are $\ovbXp$ 
  contain smooth open subschemes $\bX$ and $\bXp$
  but
  are not necessarily smooth themselves.
Since  
  $p_U$ and $\varpi$ are finite,
complements to these open subschemes are
  $$D = \ovbX\setminus \bX = \ovpU^{-1}(\infty_U)$$
  (where $\ovppU=\ovpU\circ\ovvarpi$ 
    is normalisation of $p^\prime_U\circ\varpi$
   and
    $\infty_U = \prl_U\setminus \affl_U$ ). 
  
Let's define relative analoges of 
  closed subsets $Z$ and $\pri Z$
  and points $z$ and $\pri z = \pi^{-1}(z)$.
\newcommand{\DelZ}{{\Delta_Z}}
\newcommand{\Delz}{{\Delta_z}}
\newcommand{\DelpZ}{{\Delta^\prime_Z}}
\newcommand{\Delpz}{{\Delta^\prime_z}}
\newcommand{\Deltap}{{\pri \Delta}}
\newcommand{\DelZp}{{\Delta^\prime_Z}}
\newcommand{\Delzp}{{\Delta^\prime_z}}
$$\begin{aligned}
&\mcZ&=\,&v^{-1}(Z) &\subset &\bX ,\;& 
  &\bZp&=\,&\priv^{-1}(\pri Z) &\subset &\bXp 
&\\
&\Delta&=\,&
  \Gamma_{\aff^{d-1}_K}( U \hookrightarrow X ) &\subset &\bX &
&&&&&
&\\ 
&\Delta_Z&=\,&
  \Gamma_{\aff^{d-1}_K}( Z \hookrightarrow X ) &\subset &\bX ,\;& 
&\pri\Delta_Z&=\,&
  \Gamma_{\aff^{d-1}_K}(\pri Z\hookrightarrow\pri X)&\subset&\bXp
&\\
&\Delta_z&=\,&
  \Gamma_{\aff^{d-1}_K}( z \hookrightarrow X ) &\subset &\bX ,\;& 
&\pri\Delta_z&=\,&
  \Gamma_{\aff^{d-1}_K}(\pri z\hookrightarrow\pri X) &\subset&\bXp 
&.\end{aligned}$$
By definition these are closed subschemes in $\bX$ and $\bXp$
but in fact they are closed subschemes in
  projective relative curves $\ovbX$ and $\ovbXp$
because they are finite schemes over $U$.
(Here for schemes $\mcZ$ and $\bZp$ we use 
  that $Z$ is finite over $\aff^{d-1}_K$.)

\steparhead{Description of sections of line bundles} 
In the construction of
 required quadratic spaces $(P,q_P)$ and $(H,q_H)$
following section of line bundles on 
  $\ovbX$ and $\ovbXp$ and $\ovbX\times\affl$
are used.  

\begin{equation}\begin{aligned}\label{sectblDiscr}
&\pri s&\in& \bL(\ovbXp,n \ppri D)\colon\;&
    &div\, \pri s .\bZp = \DelZp,\; 
    &div\,\pri s .{\ppri D } = 0 &,\\
&s & \in&
  \bL(\ovbX \times \affl,lnD\times\affl)\colon\; &  
    &div\,s.{\mathcal{Z}\times\affl} = \Delta_Z\times\affl,\; 
    &div\,s.{D\times\affl}=0 &,\\  
&s_0,s_1 &\in & 
    \bL(\ovbX,lnD)\colon&  
     &s\big|_{\overline{\mathcal{X}^\prime}\times 0} = s_0, \, 
     &s\big|_{\overline{\mathcal{X}^\prime}\times 1}=s_1  &,\\
&& && 
     &\varpi\colon Z(\pri s) \simeq Z(s_0)
        \simeq
          div\, s_0, \;  
     &s_1\big|_\Delta=0
&.
\end{aligned}\end{equation}

Now we assume the existence of such section
and construct the required spaces $(P,q_P)$ and $(H,q_H)$.
The construction of $s$, $s_0$ and $\overline s$ 
  will be given later.
\steparhead{Construction of quadratic spaces} 
Let 
  $$\begin{aligned}
  \pri S = Z(\pri s) &,\quad 
  S_0 = Z(s_0) &,\quad 
  S = Z(s) 
  &,\\
  P = K[\pri S] &,\quad & H = K[S]
  &.
  \end{aligned}$$
The quadratic forms  $q_H$ and $q_P$ 
  are defined by 
  trivialisation of canonical class of $\bX\times\affl$
  and  
  the function 
  $$
  F=(\frac{s}{d^n},pr_{U\times\affl}) \colon 
    \mathcal{X}^\prime \to 
    \affl\times U\times\affl
  $$
  in following sense.

Since 
  $div\,d = D = \ovbX\setminus\bX$,
morphism $F$ is finite.   
By the same reason as in sublemma 3.1.3 from \cite{AD_ShNisWtr} 
since 
  $F$ is finite morphism of varieties of the same dimension,
it is surjective.
And since 
  $F$ is 
  finite morphism 
  of smooth schemes of the same dimension 
it is flat 
(see for example 
  corollary V.3.9. and theorem II.4.7 in \cite{AltKl}).  
Then by the proposition 2.1 of \cite{OP} 
  there is an isomorphism: 
  $$q_B\colon B\simeq Hom_A(B,A),$$ 
  where 
    $B_A$ denotes algebra 
    corresponding to the morphism $F$
    (i.e. 
      $B\simeq K[\bX\times\affl]$, 
      $A \simeq K[\affl\times U\times\affl]$
      and homomorphism $A\to B$ is defined by $F$
     ).
Then let's put    
  $$
  q_P = {e_0}_*({j_0}^*(zf^*(q_B))),\quad
  q_H = e_*(zf^*(q_B))
  .$$
$$\xymatrix{
S^\prime \ar@{=}[r]\ar@{^(->}[d]_{e} 
 & S_0 \ar@{^(->}[r]\ar@{^(->}[d]_{e_0}
  & S \ar@{^(->}[r]^{zf}\ar@{^(->}[d]
   & \bX\times\affl \ar@{^(->}[d]\\
\mathcal{X}^\prime \ar[r]^{\varpi}\ar[dr]
 & \mathcal{X} \ar@{^(->}[r]^{id\times 0} \ar[d]
  & \mathcal{X}\times\affl \ar@{^(->}[r]^{0\times id} \ar[d]
   & \affl\times\mathcal{X}\times\affl \ar[d] \\
& U \ar@{^(->}[r]_{id\times 0}
 & U\times\affl \ar@{^(->}[r]_{0\times id}
  & \affl\times U\times\affl
,}$$
where
  $e_0$ and $e$ denotes embeddings of closed subschemes, 
  and
  $zf$ denotes embedding of zero fibre
    over zero section of affine line at the left side.

\steparhead{
  Checking of condition of point 3),
  i.e. equation \eqref{EtndimEx_PM_i} 
}

The condition 
  \eqref{EtndimEx_PM_i} holds since
  $div\,s^\prime\cap \mathcal{Z}^\prime =
    \Delta^\prime_Z\subset \mathcal{X}^\prime_Z$ and 
  $div\,s\cap \mathcal{Z}\times\affl = 
    \Delta_Z\times\affl\subset \mathcal{X}_Z$.

\steparhead{
  Checking of conditions of point 4),
  i.e. equations \eqref{EtndimEx_WHE_ei}
}
This quadratic forms 
  are agreed in the seems of 
  first equality of \eqref{EtndimEx_WHE_ei}  
due to 
  that fact that both of them 
  are gotten by base changes and scalar restrictions from 
  the same isomorphism $q_B$
and  
  functionality of scalar restriction 
  in respect to 
  base chances.
$${j_0}^*(q_H) = {e_0}_*( {j_0}^*(q_B) ) = 
  {\varpi}_*({\pri e}_*( {j_0}^*(q_B) )) = {\varpi}_*( q_P )
.$$  
where 
  $e_1\colon S_1\hookrightarrow \ovbX$ 
  denote embedding of closed subscheme $S_1 = Z(s_1)$.

It can happens that 
  the second equality from \eqref{EtndimEx_WHE_ei} doesn't holds.
Nevertheless
  the restriction 
    of the space $(H,q_H)$ 
    on the unit section 
    $j_1\colon\bX\times 1\hookrightarrow \bX\times\affl$ 
  splits into direct sum 
    $$\begin{aligned}
    (H_1,q_{H_1}) = {j_1}^*(H,q_H) = 
      (E,q_E) \oplus (G,q_G) ,\\
    G\simeq K[X - Z]\otimes_{K[X]}G   
    ,\end{aligned}$$  
  since 
    it's support $S_1$ splits  
      into disjoint union 
      $$S_1 = Z(s_1) = \Delta\coprod \mathbf R$$
      where 
      $\mathbf R$ id closed subscheme of $\bX - \mcZ$.

Let's give formal proof of the last statement. 
Since
  $s_1|_{\Delta}=0$,
$div\,s_1 > \Delta$.   
Let
  $div\,s_1 = \Delta + R$
    for some effective divisor $R$ on $\ovbX$,
and let
  $s_1 = \delta \cdot r$ 
  where $r\in L(\ovbX,R)$.  
Then since  
  $div\,s_1.\mcZ = \DelZ 
    = \Delta.\mcZ$, 
$R.\mcZ=0$
and 
$r$ is invertible on the subscheme $\mcZ$.
Since 
  $\Delz$ is unique closed point 
    of $\Delta$ 
  and 
  $\Delz\subset \mcZ$, 
$r$ isn't equal to zero at $\Delz$
and is invertible on $\Delta$.
Therefore
  $r$ is invertible on the subscheme $\Delta$. 
Hence
  on the open subscheme $\bX - \mcZ$ containing $\Delta$
  the sheaf of ideals $\bI(s_1)$ is equal to $\bI(\delta)$,
and so  
  $\Delta$ is connected component of $S_1$.
The rest part $\mathbf R$ of $S_1$ 
  is $Z(r)$ 
and since 
  $r$ is invertible on $\mcZ$,
$\mathbf R$ is a closed subscheme of $\bX - \mcZ$.

Thus there is isomorphism of algebras
  $K[S_1] = K[\Delta]\times K[\mathbf R]$
and it induce 
  decomposition of quadratic space ${j_1}^*(q_B)$
  that leads to decomposition of $(H_1,q_{H_1})$
    into sum of spaces $(E,q_E)$ and $(G,q_G)$
      with the supports $\Delta$ and $\mathbf R$,
      and modules $E$ and $G$ 
        are isomorphic to $K[\Delta]$ and $K[R]$ 
      respectively.       
Then since $E$ is free module of rank 1 over $k[U]$,
quadratic form $q_E$ is defined by some invertible function
$\lambda\in K[U]^*$.
Let's multiply 
  quadratic forms $q_P$ and $q_H$ 
  on the inverse function $\lambda^{-1}$ 
Or equivalently let's 
  compose $Witt$-correspondence defined by $(Pmq_P)$ and $(H,q_H)$ 
  with endomorphism of $U$ in $WCor$ defined by function $\lambda^{-1}$
  at the right side.
Then  
  equation \eqref{EtndimEx_PM_i} and 
  the first equation of \eqref{EtndimEx_WHE_ei} 
  remain true
and
  quadratic form on $E$ becomes unit.
So we get the required equality
  $$(
  H_1,q_{H_1})\simeq (1)_{\Delta}\oplus (G,q_G)
  .$$    

\steparhead{
  Construction of sections of line bundles from \eqref{sectlbDiscr}}
Now we preceed 
to construct required 
  sections $s^\prime$, $s_0$, $s_1$ and $s$
  from \eqref{sectlbDiscr}.
Like as in proofs of 
  etale excision isomorphism for curves
  and
  excision isomorphism in presious section 
in this construction 
  lemma \ref{surcon},
  direct images of divisors and
  affine linear homotopy. 
let's note that since 
  $\ovpU$ and $\varpi$ are finite, 
and since 
the sheave $\bL(\infty_U)$ is very ample, 
sheaves 
  $\bL(D)={\ovpU}^*(\bL(\infty_U))$
  $\bL(\ppri D) = \varpi^*(\bL(D))$
are ample too.

We start with construction of some section $\pri s$ 
$$
\pri s\in \bL(\ovbXp,n \ppri D)\colon\;
    div\, \pri s .\bZp = \DelZp,\; 
    div\,\pri s .{\ppri D } = 0,\;
    \ovvarpi\colon Z(\pri s) \simeq \varpi_*(Z(\pri s)) 
$$
for sufficiently large $n$. 
And firstly
  using sublemma \ref{leqovS}  
  we  
  construct it on a closed fibre of $\ovbXp$, 
  i.e. 
    on closed subscheme 
    $\ovbXp_z = \ovbXp\times_U z \subset \ovbXp$. 
Namely
  by sublemma \ref{leqovS}
  applying to morphism $\ovvarpi_z\colon \ovbXp_z \to \ovbX_z$
    (that is closed fibre of $\ovvarpi$)
  the sheaf $\bL(\ppri D_z)$
    (where $\ppri D_z = \ppri D\times_U z$)    
  closed point $\Delzp\in \ovbXp_z$ 
  and
  closed subscheme $ (\bZp_z\cup \pri D_z) - \Delzp $
    (it is exactly the set of closed points of 
      $\bZp$ and $\pri D$ distinct to $\Delzp$) 
for all $n$ larger some $\overline k$
there is a section
\newcommand{\ovsp}{\overline{s^\prime}}
$$\begin{aligned}
\ovsp\in \bL(\ovbXp_z,n \ppri D_z)\colon\qquad
   &div\,\ovsp .\bZp_z = \Delzp,\; 
    div\,\ovsp .{\ppri D_z } = 0  
   ,\qquad
    \ovvarpi_z\colon div\,\ovsp \simeq {\ovvarpi_z}_*(div\,\ovsp)\\
  \text{(or equivalently }\;
   &\ovsp\big|_{\Delzp} = 0 ,\; 
    \ovsp\big|_{\xi} \neq 0\;
      \forall\, \xi\, \in\, (\bZp_z \cup \ppri D_z) - \Delzp\;  
    \text{)} 
&.\end{aligned}$$

\begin{sublemma}[Sublemma 3.1.2 from \cite{AD_ShNisWtr}] 
\label{leqovS}
Let 
  $\pi\colon \pri X\to X$
  be finite morphisms of projective curves
    over infinite field,
  $z$
  be a closed point of $\pri X$,
  $Y$ 
  be closed subscheme of $\pri X$, $Y\not\ni z$
  and 
  $\bL$ ample invertible sheaf on $\pri X$.
Then for all $n$ larger some $k$
  there exists global section $s$ of $\bL^n$ on $\pri X$,
  such that
    $s$ is equal to zero at $z$,
    $s$ isn't equal to zero at any point of $Y$   
  and such that
    restriction of $\pi$ onto $Z(s)$ is closed embedding.  
\end{sublemma}

Then
  to lift the section $\ovsp$ to global section on $\ovbXp$
let's note that
  schemes $\bZp$ and $\ppri D$   
  are finite over local scheme $U$
  and hence they are semi-local.
So since 
  any line bundle on semi-local scheme is trivial 
  (and any divisor is prime)
  there is a section
  $$
  \delta_{\bZp\cup \ppri D}\in \bL(\bZp\cup \ppri D)\colon\quad
    div\,\delta_{\bZp\cup \ppri D} = \DelZp
  .$$
Moreover we may assume that
  $$\delta_{\bZp\cup \ppri D}\big|_{\bZp_z\cup\pri D_z} = 
      \ovsp\big|_{\bZp_z\cup\pri D_z},
  $$     
because 
  $\DelZp\cap \ovbXp_z = \Delzp$ 
  $\ovsp\neq 0$ at another points of $\bZp_z\cup\pri D_z$ 
(Here we use that 
  projection $Z \to \aff^{d-1}_K$ is unramified at $z$
 because we doesn't state in sublemma \ref{leqovS} that
  degree of zero of section $\ovsp$ at $\Delzp$ is one ).   
And now 
  by lemma \ref{surcon}
  applying to 
    the sheave $\bL(\ppri D)$ and
    closed subscheme $\ovbXp_z \cup \bZp \cup \ppri D$  
for all $n$ larger some $k$
there is a section
  $$
  \pri s \colon
    \pri s \big|_{\ovbXp_z} = \ovsp,\;
    \pri s \big|_{\bZp\cup \pri D} = \delta_{\bZp\cup \pri D}
  .$$  

Let's check the properties \eqref{sectpropEtExi_sp}.
$div\,\pri s.\bZp = \DelZp$ 
by definition of $\delta_{\bZp\cup \pri D}$. 
And
$div\,\pri s.\pri D = 0$
because $\DelZp$ doesn't intersect with $\pri D$ 
and so $\delta_{\bZp\cup \pri D}$ is invertible on $D$.

The last property states that 
  the restriction of $\ovvarpi$ onto $Z(\pri s)$ 
  is closed embedding of schemes.
And since it is true at the closed fibre (by definition of $\ovsp$),
it is true over local base $U$.   
In fact 
this property is equivalent to that
morphism of coherent sheaves
  $\epsilon_{\varphi}\colon \bO(\ovbX) \to \ovvarpi_*(\bI(\pri s))$
  induced by $\ovvarpi$ 
is surjective.
But support of its cokernel $Supp\,coker(\epsilon_{\varphi})$
is closed subscheme in relative projective scheme over local scheme.
So if  it isn't empty then its closed fibre isn't empty too. 

Next 
since $\ovvarpi(n\ppri D) = ln D$ (where $l=deg\,ovvarpi$)
there is a section
  $$
  s_0\in L(\ovbX,ln D)\colon 
    div\,s_0 = \ovvarpi(div\,\pri s)
  .$$
Then
  $$
  div\,s_0.\mcZ = \ovvarpi(div\,\pri s . \bZp) = \DelZ,\quad
  div\,s_0.D = \ovvarpi(div\,\pri s . \ppri D) = 0
  .$$

Next
since
  $\Delta\cap\mcZ = \DelZ$ and
  $\Delta\cap D = \emptyset$%
,  
  by lemma \ref{surcon}
  applied to 
  sheave $\bL(ln D)$ on $\ovbX$ and
  closed subscheme $\Delta \cup \mcZ \cup D$
for all $n$ larger some $k_1$
there is a section 
  $$
  s_1\in L(\ovbX,ln D)\colon 
    div\,s_1.\mcZ = div\,s_1.\mcZ,\quad
    div\,s_1.D = div\,s_1.D,\quad
    s_1\big|_\Delta = 0    
  .$$  

Thus for all $n > max(\overline k,k,k_1)$
there are sections $\pri s$, $s_0$ and $s_1$ described above
and finally let's put 
$$
  s = s_0\cdot (1-t) + s_1\cdot t \in L(\ovbX,lnD\times\affl)
.$$ 

\vskip 10pt
\steparhead{Point b)}
\emph{b)}
\newcommand{\hbase}{{\pri U\times\affl}}
\newcommand{\mbase}{{\pri U}}
\newcommand{\pbase}{U}
Now to construct morphisms between $\pri U$ and $\pri X$
we consider 
  fibred products
  $$\bXpp = \bXp\times_U \pri U \; \ovbXpp = \ovbXp\times_U \pri U$$ 
  and their 
  closed subschemes
\newcommand{\Delp}{{\pri\Delta}}
\newcommand{\Delzpp}{{\ppri\Delta_z}}
\newcommand{\DelZpp}{{\ppri{\Delta_Z}}}  
  \begin{gather*}
  \bZpp = {\pri p_U}^{-1}(\pri Z) = \vartheta^{-1}(\bZp) ,\;
  \Delp = \Gamma_{\aff^{d-1}_K}( \pri U \hookrightarrow \pri X ),\;
  \DelZpp 
   = \Gamma_{\aff^{d-1}_K}( \pri Z \hookrightarrow \pri X )
    = \vartheta^{-1}(\DelZp),\;\\
  \Delzpp 
   = \Gamma_{\aff^{d-1}_K}( \pri z \hookrightarrow \pri X )
    =  \vartheta^{-1}(\Delzp)
  \end{gather*}    
\newcommand{\ovvartheta}{{\overline{\vartheta}}}
\newcommand{\pUp}{{p^\prime_U}} 
\newcommand{\ovpUp}{{\overline{p^\prime_U}}}

\begin{equation}\label{compCurprime}
\xymatrix{
\DelZp\ar@{^(->}[r]\ar[dr]
 & \bZpp\ar@{(->}[r]\ar[d]
  & \bXpp \ar@{(->}[r]\ar[d]^{\vartheta}
   & \ovbXpp \ar[r]\ar[d]^{\ovvartheta}
    & \pri U\ar[d]^{\pi}\\
& \bZp\ar@{(->}[r]
 & \bXp \ar@{(->}[r]\ar[d]^{\pUp}
  & \ovbXp\ar[r]\ar[d]^{\ovpUp}
   & U\ar[d]^{f^U}\\
&& \aff^n\ar@{^(->}[r]
  & \prl\times\aff^{n-1}\ar[r]^{pr_{\pro}}
   & \aff^{n-1}}.
\end{equation}

\steparhead{Description of quadratic spaces}
To find required morphisms
  $\Psi$, $\pri\Omega$, $\Xi$  
  from remark \ref{rem_EtExDiad} 
it is enough   
to construct following data
 \begin{list}{}{\leftmargin 0pt \itemsep 3pt}

\item[\em 1)] 
  The quadratic space $(P,q_P)$ in the category $Proj(pr_U)$. 
  I.e. 
    $P\in K[\mathcal X^\prime]-mod$ 
      that is finitely generated over $K[U]$ 
    and 
    $K[\mathcal{X}^\prime]$-linear 
      isomorphism $P\simeq Hom_{K[U]}(P,K[U]),$
\item[\em 2)] 
  The quadratic space $(H,q_H)$ 
    in the category $Proj(pr_{U^\prime\times\affl})$. 
  I.e. 
    $H\in K[\mathcal{X}^{\prime\prime}\times \affl]-mod$ 
      that is finitely generated over $K[U^\prime\times\affl]$ 
    and 
    $K[\mathcal{X}^{\prime\prime}\times\affl]$-linear 
      isomorphism
      $H \simeq Hom_{K[U^\prime\times\affl]}(H,K[U^\prime\times\affl]),$
\end{list} 
such that
\begin{list}{}{\leftmargin 0pt \itemsep 3pt}
\item[\em 3)] 
canonical homomorphisms
\begin{equation}\begin{aligned}\label{EtndimEx_PM_s}
&P&\otimes_{K[U]}&K[U\hm-Z]& 
  &\to& 
    K[X^\prime\hm-Z^\prime] &\otimes_{K[X]}&
      &P&  \otimes_{K[U]}&K[U\hm-Z]&,\\
&H&\otimes_{K[U^\prime]}&K[U^\prime\hm-Z^\prime]&
  &\to& 
     K[X^\prime\hm-Z^\prime] &\otimes_{K[X^\prime]}&
      &H&   \otimes_{K[U^\prime]}&K[U^\prime\hm-Z^\prime]&
\end{aligned}\end{equation} 
are isomorphisms, 
\item[\em 4)]
There exists isomorphisms of quadratic spaces
  \begin{equation}\begin{aligned}\label{EtndimEx_WHE_es}
  {j_0}^*(q_H)\simeq {(id\times\pi)}^*(q_P),\\
  {j_1}^*(q_H)\simeq (1)_{\Delta^\prime}\oplus q_G.
  \end{aligned}\end{equation}
  Here 
    $(1)_{\Delta^\prime}$ denotes 
      unit quadratic form on $K[\Delta^\prime]$ 
    (i.e. the form that is gotten from unit by isomorphism
        $K[\Delta^\prime]\simeq K[U^\prime]$), 
  and 
    $G$ is 
    $K[\mathcal{X}^\prime]$-module equipped with quadratic form $q_G$ 
    such that
      $$G\simeq K[X^\prime\hm-Z^\prime]\otimes_{K[X^\prime]}G.$$
\end{list}

\steparhead{Construction of sections $\pri s$, $s_0$, $s_1$, $s$ }
Then since 
  $\bL(\infty_U)$ is ample and
  $\ovpUp$ is finite, 
$\bL(\pri D) = {\ovpUp}^*(\bL(\infty_U))$ is ample. 
Similar the sheaf $\bL(\pppri D)$ is ample
  because $\vartheta$ is finite and $\pppri D= \vartheta^{-1}(\pri D)$. 

As usual
  the construction is based on  
  lemma \ref{surcon}.
We will consequently 
  find the restrictions of section on some closed subschemes
  and
  continue them to whole schemes by using of lemma \ref{surcon}
Thus to execute the algorithm we should
  previously choose sufficiently large $n$ 
  such that 
    all restriction homomorphisms of 
      groups of section of sheaves 
    were surjective.

Let's give first short description of this process as follows
$$\begin{aligned}\Gamma( \bXpp , \bL(n\pri D))\; &\ni& \pri s \colon \;\;&\; 
  \pri s\big|_{\bZp\coprod\pri D} = \delta_{{\bZp\cup\pri D}} ,\quad
&&\\
\Gamma( \bXpp , \bL(n\ppri D) )\; &\ni &
  s_0 \;=\; &{\pi^{\pri X}}^*(\pri s)
&&\\
\Gamma( \bXpp , \bL(n\ppri D) )\; &\ni& s_1 \colon \;\;&\;
  s_1\big|_{ \bZpp \coprod \pppri D } 
     = s_0\big|_{ \bZpp \coprod \pppri D } ,\quad 
  s_1\big|_{ \Deltap } = 0   
&&\\  
\Gamma( \bXpp \times \affl, \bL(n(\ppri D\times \affl) )\; &\ni &
  s \;=\; &(1-t)\cdot s_0 + t\cdot s_1      
&&.\end{aligned}$$ 
Now let's explain it in detail.
In first row we apply lemma \ref{surcon}  
  to line bundle $\bL(D)$   
  and closed subscheme $\bZp\coprod\pri D \subset \ovbXp$.
Scheme $\bZp\coprod\pri D$ 
is quasi-finite projective and hence projective 
  over local scheme $U$.
Hence it is semi-local.
Let's consider $\DelZp$ as divisor in $\bZp$ 
  (and in $\bZp\coprod\pri D$).
Since any divisor on semi-local scheme is prime
         (and any line bundle is trivial)              
  there is some  
  section 
  $$\delta_{\bZp\cup\pri D}\,\in\, 
    \bL(n\pri D)\big|_{\bZp\coprod\pri D}\colon\quad 
    div\,\delta_{\bZp\cup\pri D}
     = \DelZp
  $$
and we set $\pri s$ to be equal to this section on $\bZp\coprod\pri D$.  
(In particular it implies that $\pri s $ is invertible on $\pri D$.)   

Next we define section $s_0$ explicitly as inverse image of $\pri s$.
Further we  apply lemma \ref{surcon} 
  to the same bundle and 
  to closed subscheme $\bZp\cup\Delp\coprod \pri D$.
To check that conditions are compatible 
it is enough to note that 
  $\Delp \cap \bZp = \ppri\Delta_Z $
and that
  $\DelZpp = \vartheta^{-1}(\DelZp)$
  (and hence $s_0\big|_{\DelZp} = 0$).
Finally we explicitly define $s$ as homotopy o $s_0$ and $s_1$. 
Then it is immediate follows from definition that 

$$\begin{aligned}
&s^\prime &\in                                      \;& 
  \Gamma(\bXp,\bL(nD^\prime)) \colon  \;&
 &div\,\pri s.\bZp = \pri\Delta_{z}           ,\, 
 &div\,\pri s.D^\prime           = 0                             &,\\
&s &\in                                             \;& 
  \Gamma(\bXpp\times\affl,
         \bL(nD^{\prime\prime\prime}))\colon        \;&
 &div\,s.(\bZpp\times\affl)
       = \ppri\Delta_{z}\times\affl                       ,\, 
 &div\,s.(\pppri D\times\affl) = 0                 &,\\
&&&& &s\big|_{\bXpp\times 0} = s_0   ,\, 
 s\big|_{\bXpp\times 1} = s_1 &        &,\\
&s_0 &\in                                           \;& 
  \Gamma(\bXpp,
         \bL(nD^{\prime\prime\prime}))\colon        \;&
 &s_0 = {\pi^{X^\prime}}^{*}(s^\prime)&                      &,\\
&s_1 &\in                                           \;& 
  \Gamma(\bXpp,
         \bL(nD^{\prime\prime\prime}))\colon        \;&
 &s_1|_{\Delta^{\prime\prime}_{z}} = 0&                      &.
\end{aligned}$$

\steparhead{Definition of quadratic spaces}
The construction spaces $(P,q_P)$ and $(H,q_H)$
is similar to that one in previous section. 
In short, 
  divisors of sections $s^\prime$ and $s$  
  will be supports of 
  required modules
and
  to define quadratic forms,
  we variate zero divisors of $s^\prime$ and $s$
    to a families of schemes 
    that constitutes together smooth schemes
    with trivialized canonical class.
 
In detail 
  we start with three following 
  morphisms of relative projective curves
  corresponding to the 
  sections 
    $s^\prime$, $s_0$ and $s$.  
\newcommand{\ovF}{{\overline{F}}}
\newcommand{\ovFp}{{\overline{F^\prime}}}
\newcommand{\movF}{{\overline{F_0}}}
$$\begin{aligned}
 &\ovFp &=\;& 
   ([\pri s:{\pri d}^n],pr_{U})
 &\colon & 
   \ovbXp \to \prl\times U,\\ 
 &\movF &=\;& 
   ([s_0:{\pppri d}^n],pr_{\pri U})
 &\colon & 
   \ovbXpp\to \prl\times \pri U,\\
 &\ovF &=\;& 
   ([s:{\pppri d}^n],pr_{\pri U\times\affl})
 &\colon & 
   \ovbXpp\times\affl \to 
     \prl\times \pri U\times\affl
\end{aligned}$$

The fibres of these morphisms 
  over zero section of  
  affine line (at the left side)
are 
  $S$, $S_0$, and $\pri S$.     
Then 
  to get a morphisms of
  smooth varieties
we consider 
  neighbourhoods of 
  these zero fibres of 
  $\ovF$, $\movF$ and $\ovFp$. 
I.e.
we consider
  closed subscheme of $\prl\times \pbase$
    that is image of $\pri D$ along $\ovFp$
  and its preimages in $\prl\times \mbase$ and $\prl\times \hbase$
$$
\pri C = \ovFp(\pri D) ,\;
C_0 = \movF(\pppri D) = (id_\prl\times\pi)^{-1}(\pri C) ,\;
C = \ovF(\pppri D \times\affl) = (pr_h\times id_{\mbase})^{-1}(C_0)
.$$ 
And then
  we consider base changes of 
  $\ovF$, $\movF$ and $\ovFp$
  over open subschemes that are complements to this closed subschemes
\newcommand{\pbV}{ {\mathcal V^\prime} }
\newcommand{\mbV}{ {\mathcal V_0} }
\newcommand{\hbV}{ {\mathcal V} }
$$\begin{aligned}
&\pri F &=& \frac{\pri s}{ {\pppri d}^n } 
  &\colon& \pbV \to \prl_\pbase - \pri C ,\;&
  \pbV &=& \ovbX\times\affl - \ovFp^{-1}(\pri C) &,\\
&F_0 &=& \frac{s_0}{ {\pppri d}^n } 
  &\colon& \mbV \to \prl_\mbase - C_0 ,\;&
  \mbV &=& \ovbXpp - \movF^{-1}(C_0) = \pbV \times_U \pri U &,\\
&F &=& \frac{s}{ {\pri d}^n } 
  &\colon& \hbV \to \prl_\hbase - C,\;&
  \hbV &=& \ovbXp - \ovF^{-1}(C) = \mbV\times\affl 
&.\end{aligned}$$
Since
  $\pri s \big|_{\pri D}$, $s_0 \big|_{\pppri D}$ and $s\big|_{\pppri D \times\affl}$ 
  are invertible%
,  
  $\pbV \subset \bXp$,
  $\mbV \subset \bXpp$,
  $\hbV \subset \bXpp\times\affl$
and hence
  schemes $\pbV$, $\mbV$, $\hbV$ are smooth.
And  by the same reason
  zero sections of relative projective lines 
  doesn't intersect with 
  $\pri C$, $C_0$ and $C$.

So we get following commutative diagram with Cartesian squares 
$$\xymatrixrowsep{0.1in}\xymatrixcolsep{0.12in}\xymatrix{
&&& \bXpp\times\affl
    \ar[ddddl]
     &&& \bXpp       
         \ar@{^(->}[lll]
         \ar@{^(->}[rrr]
         \ar[ddddl]
          &&& \bXp
              \ar[ddddl] \\
 && S 
    \ar@{^(->}[dl] 
    \ar[ddd]_{}
    \ar@{^(->}[ru]^{cs}
     &&& S_0 
         \ar@{^(->}[dl]_{}
         \ar@{^(->}[lll]
         \ar@{^(->}[rrr]
         \ar[ddd] 
         \ar@{^(->}[ru]^{cs_0} 
          &&& \pri S 
              \ar@{^(->}[dl]
              \ar[ddd]
              \ar@{^(->}[ru]^{\pri cs} &\\
  & \hbV 
    \ar@{^(->}[dl] 
    \ar[ddd]_{F} 
     &&& \mbV 
         \ar@{^(->}[dl] 
         \ar[ddd]_{F_0} 
         \ar@{(->}[lll]
         \ar@{(->}[rrr] 
         &&& \pbV
             \ar@{^(->}[dl] 
             \ar[ddd]_{\pri F} && \\
    \ovbXpp\times\affl
    \ar[ddd]_{\ovF}  
     &&& \ovbXpp
         \ar[ddd]_{\movF} 
         \ar@{(->}[lll]
         \ar@{(->}[rrr] 
          &&& \ovbXp 
              \ar[ddd]_{\ovFp} &&& \\
 && \hbase 
    \ar@{^(->}[dl]_{zf}
     &&& \mbase  
         \ar@{^(->}[dl]_{zf_0} 
         \ar@{^(->}[lll]^{ id_{\pri U}\times 0 }  
         \ar@{^(->}[rrr]_{\pi}
          &&& \pbase
              \ar@{^(->}[dl]_{\pri zf} &\\
  & \prl_\hbase - C 
    \ar@{^(->}[dl]
     &&& \prl_\mbase - C_0
         \ar@{^(->}[dl] 
         \ar@{^(->}[lll]^{ id_{\affl\times \mbase}\times 0 } 
         \ar[rrr]_{ id_{\affl}\times \pi } 
          &&& \prl_\pbase - \pri C 
              \ar@{^(->}[dl] &&\\
  \prl_\hbase  
   &&& \prl_\mbase  
       \ar@{^(->}[lll]^{id_{\prl\times\mbase}\times 0} 
       \ar[rrr]_{ id_\prl\times{\pi} }  
        &&& \prl_\pbase  &&&
}$$ 
where
  $zf$, $zf_0$ and $\pri zf$
  denotes zero sections of the affine line at the left side
  and
  $cs$, $cs_0$ and $\pri cs$
  denotes closed embeddings of subschemes.

Morphisms 
  $\overline{F}$, $\overline{F_0}$, $\overline{F^\prime}$
 are finite as quasi-finite  projective morphisms.
Hence  
$F$, $F_0 $, $F^\prime$ are finite too. 
And since 
  $\bXpp$, 
  $\affl\times \mbase$, 
  $\bXp $ and 
  $\aff\times \pbase$ 
  are
  essential smooth schemes of the same dimension $d+1$, 
morphisms $F$, $F_0 $, $F^\prime$
are flat.  

Then if we denote by 
  $B_A$, ${B_0}_{A_0}$ and $B^\prime_{A^\prime}$, 
the algebras corresponding to the 
  finite flat morphisms of essential smooth schemes 
  $F$, $F_0$, and $F^\prime$, 
then proposition  2.1 from \cite{OP}, 
 provides isomorphisms
\begin{gather*}
 q_B\colon B\simeq Hom_A(B,A),\,
 q_{B_0}\colon B_0\simeq Hom(B_0,A_0),\, 
 q_{B^\prime}\colon B^\prime\simeq Hom(B^\prime,A^\prime \colon\\
   {id_{\affl\times U^\prime}\times 0}^*(q_B)     = q_{B_0} ,\, 
   {id_\affl\times\pi^{X^\prime}}^*(q_{B^\prime}) = q_{B_0} 
.\end{gather*}
Now 
applying
  base changes along zero sections 
    $zf$, $zf_0$ and $\pri zf$
  and
  restrictions of scalars along closed embeddings 
    $cs$, $cs_0$ and $\pri cs$  
we define quadratic forms 
$$\begin{aligned}
&(K[S],q_S) 
  &=\;& zf^*(B,q_B) 
,\;&       
(K[S_0],q_{S_0})
  &=& {zf_0}^*(B_0,q_{B_0}) \;\quad
,\;\;&  
(K[\pri S],q_{\pri S})
  &=\;& {\pri zf}^*(\pri B,q_{\pri B})  
,    
\\
&(H,q_H) 
  &=\;& cs_*(K[S],q_S) 
,\;& 
(H_0,q_{H_0}) 
  &=& {cs_0}_*(K[S_0],q_{S_0}) 
,\;\;& 
(P,q_P) 
  &=\;& {\pri cs}_*(K[\pri S],q_{\pri S}) 
.\end{aligned}$$
In particular it means that 
$$
H = K[S]_{K[\bXpp\times\affl]} ,\quad P = K[\pri S]_{K[\bXp]}
$$
and 
simultaneously with defining of quadratic forms
we have prove that
this modules are projective finitely generated 
  over $K[\hbase]$ and $K[\pbase]$.

\steparhead{Checking of point 3,i.e. equalities \eqref{EtndimEx_PM_s}}
Condition \eqref{EtndimEx_PM_s} 
holds, 
because
$$
div\,s\cap (\bZpp\times\affl)
  = \ppri\Delta_z\times\affl 
  \subset 
    \bXpp_{\pri Z}\times\affl
,\quad
div\,\pri s\cap \bZp
  = \pri\Delta_z
  \subset 
    \bXp_{Z}
$$
where 
  $\bXp_{Z}\to Z$ denoted 
  base change of $\bXp\to U$ 
  along the embedding $Z\hookrightarrow U$
  and  
  $\bXpp_{\pri Z}\to\pri Z $ denoted 
  base change of $\bXpp\to \pri U $ 
  along the embedding $\pri Z\hookrightarrow \pri U$.

\steparhead{Checking of point 4,i.e. equalities \eqref{EtndimEx_WHE_es}}
The first equality of \eqref{EtndimEx_WHE_es} holds
because of 
 existence of form $q_{H_0}$
 and
 functoriality of scalar restrictions 
  in respect to base changes
$$
{{U^\prime}\times 0}^*(q_H) = q_0 = {\pi^{X^\prime}}^*(q_P)
.$$

The second equality from \eqref{EtndimEx_WHE_es} 
doesn't necessary holds in full power for defined above spaces.  
But 
  due to properties of $s_1 = s\big|_{\bXpp\times 1}$
it is true that 
  $(H,q_H)$ splits 
  after the restriction on $\bXpp\times 1$
  into direct sum of some spaces
  \begin{equation}\label{splitH1}\begin{aligned}
  (H_1,q_{H_1}) = {j_1}^*(H,q_H) \simeq 
    (\Lambda,q_\Lambda) \oplus (G,q_G)\colon &\\
  \Lambda \simeq K[\ppri\Delta]_{K[\bXpp]},\;
  \lambda\in K[\ppri\Delta]^* = K[\mbase]^*,\quad  
  K[]\otimes_K[] G \simeq G   
  &.\end{aligned}\end{equation}
Because
  $$H_1 = K[S_1]_{K[\bXpp]},\; S_1 = Z(s_1)$$
  and
  $q_{H_1}$ can be considered $K[S_1]$-linear isomorphism.
But
  the equalities 
    $$s_1|_{\ppri\Delta}=0,\; div\,s_1.\bZpp = \Delta_{Z^\prime}$$
  implies that
  $$S_1 =  Z(s_1) = \Delta \coprod \mathbf R $$  
  where 
    $\mathbf R$ is closed subscheme of $\bXpp-\bZpp$

The formal proof of last fact
can be done for example like as in similar place in point a).

To satisfy the second condition of \eqref{EtndimEx_WHE_es} 
we modify constructed quadratic forms in such way
  to make 
  morphism of pair $(U^\prime,Z^\prime)\to (X^\prime,Z^\prime)$ 
    defined in $WCor$ 
    by the first summand $(\Lambda,\lambda)$ from \eqref{splitH1} 
  to be equal to
    the morphism defined by space $(\Lambda,1)$. 
Namely 
  we  
    choose a function
  $$\pri\lambda \in K[\pbase]^* \colon 
   \pri\lambda(z)=\lambda(z^\prime)$$
and  
  multiply quadratic forms $q_H$, $q_P$ 
  and other forms in the construction
on its inverse\label{inverttinverse} $\pri\lambda^{-1}$
using that fact that all bases of 
  considered spaces
  (i.e. schemes $\hbase$, $\mbase$  and $pbase$) 
  are schemes over $\pbase$. 
After such multiplication  
the function $lambda$ 
  defines quadratic form 
  on the first summand of $(H_1,q_1)$
becomes equal to 1 at the closed point of $U^\prime$.
Then 
  because of 
    proved in point a) 
    existence of right inverse 
like as in  
  the proof of 
  surjectivity of etale excision on curves
  (see sublemma 3.1.4 from  \cite{AD_ShNisWtr}) 
morphisms of pairs 
 $(U^\prime,Z^\prime) \to (X^\prime,Z^\prime)$
  in $WCor$ 
 defined by  quadratic spaces
   $(\Lambda,\lambda)$ 
   and 
   $(\Lambda,1)$ 
 are equal.

{\em Lemma is proved}

\section{Strictly homotopy invariance}
\label{sec_DHI}

It was proved in \cite{AD_WtrSh}
that Nisnevich sheafication of 
homotopy invariant presheave with $Witt$-transfer
is homotopy invariant. 

At this section 
  Zariski excision and etale excision isomorphism
    proved in previous sections
  are applied to the prove of \label{toveri-of}
  strictly homotopy invariance of
    homotopy invariant sheaves with $Witt$-transfers.

\begin{theorem}\label{DHI}
Nisnevich sheafication
  of the homotopy invariant presheave with $Witt$-transfers
is strictly homotopy invariant.
I.e. the presheaves of 
  Nisnevich cohomologies $H^i_{Nis}( F_{Nis}), i\geqslant 0$ 
  of the sheafication 
    of homotopy invariant presheave with $Witt$-transfers $\bF$
  are homotopy invariant.

\end{theorem}

We start with
  particular case of the theorem 
  for spectres of fields
that states that
  cohomology of
  affine line over any point in $Sm_k$
are isomorphic to the cohomology of this point. 

\begin{theorem}
Let 
  $\bF$ be homotopy invariant presheave with $Witt$-transfers   
  and
  $K=k(X)$ be a field of functions of 
  some smooth variety $X$ over the base field $k$.
Then 
  \begin{equation*}\begin{aligned}
  &\bF_{Nis}(\affl_K) \simeq \bF(\affl_K) & &\\
  &H^i_{Nis}(\affl_K,\bF_{Nis}) \simeq 0 ,\quad & i>0 &. 
   \end{aligned}\end{equation*}

\end{theorem}
{\em Proof of theorem.} 
The firs equality for global section $\bF(\affl_K)$
is a particular case of homotopy invariance of $\bF_{Nis}$ proved in \cite{AD_WtrSh} (theorem 4).

To compute higher Nisnevich cohomologies
let's consider following
flat resolvent of the restriction of $\bF_{Nis}$ on small etale site over $\affl_K$
\begin{equation}\label{fresIPre} 
{\bF_{Nis}}\big|_{\affl_K}\rightarrow 
\eta_*(\bF(\eta)\xrightarrow{d}
\sum\limits_{z\in MaxSp(\mathcal U)} 
  z_*(\coker({\bF(U^h_z)\to {\bF(U^h_z-z)}})),
\end{equation}
  where 
  $\eta$ is generic point of $\affl_K$,
  $\eta_*$ denotes direct image along the embedding $\eta\to\affl_K$, 
  $MaxSp(\mathcal U)$ denotes the set of closed points $z$ in any element $\mathcal U$ of small etale site, 
  $z_*$ denotes direct image along the embedding $z\to\mathcal U$
  and $U^h_z$ denotes corresponding Hensel local scheme 
    (that is spectrum of Henselisaction of local ring at $z$).
The first arrow is injective by the injectivity theorem for presheaves with $Witt$-transfers
    proved by K.Chepurkin in his diploma work. 
The second is surjective  
  because $dim \mathcal U= 1$ for all $\mathcal U$ in small etale site of $\affl_K$ 
  and so $U^h_z-z$ is generic point of $U^h_z$.
And exactness in middle term follows from definition of Nisnevich sheafication.

Then since this is a flat resolvent of length 2,
 $H^0_{Nis}(U)=\ker(d(U))$, $H^1_{Nis}(U)=\coker(d(U))$, 
 and higher cohomologies are zero. 

To show that 
  $H^1_{Nis}(U)$ are zero 
it is enough 
to check exactness in last term 
of the sequence combined
  by global sections of $\bF$ and 
  global sections of the resolvent \eqref{fresIPre}
\begin{equation}\label{CG_F}
\bF(U)\xrightarrow{i}
\bF(\eta)\xrightarrow{d^1}
\sum\limits_{z\in MaxSp(U)}\coker({\bF(U^h_z)}\to {\bF(U^h_z-z)})
.\end{equation} 
But the second arrow is surjective due to  
  surjectivity of the excision homomorphism 
  proved in theorem 2 of \cite{AD_WtrSh}. 
{\em Theorem is proved}

{\em Proof of the theorem \ref{DHI}.} 
The proof is based on 
  the main theorems form \cite{AD_WtrSh} and 
  theorems \ref{locrelAfZarEx} and \ref{ndimEtEx}
and
it is similar to the proof of 
  the same result for sheaves with $Cor$-transfers 
  form \cite{VSF_CTMht}. 

We will prove it by induction on $i$ 
  for all homotopy invariant sheaves with $Witt$-transfers
  simultaneously. 
The base of induction is $i=0$ and it obviously holds.
Suppose that the statement of theorem holds for $i-1$.

Let $p\colon X\times\affl\to X$ denotes  projection
and $i\colon X\to X\times\affl$ denotes embedding of zero section.
The composition 
  $i^*\circ p^*\colon H^i_{Nis}(X) \to H^i_{Nis}(X)$
is identity.
So to prove that  
  $$p^*\colon  H^i_{Nis}(X)\to H^i_{Nis}(X\times\affl)$$ is isomorphism
is equivalent to prove that  
  kernel of $i^*$ is zero. 

Let $$a\in H^i_{Nis}(X\times\affl,\bF),\; i^*(a)=0.$$
Since 
  $H^i_{Nis}(\eta\times\affl)=0$ 
    (where $\eta$ is generic point of $X$)
, 
  $a\big|_{\eta\times\affl}=0$ 
and hence 
  $a\big|_{U\times\affl}=0$ 
    for some open affine subscheme $U\subset X$.
Let $Z_1 = sign\, Z$ is the subset of singular points of $Z$.
Let $U_1 = X\setminus Z_1$.
By the following lemma applied to $U_1$ and $U$ 
  $a\big|_{U_1\times\affl}=0$.
And by induction we can prove that 
$$
a\big|_{U_i\times\affl}=0,\;\text{where}\;
U_i = X\setminus Z_i,\;Z_i = sing\, Z_{i-1}.$$
Since $k$ is perfect
  $$dim\, Z> dim\, Z_1 >\dots > dim\, Z_i >\dots .$$
Hence for some finite $i$ 
  $Z_i = \emptyset$ and $U_i = X$.  
Thus we get that
  $a\big|_{X\times\affl}=0$ for any $a\in ker i^*$.
  
\begin{lemma}
Let $U$ be open subscheme of smooth affine scheme $X$
and $Z = X\setminus U$ is smooth.
Let 
  $$a\in H^i_{Nis}(X\times\affl,\bF), 
    i^*(a)=0, a\big|_{U\times\affl}=0$$
  (where $i$ denotes zero section $X\times Z\times\affl$). 
Then $a=0.$
\end{lemma}
{\em Proof of lemma.}
Let's consider short exact sequence of presheaves
  \begin{equation}\label{ZXUseq}
  \bF\xrightarrow{\epsilon} j_*(j^*(\bF))\to coker \epsilon
  \end{equation}
  where $j$ denotes embedding of $U\times \affl$ into $X\times\affl$.
Using theorems \ref{ndimEtEx} and \ref{locrelAfZarEx} we get
  $$coker \epsilon(V)\stackrel{def}{\simeq}
    \frac{\bF(V\setminus Z}{\bF(V)}\simeq
    \frac{\bF(\mathbb G_m\times(Z\cap U))}{\bF(\affl\times(Z\cap U))}
    \stackrel{def}{=}
    \bF[1](Z\cap U)
  .$$
Thus there is isomorphism of presheaves 
$coker\epsilon\simeq i_*(\bF[1])$ 
  where 
    $i$ denotes 
    the embedding of $Z\times\affl$ into $X\times\affl$ and 
    $\bF[1]$ 
    is regarded as presheave on $Z\times\affl$.    
The sheave   
$\mathcal Hom(\mathbb G_m,\bF)$  
splits into direct sum $\bF\oplus \bF[1]$
because of existence of 
  the projection $\mathbb G_m\to pt$ and 
  unit section $pt\to \mathbb G_m$.
Hence $\bF[1]$ is the Nisnevich sheave
and \eqref{ZXUseq} is the exact sequence of Nisnevich sheaves.
Thus the sequence \eqref{ZXUseq} induce the long exact sequence of Nisnevich cohomology
$$\dots H^{i-1}_{Nis}( X\times\affl,i_*(\bF[1]) )  
  \stackrel{\delta}{\to}
        H^i_{Nis}    ( X\times\affl ,\bF ) \to 
        H^i_{Nis}    ( X\times\affl, j_*(j^*(\bF) )  \to dots.$$   
Note that
$H^{i-1}_{Nis}( X\times\affl,i_*(\bF[1]) )\simeq H^{i-1}(Z\times \affl, \bF[1])$and
$$H^i_{Nis}( X\times\affl, j_*(j^*(\bF) ) \simeq H^i_{Nis}(U\times\affl, \bF)).$$
We get following commutative diagram 
$$\xymatrix{
H^{i-1}(Z\times \affl, \bF[1])
\ar@<0.5ex>[d]^{i^*} \ar[r] & 
       H^i_{Nis}(V \times\affl, \bF)) 
       \ar@<0.5ex>^[d]{i^*} \ar[r] &
                     H^i_{Nis}(U \times\affl, \bF))  
                     \ar@<0.5ex>[d]^{i^*}\\
H^{i-1}(Z, \bF[1])
\ar@<-0.5ex>[u]_{p_*}\ar[r]  &
       H^i_{Nis}(V, \bF)) 
       \ar@<-0.5ex>[u]_{p_*}\ar[r] &
                     H^i_{Nis}(U, \bF))
                     \ar@<-0.5ex>[u]_{p_*}
}$$
Since $a\big|_{U\times\affl} = 0$,
then 
  $$
  a\big|_{V\times\affl} = \delta(b),\; 
  b\in H^{i-1}(V\times\affl,\bF[1])
  .$$
Since 
  the sheave $\bF[1]$ 
  is direct summand of $\mathcal Hom(\mathbb G_m,\bF)$,
it is homotopy invariant sheave with $Witt$-transfers, 
and by induction assumption 
  $$H^{i-1}(V\times\affl,\;[1])\simeq H^{i-1}(X,\bF[1]).$$
Hence $b=p^*(i^*(b))$ 
and
since 
$\delta$ commutes with $p^*$ and $i^*$, 
$$a\big|_{V\times\affl}=p^*(i^*(a\big|_{V\times\affl})=0.$$  
{\em Lemma is proved.}

{\em Theorem is proved.}

\section{Preserving of $Witt$-transfers}\label{secWtrNTop}
\label{secWtrNTop}

\newcommand{\DPre}{{D^-(Pre)}}
\newcommand{\DShN}{{D^-(ShNis)}}
\newcommand{\DShNis}{{D^-(ShNis)}}
\newcommand{\DPreWtr}{{D^-(PreWtr)}}
\newcommand{\DShNWtr}{{D^-(\ShNW)}}

In current section we prove that
  Nisnevich sheafication and cohomology
  of presheave with $Witt$-transfers
  has $Witt$-transfers.

The key fact that provides such behaviour 
is that 
  any $Witt$-correspondence 
    form locale Henzel scheme
    to any other scheme 
  is equal to 
  the sum of 
  $Witt$-correspondence between 
    locale Henzel schemes. 

\begin{theorem}
\label{TrForTheSheaf}
There is an 
  unique structure 
      of presheave with $Witt$-transfers
    on the Nisnevich sheafication $\bF_{Nis}$
  such that 
  $\varepsilon: \bF\to\bFN$ 
  is homomorphism of presheaves 
    with $Witt$-transfers. 

Moreover it is natural i.e.
  there is 
  a functor $PreWtr\to \ShNW$ that 
    sends any presheave 
    to it's Nisnevich sheafication.
And 
  this functor is 
  left adjoin 
  to the embedding $\ShNW\to PreWtr$. 
\end{theorem}

Before start the proof 
we give discussion of properties of $Witt$-correspondence 
 in respect to Nisnevich coverings 
and give some useful definitions.  

\newcommand{\QS}[1]{QuadSpace(#1)}  
\begin{definition}
\label{defSuppQC}
Let 
  $(P,q_P)$ be quadratic space 
  in the category $Proj(X,Y)$.
We call by 
  \textbf{support} of $(P,q_P)$ 
  a closed subscheme of $X\times Y$
    corresponding to the 
    $Ann_{k[X\times Y]}\,P \subset k[X\times Y]$,
    where 
      $Ann_{k[X\times Y]}\,P\subset k[X\times Y]$ denotes 
      annihilator ideal of $P$ as $k[X\times Y]$-module, 
  i.e.
  $$Supp\,(P,q_P) = Spec k[X\times Y]/Ann_{k[X\times Y]}\,P.$$  
\end{definition}

\begin{definition}\label{defGoodLiftQS}
Let 
  $X$, $Y$ be smooth affine schemes and
  $u\colon U\to X$, $v\colon V\to Y$ be Nisnevich coverings. 
Let 
  $\Phi$ be quadratic space 'between' $X$, and $Y$ and 
  $\Psi$ its lift along $v$ and $u$, i.e.
  \begin{multline*}
  \Phi\in \QS{Proj(X,Y)},\;
  \Psi\in \QS{Proj(U,V)}\colon\\
   u^*(\Phi) = v_*(\Psi) \in \QS{Proj(U,Y)}.
  \end{multline*}
Then
  $\Psi$ is called a \textbf{'good' lift} 
  if 
    the morphism $v\times id_U\colon V\times U\to Y\times U$
    induces isomorphism $$Supp\,\Psi \simeq Supp\,u^*(\Phi).$$
\end{definition}

Any quadratic space is well defined 
  over the function ring of its support.
Now let's involve denotation for 
such induced quadratic space 
  on closed subscheme containing support.  
\begin{definition}[quadratic space $\Phi^Z$]
\label{LocQStoSup}
Let
  $p\colon S\to X$ be regular map of affine varieties,
  $\Phi=(P,q_P)$ be quadratic space in ${Proj(p)}$
  and 
  $i\colon Z\subset S$ be closed embedding
  and $Supp\Phi\subset Z$.
Then
  module $P$ and isomorphism $q_P$ are
    well defined over $k[Z]$.  
So they define
    a quadratic space $\Phi^{Z}$ in 
    the category $Proj(p\big|_{Z})$
    corresponding to restriction of $p$ onto $Z$
    $p\big|_{Z} \colon \pri S \to X$,
    such that
    $\Phi = i_*(\Phi^Z)$.   
\end{definition}

\begin{remark}\label{remGoodLiftQS}
'Good' lifts are such lifts 
  that are defined by 
  lifting of the support
  in following sense. 

{\sloppy 
Let
  $u\colon U\to X$ and $v\colon V\to Y$ 
  be Nisnevich coverings of affine varieties,
  $\Phi\in \QS{Proj(X,Y)}$, 
  $Z = Supp\, u^*(\Phi) $ and
  $\Phi^\prime = u^*(\Phi)^Z$
  (in sense of definition \ref{LocQStoSup}). 
Then 
  $\Psi\in \QS{Proj(U,V)}$ is 'good' lift of $\Phi$ 
if and only if     
  there is a lift 
    $l\colon Z \to V$
    of the morphism of projection 
    $ Z \to Y$
  such that
    $\Psi = g_*({\pri\Phi})$
    where 
    $g = l\times id_U\colon Z\hookrightarrow U\times V$
    (that is closed embedding because
      it is lift of closed embedding of $Z$ into $U\times Y$
      ). 

}      

In fact
  $l\colon Z \to V$ can be defined as 
    compositions of 
    canonical projections $Z\to V$
    with inverse to isomorphisms  $Z_\Psi \simeq Z$ 
      from definition of 'good' lift.
And conversely 
  if $l$ is such lift,
  then $Z\simeq g(Z)=  Supp\,\Psi$. 
$$\xymatrix{
V\ar[rr]^v & & Y
\\
Z_{\Psi}\ar[u]\ar[d]\ar@{=}[r] & 
Z\ar[ur]\ar[dl]\ar[ul]_l &
Z_\Phi\ar[u]\ar[d]
\\
U\ar[rr]_u & & X
}$$

\end{remark}

\begin{lemma}\label{finSuppQC}
For 
  any affine $X$, $Y$ and 
  any quadratic space $\Phi$ in $Proj(X,Y)$ 
$Supp\,\Phi$ is finite over $X$.
\end{lemma}
{\em Proof of the lemma.}
Let $\Phi = (P,q_P)$.
By definition 
  $$Supp\,\Phi = Spec k[X\times Y]/ Ann\,P.$$
Let      
  $m_1,\dots,m_n$ be finite set of generators of $P$ over $k[X]$
then
  $$Ann\,P=\bigcap\limits_{i=1}^n Ann\,m_i$$
and
  the homomorphism of $k[X]$-modules
  $$(1,\dots 1)\colon R/Ann\,P\to \bigoplus\limits_{i=1}^n R/Ann\,m_i$$
  is embedding.
The composition of 
  this homomorphism with 
  embeddings   $R/Ann\,m_i\hookrightarrow P$ 
gives us 
  the embedding of $R/Ann\,P$ into $\bigoplus\limits_{i=1}^n P$.
So $R/Ann\,P$ is isomorphic to the 
  submodule of the finitely generated $k[X]$-module
and since $k[X]$ is Noetherian
  $R/Ann\,P$ is finitely generated over $k[X]$. 
{\em Lemma is proved.}

\begin{corollary}\label{trNtopQC}
For
  any 
    local Hensel $U^h$,
    smooth affine $Y$ and 
  any $\Phi\in \QS{Proj(U^h,Y)}$
$$\Phi = \sum\limits_{i=1\dots n}\Phi_i$$
for some
  $\Phi\in \QS{Proj(U^h ,Y)}$ 
  such that each
  $\Phi_i$ is represented by 
    quadratic space $(P_i,q_{P_i})$
    such that its support is local Hensel scheme. 
\end{corollary}
{\em Proof of the lemma.}
This immediately follows from previous because
finite subscheme over local Hensel scheme splits into disjoint union of local Hensel subschemes.
{\em Lemma is proved}

\begin{lemma}\label{goodliftQS} 
{\sloppy
For 
  any smooth affine $X$, $Y$, 
  Nisnevich covering $v\colon V\to Y$ and 
  any $\Phi\in \QS{Proj(X,Y)}$
there are
  some Nisnevich covering 
    $u\colon U\to X$
  and $\Psi\in \QS{Proj(U,V)}$
  that is
  a 'good' lift of $\Phi$.

}  
\end{lemma}
{\em Proof of the lemma.}
We will construct covering $U$ of $X$ 
as disjoint union of 
  Nisnevich neighbourhood for all points $x\in X$.
I.e. $U\simeq \coprod\limits_{x\in X} U_x$
and for all $x\in X$ 
  there is a lift $l_x\colon x\to U_x$
  such that 
    composition $u\circ l_x\colon x\to X$ is equal to 
    the embedding of $x$ into $X$.

So let $x\in X$ be any point.
Let $u_x^h\colon U^h_x\to X$ be 
  Hensel neighbourhood of $X$ at $x$, i.e. 
  $k[U_x]$ is henselisation $O^h_{X,x}$  
    of local ring $O_{X,x}$ at $x$.

By lemma \ref{trNtopQC}
 support of $\Phi$ over $U^h_x$ 
 splits into finite disjoint union
   of Hensel local schemes 
   $$
   Z^h_x = Supp\;({(u^h_x)}^*(\Phi))
    = \coprod\limits_{y_i} Z^h_{x,y_i},$$
   where points $y_i\in Y$ are 
   the image along the projection on $Y$
     of the closed points 
       in the fibre $Z^h_x$ over $x$.

Since all $Z^h_{x,y_i}$ 
  are Hensel local and $V\to Y$ is Nisnevich covering
there are 
  lifts $l_{x,y_i}\colon Z_{x,y_i}\to V$ 
  of projections $Z_{x,y_i}\to Y$.
And this defines 
  lifts 
  $$g_{x,y_i}\colon Z_{x,y_i}\hookrightarrow U^h_x\times V$$
  of closed embeddings 
  $Z_{x,y_i}\hookrightarrow U^h_x\times Y$ .
  
Quadratic space
  $\Phi^h_x = {(u^h_x)}^*(\Phi)$ splits into direct sum
  $\Phi^h_x = \sum\limits_i \Phi^h_{x,y_i}$
  of quadratic spaces with supports $Z_{x,y_i}$.  
By remark \ref{remGoodLiftQS}
 lifts $l_{x,y_i}$ define
  a 'good' lift  
  of quadratic space $\Phi^h_x$
  $$
  \Psi^h_x = \sum\limits_i {g_{x,y_i}}_*(\Phi^h_{x,y_i})
   \in \QS{Proj(U^h_x, V)},
  $$
  i.e.
  \begin{equation}\label{glHenloc}
  {(u^h_x)}^*(\Phi) = v_*(\Psi^h_x),\quad  
  Supp\,\Psi^h_x \simeq  Z^h_x.  
  \end{equation}

Any quadratic space in $Proj(U^h_x, V)$ is 
germ of some quadratic space 
well defined over some affine scheme $U_x$.
So there are 
  affine schemes $U_x$ and
  quadratic spaces $\Psi_x$, 
  such that 
  $x$ is closed point of $U_x$, 
  $el_x\colon U^h_x\to U_x$ 
    is Nisnevich neighbourhood of $U_x$ at $x$,
  and 
$$
\Psi_x \in \QS{Proj(U_x,V)}\colon \Psi^h_x = {el_x}^*(\Psi_x)
.$$ 
To make the equalities
$$v_*(\Psi_x)\simeq u_x^*(\Phi),\quad  Supp\,\Psi_x \simeq Z_x$$
holds, 
it is enough to change $U_x$ to some its open subscheme,
because 
  over Nisnevich neighbourhood $U^h_x$ 
  equalities \eqref{glHenloc} holds.

{\em Lemma is proved.}

\begin{lemma}\label{rem_propGoodLift}
'Good' lifts are closed under 
 base changes and 
 shredding.
I.e. 
  if 
    $X_1$, $X_2$, $Y$ are any varieties 
    $f\colon X_1\to X_2$ is regular map 
    and 
    $v\colon V\to Y$,
    $u\colon U\to X_2$ and
    $\pri u\colon \pri U\to U$  
    are Nisnevich coverings
  then
    for any     
      $\Phi\in \QS{Proj(X_2,Y}$
      and its 'good' lift
      $\Psi\in \QS{Proj(U,V)}$,
    ${\pri u}^*\Psi$ is 'good' lift of $\Phi$ 
    and
    ${f_U}^*\Psi\in\QS{Proj(U\times_{X_2}X_1,V)}$ 
    is 'good' lift of $f^*(\Phi)$.
$$\xymatrix{
& V\ar[r]^v& Y\\
\pri U\ar[r]^{\pri u}\ar[ur]^{{\pri u}^*\Psi} & 
U\ar[r]^u\ar[u]_{\Psi} & X_2\ar[u]_{\Phi}\\
& U_{X_1}\ar[u]^{f_U}\ar[r]_{X_1\times_{X_2}u}&
X_2\ar[u]_f 
}$$      
  
And it is closed under compositions in week sense, i.e.
if 
  $X_1$, $X_2$, $Y$ are any varieties
  and $v\colon V\to Y$ is Nisnevich covering
then 
  for any 
    $\Phi_1\in \QS{Proj(X_1,X_2}$, 
    $\Phi_2\in \QS{Proj(X_2,Y}$
  there are
    coverings 
      $u_1\colon U_1\colon X_1$,
      $u_2\colon U_2\colon X_2$
    and
    'good' lifts
      $\Psi_1\in \QS{Proj(U_1,U_2}$,   
      $\Psi_2\in \QS{Proj(U_2,V}$,
  such that 
    $\Psi_2\circ\Psi_1$ 
    is 'good' covering of 
    $\Phi_2\circ\Phi_1$.
$$\xymatrix{
U_1\ar[d]^{u_1}\ar[r]^{\Psi_1} & 
U_2\ar[d]^{u_2}\ar[r]^{\Psi_2} & 
V\ar[d]^v
\\
X_1\ar[r]^{\Phi_1} & X_2\ar[r]^{\Phi21} & Y
}$$

\end{lemma}

Instead of  
the proof let's give following remark
\begin{remark}
Statement about pomposition with Nistevich covering  $\pri u\colon\pri U\to U$ and morphisms $f\colon X_2\to X_1$ 
follows immideately from definition 
and 
compatability with composition
can be proved similar to lemma \ref{trNtopQC}
by lifting of the diagram of composition firestly at local Nisnevich neigbourhhos
and then continuation of it to some covering.
Let's note in addition
 that to get a 'good' lift over local Hensel scheme
  we use fixed choise  
of lifts of images points along covering of $Y$ ,  
i.e. points $y_i$ in proof of lemma \ref{trNtopQC}
for all 'preimage' points in $X_1$ and $X_2$.

\end{remark}

\begin{lemma}\label{liftSum}
For any $\Phi_1,\Phi_2\in \QS{Proj(X,Y)}$ 
there are 
  Nisnevich coverings $
    u\colon U\to X$, $v\colon V\to Y$
  and 
  'good' lifts
    $\Psi_1,\Psi_2\in \QS{Proj(U,V)}$
  of 
    $\Phi_1,\Phi_2$  
such that
  $\Psi_1\oplus\Psi_2$ is 'good' lift of $\Phi_1\oplus\Phi_2$.

\end{lemma}
{\em Proof of the lemma.}
By lemma \ref{goodliftQS}
  for some coverings $u\colon U\to X$, $v\colon V\to Y$
  there is a 'good' lift 
  $\Psi_+\in QuadS(Proj(U,V))$ of $\Phi_+ = \Phi_1\oplus\Phi_2$
and by terms of the remark \ref{remGoodLiftQS} 
  it corresponds to some lift 
  $l_+\colon Supp\,u^*(\Phi_+) \hookrightarrow U\times V$. 
Since 
  $Supp\,u^*(\Phi_1)\subset Supp\,u^*(\Phi_+)$
quadratic spaces $u^*(\Phi_1)$, $u^*(\Phi_2)$ 
are well defined over $k[Supp\,u^*(\Phi_+)]$,
i.e. they defines some 
  $\Phi_1^\prime,\Phi_2^\prime\in QuadS(Proj(s_+))$
  where $s_+$ is projection of $Supp\,u^*(\Phi_+)$ to $U$.
Then 
  ${l_+}_*(\Phi_1^\prime)$ and ${l_+}_*(\Phi_1^\prime)$
  are the required 'good' lifts.      
{\em Lemma is proved.}

\begin{lemma}\label{metabolicGoodLift}
Let 
  $\Psi\in \QS{Proj(U,V)}$ be a 'good' lift of 
  $\Phi\in \QS{Proj(X,Y)}$
    along Nisnevich coverings 
    $u\colon U\to X$, $v\colon V\to Y$.
Then 
  if $\Phi$ is metabolic then $\Psi$ is metabolic too.     
\end{lemma}
{\em Proof of the lemma.}
Since $\Phi$ is metabolic 
$u^*(\Phi)$ is metabolic.
Then 
  $\pri \Phi  = u^*(\Phi)^Z$
  is metabolic
  (where $Z = Supp\, u^*(\Phi)$).
And since in terms of  the remark \ref{remGoodLiftQS} 
$\Psi  = g^*(\pri\Phi)$,
$\Psi$ is metabolic.    
{\em Lemma is proved.}

\begin{definition}\label{productGoodLiftQS}
Let's define a product for 'good' lifts of quadratic spaces,
i.e. an operation 
  that 
  for two 'good' lifts  
  $$\Psi_1,\Psi_2\in QuadS(Proj(U,V))$$ 
    of quadratic space $\Phi\in Quad(Proj(X,Y))$
    along Nisnevich coverings $u\colon U\to X$ and $v\colon V\to Y$,
  gives
    a 'good' lift 
    $$\Psi_1\times_{\Phi}\Psi_2\in QuadS(Proj(U,V\times_Y V)) $$
    along coverings $u$ and $v\times_Y v$ 
    such that
    $$pr_i\circ \Psi_1\times_{\Phi}\Psi_2= \Psi_i,\;i=1,2$$
  $$
  \times_{\Phi}\colon
    (\Psi_1,\Psi_2)\mapsto \Psi_1\times_\Phi \Psi_2  $$
  $$\xymatrix{
  V\times_Y V \ar@<-0.5ex>[r]_{pr_1}\ar@<-0.5ex>[r]^{pr_2} &
    V\ar[r]^v & 
      Y\\
  & U\ar[r]_u
     \ar@<-0.5ex>[u]_{\Psi_1}\ar@<-0.5ex>[u]^{\Psi_2}
     \ar[ul]^{\Psi_1\times_{\Phi}\Psi_2} & 
      X\ar[u]_{\Phi}
  }$$    
  
To construct it let's use 
  descriptions of 'good' lifts 
  giver in remark \ref{remGoodLiftQS}.
Let 
  $Z = Supp\,(\Phi\circ u)\subset U\times Y$,
  $l_i\colon Z\to V$ 
    be its lifts corresponding to $\Psi_1$ and $\Psi_2$.
Then 
  the product 
    $\Psi_1\times_{\Phi,V} \Psi_2$    
    of lifts of quadratic spaces
  is defined by the product of lift of their supports 
  $(l_1,l_2)\colon Z\to V\times_Y V$.
  I.e. 
  $\Psi_1\times_{\Phi}\Psi_2$            
  is 
  direct image 
    of $\Phi^\prime$ 
    along $(l_1,l_2) \colon Z \to V\times_Y V $. 
  where
  $\Phi^\prime = \Phi\circ u\in QuadS(Proj(U,Y))$.
    
$$\xymatrix{
V\times_Y V \ar@<-0.5ex>[r]_{pr_1}\ar@<-0.5ex>[r]^{pr_2} &
  V\ar[r]^v & 
    Y\\
& Z\ar[d]|p
   \ar@<-0.5ex>[u]_{l_1}\ar@<0.5ex>[u]^{l_2}
   \ar[lu]|{(l_1, l_2)}
   \ar[r] & 
     Supp\,\Phi\ar[u]\ar@<-0.5ex>[d]\\   
& U\ar[r]_u
   \ar@<-1ex>[u]_{\Psi_1}\ar@<+1ex>[u]^{\Psi_2}
   \ar[ulu]^{\Psi_1\times_{\Phi}\Psi_2} & 
      X\ar@<-0.5ex>[u]_{\Phi}
}$$ 
\end{definition}

\vspace{10pt} 
{\em Proof of the theorem \ref{TrForTheSheaf}.}

Let $a\in \bF_{Nis}(Y)$. 
Let $a$ be represented by $\tilde a\in \bF(V)$, 
  i.e. $\epsilon(\tilde a)=v^*(a)$
  for some Nisnevich covering $v\colon V \to Y$.    
Let
  $\Phi$ be some quadratic space in $Proj(X,Y)$.
By lemma \ref{goodliftQS}  
  for some $u\colon U\to X$
exists a 'good' lift $\Psi\in QuadS(Proj(U,V))$.
If 
  sought-for structure 
    of presheave with $Witt$-transfers on $\bF_{Nis}$
  exists
then 
  $u^*(\Phi^*(a))$ should be equal to $\epsilon(\Psi^*(\tilde a)$.
So we want to put $\Phi^*(a)$
to be the element $\bF_{Nis}(X)$
represented by the element $\Psi^*(\tilde a)$ 
  in the group of sections of presheave $\bF(U)$.

\newcommand{\prUf}{{pr_{u,1}}}
\newcommand{\prUs}{{pr_{u,2}}}
\newcommand{\prUi}{{pr_{u,i}}}
\newcommand{\prVf}{{pr_{v,1}}}
\newcommand{\prVs}{{pr_{v,2}}}
\newcommand{\prVi}{{pr_{v,i}}}
To do it 
first of all we should check that
$\Psi^*(\tilde a)$ defines the section of the sheaf,
i.e. that 
  $${\prUf}^*(\Psi^*(\tilde a)) = {\prUs}^*(\Psi^*(\tilde a)),$$
where $\prUi\colon U\times_X U\to U$ are canonical projections.
By the remark \ref{rem_propGoodLift} 
the compositions 
 $${\prUi}^*(\Psi)\in QuadS(Proj(U\times_X U,V)),\; i=1,2$$
are two 'good' lifts of $\Phi$
  along Nisnevich coverings 
    $v\colon V\to Y$ and $u^2\colon U\times_X U\to X$.         
So due to construction from definition \ref{productGoodLiftQS}
 there is  
  $$\Psi_3 = 
   {\prUf}^*(\Psi)\times_{\Phi} {\prUs}^*(\Psi) 
   \in QuadS(Proj(U\times_X U,V\times_Y V))\colon\quad
    {\prVi}_*\Psi_3 = {\prUf}^*(\Psi).$$
$$\xymatrix{
V\times_Y V
  \ar@<0.5ex>[r]^{\prVs}\ar@<-0.5ex>[r]_{\prVf}
     &V\ar[r]^v &Y \\
U\times_X U
  \ar[u]^{\Psi_3}
  \ar@<0.5ex>[r]^{\prUs}\ar@<-0.5ex>[r]_{\prUf} 
    &U
      \ar[u]^{\Psi}
      \ar[r]_u 
        &X\ar[u]_{\Phi}
}$$
Then 
  \begin{equation*}
  {\prUf}^*(\Psi^*(\tilde a)) =
  \Psi_3^*({\prVf}^*(\tilde a)) = 
  \Psi_3^*({\prVs}^*(\tilde a)) = 
  {\prUs}^*(\Psi^*(\tilde a)).
  \end{equation*}

Thus for any $\Phi\in QuadS(Proj(X,Y))$
we start with a section $a\in \bF_{Nis}(Y)$
and
construct a section $\pri a\in \bF_{Nis}(X)$ 
  of the sheaf $\bF_{Nis}$ on $X$.

Let's summarize 
  additional data used in this construction.
For any set 
 $$\begin{aligned}\mathcal D =
 (\;v\colon V\to Y,\tilde a\in \bF(V),u\colon U\to X,
  \Psi\in \QS{Proj(U,V)}
  \;)\colon \\
  v^*(a)=\epsilon(\tilde a),\quad
  u^*(\Phi) = v_*(\Psi),\quad
  \Psi \text{ is 'good' lift of } \Phi
 \end{aligned}$$
we construct    
  $$
  \Phi^*_{\mathcal D}(a) \in \bF_{Nis}(X)\colon\quad
    u^*(\Phi^*_{\mathcal D}(a)) = \epsilon(\Psi(v^*( a )))\quad
    \; (\Phi^*_{\mathcal D}(a) = \pri a)
    .
    $$

So to get well defined map 
  $\Phi^*\colon \bF(Y)\to \bF(X)$
we should check that 
  this construction doesn't depend on 
  additional data $\mathfrak D$.

Firstly we check \textbf{independence on 
  the choice of the lift} $\Psi$
   for fixed $u$, $v$ and  $\tilde a$.   
Let $\Psi_1,\Psi_2\in QuadS(Proj(U,V))$ are two 'good' lifts of $\Phi$.
The definition \ref{productGoodLiftQS} provides
a quadratic space
  \begin{gather*}
  \Psi_1\times_{\Phi}\Psi_2\in QuadS(Proj(U,V\times_Y V))\colon\\
  {pr_i}_*((\Psi_1\times_{\Phi}\Psi_2)) = \Psi_i,\,i=1,2
  \end{gather*}
  where $pr_i\colon V\times_Y V\to V$ are canonical projections.
$$\xymatrix{
V\times_Y V
  \ar@<0.5ex>[r]^{pr_1}\ar@<-0.5ex>[r]_{pr_2}
  & V
    \ar[r]
    & Y\\
& U
  \ar[lu]
  \ar@<0.5ex>[u]^{\Psi_1}\ar@<-0.5ex>[u]_{\Psi_2}
  \ar[r]^u
  & X\ar[u]_{\Phi}    
}$$

Then 
  \begin{equation*}
  \Psi_1^*(\tilde a) = 
  (\Psi_1\times_{\Phi}\Psi_2)^*(pr_1^*(\tilde a)) = 
  (\Psi_1\times_{\Phi}\Psi_2)^*(pr_2^*(\tilde a)) = 
  \Psi_2^*(\tilde a)
  \end{equation*}
and
  $$  
  u^*\Phi^*_{\Psi_1}(a) =
  \epsilon(\Psi_1^*(\tilde a)) =
  \epsilon(\Psi_2^*(\tilde a)) =
  u^*\Phi^*_{\Psi_2}(a). 
  $$         
So
  $$\Phi^*_{\Psi_1}(a) = \Phi^*_{\Psi_2}(a).$$  
 
Next we check \textbf{independence on the covering} $u$. 
Let 
  $u_1\colon U_1\to X$, $u_2\colon U_2\to X$
    are two Nisnevich covering and
  $\Psi_i\in QuadS(Proj(U_i,Y)),\,i=1,2$
    are 'good' lifts of $\Phi$.
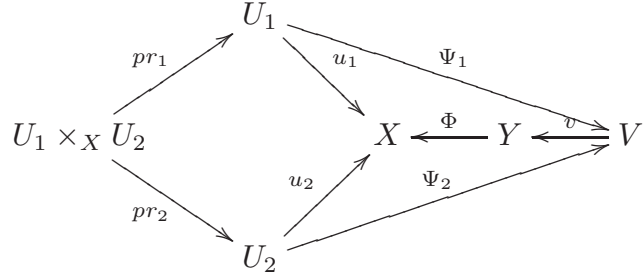
$$\xymatrix{
& U_1
  \ar[rd]^{u_1}
  \ar[rrrd]^{\Psi_1}
  & &&\\
U_1\times_X U_2
  \ar[ru]^{pr_1}\ar[rd]_{pr_2}
  & & X & Y\ar[l]_{\Phi} & V\ar[l]_v\\
& U_2
  \ar[ru]^{u_2}
  \ar[rrru]^{\Psi_2}
   & &&
}$$ 
The compositions 
  $$pr_1^*(\Psi_1),\,pr_2^*(\Psi_2) 
  \in QuadS(Proj(U_1\times_X U_2,V))$$
are two 'good' lifts of $\Phi$ along covering $U_1\times_X U_2\to X$
and 
  $\Phi_{U_1\times U_2,pr_i^*(\Psi_i)}^*(a)=\Phi_{U_i,\Psi_i}^*(a)$
because
  $(pr_i^*(\Psi_i))^*(\tilde a) = pr_i^*(\Psi_i^*(\tilde a))$.  
Thus by independence on choice of lift for fixed covering
$$\Phi^*_{U_1}(a) = 
  \Phi^*_{U_1\times U_2,pr_i^*(\Psi_i)}(a) 
  \Phi^*_{U_2}(a).$$

The \textbf{independence on 
  the choice 
   of the covering  $V$ and 
   representation $\tilde a$}
we check in two steps.
First we note that
  for any Nisnevich covering $v^\prime\colon V^\prime \to V$
  $$
  \Phi^*_{V^\prime,{v^\prime}^*(\tilde a)}(a) = \Phi^*_{V,\tilde a}(a)
  .$$ 
Because
  for some Nisnevich covering $u^\prime\colon U^\prime\to U$
  exists a 'good' lift $\Psi^\prime$ of $\Psi$ along $v^\prime$.
And 
  by commutativity of the diagram
  $$\xymatrix{
  V^\prime\ar[r]^{v^\prime} & V\ar[r]^v & Y\\
  U^\prime
  \ar[r]_{u^\prime}\ar[u]^{\Psi^\prime} & 
    U
    \ar[r]_u\ar[u]_{\Psi} & 
      X\ar[u]_{\Phi}
  }$$        
  $$
  (u\circ u^\prime)^*(\Phi^*_{V,\tilde a}(a)) = 
  \epsilon({u^\prime}^*(\Psi^*(\tilde a)) = 
  \epsilon({\Psi^\prime}^*({v^\prime}^*(\tilde a))) = 
  (u\circ u^\prime)^*(\Phi^*_{V^\prime, {v^\prime}^*(\tilde a)}(a)).
  $$ 

Secondly for any two representations of section $a$ 
  along  two Nisnevich coverings $v_i\colon V_i\to Y,\, i=1,2$
  $$
  {\tilde a}_1\in \bF(V_1),\; {\tilde a}_2\in \bF(V_2)\colon\quad
  {v_1}^*(a)= \epsilon({\tilde a}_1),\; 
  {v_2}^*(a) =  \epsilon({\tilde a}_2)
  $$ 
there is 
  common shredding 
  \begin{gather*}\xymatrix{
  V\ar[r]\ar[d] & V_1\ar[d] \\
  V_2\ar[r] & Y }\\
  {v^\prime_1}^*({\tilde a}_1) = 
    \tilde a = {v^\prime_2}^*({\tilde a}_2).
  \end{gather*}
So by discussion in previous paragraph 
  $$\Phi^*_{V_1,{\tilde a}_1}(a) = 
    \Phi^*_{V,{\tilde a}}(a) = 
      \Phi^*_{V_2,{\tilde a}_2}(a) .$$

Thus 
  for any smooth affine $X$ and $Y$ 
  and $\Phi\in QuadS(Proj(X,Y))$
we get a well defined map 
  $\Phi^*\colon \bF_{Nis}(Y) \to \bF_{Nis}(X)$.
To finish the proof we should check following.
\begin{itemize}
\item[1)]
  This maps defines
     additive homomorphisms 
     $WCor(X,Y) \to Hom(\bF(X),\bF(Y)$.
\item[2)]
  This homomorphisms
     combines the functor $\bF_{Nis}\colon WCor\to Ab$
     such that 
       its composition with canonical functor $Sm\to WCor$
       is naturally equal to $\bF_{Nis}\colon Sm\to Ab$.
\item[3)] 
  This construction
    is natural on the presheave $\bF$. 
  I.e
  for a morphism of presheaves $s\colon \bF_1\to\bF_2$            
  the morphism $s_{Nis}\colon{\bF_1}_{Nis}\to {\bF_2}_{Nis}$
  becomes a natural homomorphism of functors from $WCor$ to $Ab$. 
\end{itemize}

Or more precisely:
\begin{itemize}
\item[1.1)]
  $\Phi^*$ is additive for any $\Phi\in QuadS(Proj(X,Y))$,
\item[1.2)]
  $(\Phi_1\oplus \Phi_2)^* = \Phi_1^* + \Phi_2^*$
  for any $\Phi_1,\Phi_2\in QuadS(Proj(X,Y))$,
\item[1.3)]
  $\Phi^* = 0$ for any metabolic quadratic space $\Phi$. 
\item[2.1)]     
  $\Phi_1^*\circ\Phi^*_2 = (\Phi_2\circ\Phi_2)^*$ 
    for $\Phi_1\in WCor(X,Y)$ and $\Phi_2 \in WCor(Y,Z)$,
\item[2.2)]
  for $Witt$-correspondence $\Phi\in WCor(X,Y)$
  defined by regular map $f\colon X\to Y$.
    $\Phi^* = f^*.$
\item[3)]
  The diagram
  $$\xymatrix{
  \bF_2(X)\ar[r]^{\Phi^*}                 & \bF_2(Y) \\
  \bF_1(X)\ar[r]^{\Phi^*}\ar[u]^{s_{Nis}} & \bF_1(Y) \ar[u]_{s_{Nis}} 
  }$$
    for any
      $\Phi\in WCor(X,Y)$ and
      morphism of presheaves $s\colon \bF_1\to\bF_2$.
\end{itemize}        

We don't not give detailed checking for all points 
but write down it for point 1)
and explain general scheme of such proof.

Let $a_1,a_2\in \bF_{Nis}(Y)$.
Let's choose 
  a covering $v\colon V\to Y$
  such that $v^*(a_i)=\epsilon({\tilde a}_i),\,i=1,2$.
It can be done by choosing of the coverings
  $v_i\colon V_i\to Y$ 
  such that $v_i^*(a_i)=\epsilon({\tilde a}_i),\,i=1,2$
  and 
  putting $V$ to be a product of $V_i$.
Then by lemma \ref{goodliftQS}  
  for some Nisnevich cover $u\colon U\to X$ there is
  a 'good' lift $\Psi\in WCor(U,V)$ of $\Phi$
and 
we can use it to define 
  $\Phi^*(a_1)$, $\Phi^*(a_2)$ and $\Phi^*(a_1+a_2)$
  to be a sections of $\bF_{Nis}$ over $X$ that are
    represented by sections
  $\Psi^*({\tilde a}_1)$, $\Psi^*({\tilde a}_2)$ and 
    $\Psi^*({\tilde a}_1+{\tilde a}_2)$
  of $\bF$ over $U$.
Thus 
  $$u^*(\Phi^*(a_1)+\Phi^*(a_2)) =
    \epsilon(\Psi({\tilde a}_1)+\Psi({\tilde a}_2)) =
    \epsilon(\Psi({\tilde a}_1+{\tilde a}_2)) = 
    u^*(\Phi^*(a_1+a_2)).$$ 

Similary for
  any property
  of reverse images for presheaves 
  $\Phi^*\colon\bF(X)\to\bF(Y)$
we can   
  transfer it to the property 
  of reverse images
  $\Phi^*\colon\bF_{Nis}(X)\to\bF_{Nis}(Y)$
using 
  given upper construction of $\Phi^*$ and 
  proved independence on choice of additional data
by the following scheme.
    
First we 
  lift the diagram of 
    corresponding property for preshaves
  by consequently choosing 
   compatible 
     Nisnevich coverings for all schemes and
     lifts for all objects in the formulation of the property
       like as morphisms and sections.
The compatibility means that
  the property holds for the listed objects.
Then we use that fact that 
  lifted objects uniquely determine 
  corresponding objects for the sheaves.

Also let's note that point 1.2) uses lemma \ref{liftSum}
and 1.3) uses lemma \ref{metabolicGoodLift}.       
{\em The theorem \ref{TrForTheSheaf} is proved.}      
\vspace{10pt} 

Now we go down to consideration of Nisnevich cohomology groups of
the seheave with transfers.
One of our aims is to prove 'strictly' variant of theorem \ref{TrForTheSheaf}, i.e.
\begin{theorem}\label{StrShWtrTr}
For any Nisnevich sheave $\bF$
presheaves of cohomologyes $H^i_{Nis}(\bF)$
are equipped with canonical structure of presheave with $Witt$-transfers.
\end{theorem}
In fact we get it as immediate corollary of 
following theorem. 

\begin{theorem}\label{ReprChWtr}
There is natural isomorphism
$$Ext^i_{\ShNW}(\mathbb Z_{Wtr}(X),\bF) \simeq H^i_{Nis}(X,\bF)$$
for all smooth affine $X$ and sheaves with $Witt$-transfers $\bF$.
\end{theorem}
Note that formulation of the last theorem use that category $\ShNW$
is abelian that follows from the theorem \ref{TrForTheSheaf}.

The last theorem in turn is  
particular case of 
  adjacency isomorphism
  of derived functors of adjoin functors 
$$ShNis\rightleftarrows \ShNW$$ 
  that equips the sheaves 
  with the structure of sheave with $Witt$-transfers
  and forgetfull functor,
and  
  the proof of theorem \ref{ReprChWtr} is given in the end of the section. 

Let's involve following notations.

\begin{definition}\label{defWtrFunc}
Let's denote by 
$$
{Wtr_{Pre}}\colon Pre \to PreWtr,\quad
 Wtr_{ShN}\colon ShNis\to \ShNW
$$ 
the left Kan extension functor 
  along the functor 
    from additivisation of category $Sm_k$ into
     category $WCor_k$,
i.e. the functor that equips any presheave
with $Witt$-transfers by universal way, 
and
the the functor that is composition of 
  the functor 
  of the embedding $ShNis\to Pre$,
  the functor $Wtr_{Pre}$
  and 
  Nisnevich sheafication functor.

Then let's denote by
\begin{gather*}
L(Wtr_{Pre})\colon D^-(Pre)   \to D^-(PreWtr),\quad
L(Wtr_{ShN})\colon D^-(NisSh) \to D^-(\ShNW)
\end{gather*}
the left derived functors of $Wtr_{Pre}$ and $Wtr_{ShN}$. 

\end{definition}

\begin{remark}\label{adjWtr}
The functors
  $Wtr_{Pre}$ and $Wtr_{ShN}$ are left adjoin
  to the forgetful functors
  $F_{Wtr}\colon PreWtr\to Pre$ and 
  $F_{Wtr}\colon \ShNW\to ShNis$.
\end{remark}

\begin{remark}\label{adjLdWtr}
Forgetful functor $F_{Wtr}\colon \ShNW\colon ShNis$
is exact, hence it induce a functor between derived categories
that is booth left and right derived functor for $F_{Wtr}$.
Let's denote it by the same symbol.
Then 
since
left and right derived functors 
to the left and right adjoin functors 
are adjoin again,
there is an adjacency  
\begin{equation}\label{ShNDWtrAdj}\begin{aligned}
L{Wtr}\colon D^-(NisSh) &\rightleftarrows& D^-(\ShNW)\colon F_{Wtr} 
&\; &\\
L_{Wtr}&\dashv&F_{Wtr}\qquad\qquad\qquad\qquad
&.\;& 
\end{aligned}\end{equation}
\end{remark}
 
{\sloppy
Now to apply the adjucency \eqref{ShNDWtrAdj}
we want partly calculate the functor 
$L_{Wtr}\colon \DShNWtr\to \DShNis$.
In fact it happens that 
$L(Wtr_{ShN})$ coincides with $L(Wtr_{Pre})$.
because  
derived functor $L(Wtr_{Pre})$ is exact 
in respect to
exact sequences
that goes from Nisnevich topology structure.

}
More precisely it means the following. 
Let's identify by standard equivalence
the categories
$\DShNis$ and $\DShNWtr$
with localisations of categories
$\DPre$ and $\DPreWtr$
by Nisnevich sheaf quasi-isomorphisms 
that are the following.  

\begin{definition}
The morphism $q\colon A^\bullet \to B^\bullet$
  in category $D^-(Pre)$ or $D^-(PreWtr)$
is called sheaf quasi-isomorphism if 
homomorphism 
  $\underline h^i(q)\colon 
    \underline h^i(A^\bullet) \to \underline h^i(B^\bullet)$
  of Nisnevich sheaf cohomology, 
    i.e. Nisnevich shafication of  cohomology presheaf 
are isomorphisms.

The complex $A^\bullet$ 
  in category $D^-(Pre)$ or $D^-(PreWtr)$ 
is called sheaf acyclic if
all Nisnevich sheaves $\underline h^i(A^\bullet)$ 
are zero.
\end{definition}

Notice that this identification uses theorem \ref{TrForTheSheaf},
because by the definition category $\ShNW$
is a full subcategory of $PreWtr$
but due to theorem \ref{TrForTheSheaf} it 
is equal to 
  localisation of  $PreWtr$ at (Nisnevich-)local isomorphisms 
and there are commutative diagrams
\newcommand{\WtrNis}{Wtr_{Nis}}
\newcommand{\WtrPre}{Wtr_{Pre}}
\newcommand{\lNis}{l_{Nis}} 
\newcommand{\LWtr}{L(Wtr)}
\newcommand{\LPWtr}{L_{Pre}(Wtr)}
\newcommand{\LNWtr}{L_{Nis}(Wtr)}
\newcommand{\LWtrP}{L_{Pre}(Wtr)}
\newcommand{\LWtrN}{L_{Nis}(Wtr)}
\newcommand{\RFWtr}{F_{Wtr}}
\newcommand{\RFWtrP}{F_{Wtr}}
\newcommand{\RFWtrN}{F_{Wtr}}
\newcommand{\FWtr}{F_{Wtr}}
\newcommand{\FWtrN}{F_{Wtr}}
\newcommand{\FNis}{F_{Nis}}
\newcommand{\FNisW}{F_{Nis}}
\newcommand{\RFNis}{h_{Nis}}
\newcommand{\RFNisW}{h_{Nis}}
 \newcommand{\DShNW}{D^{\ShNW}}  
$$\xymatrix{
\ShNW\ar@<1ex>[r]^{\FNis}   \ar@<1ex>[d]^{\FWtrN} & 
PreWtr\ar@<1ex>[l]^{Nis}  \ar@<1ex>[d]^\FWtr 
 &&
\DShNW\ar@<1ex>[r]^{\RFNisW}\ar@<1ex>[d]^\RFWtrN & 
\DPreWtr\ar@<1ex>[l]^{\lNis}\ar@<1ex>[d]^\FWtr 
\\
ShNis\ar@<1ex>[u]^{\WtrNis}\ar@<1ex>[r]^\FNis & 
Pre \ar@<1ex>[u]^{\WtrPre}\ar@<1ex>[l]^{Nis} 
&&
\DShNis\ar@<1ex>[u]^{\LNWtr}\ar@<1ex>[r]^\RFNis & 
\DPre \ar@<1ex>[u]^{\LPWtr}\ar@<1ex>[l]^{\lNis}
}$$

\begin{theorem}\label{WtrNTopEx}
The functor
$L(Wtr_{Pre})\colon D^-(Pre)\to D^-(PreWtr)$
preserves Nisnevich sheaf quasi-isomorphisms.
\end{theorem}
{\em Proof of the theorem.\par}

As mentioned above 
this result follows 
  from that fact that
any $Witt$-correspondence from local Hensel scheme
splits into the sum of
  $Witt$-correspondence between local Hensel schemes.
To prove this theorem we use following modification of it.
\begin{lemma}\label{splitWCor}
For 
  any local Hensel scheme $U^h$ and affine $Y$
  there is a canonical decomposition
    $$WCor(U^h,Y) = \bigoplus\limits_z WCor_z(U^h,Y)$$
    where $z$ runs throw all closed points of $Y_x$
    and $x$ denotes closed point of $U^h$.
And
  component $WCor_z(U^h,Y)$ depends only on
  etale neighbourhood 
    of closed point $y\in Y$ that is 
       projection on $Y$ of the point $z$.
I.e.    
  for any etale neighbourhood $$v\colon(V,y)\to (Y,y),$$ 
  morphism $v$ induce isomorphism  
  $$WCor_z(U^h,V)\simeq WCor_z(U^h,Y).$$           
\end{lemma}
{\em Proof of the lemma.}\par
In fact in proof of the lemma \ref{finSuppQC} and lemma \ref{trNtopQC}.
it was shown that
  support of any object in $Proj(U^h,Y)$
  as closed subscheme in $U^h\times Y$
  splits into disjoint union of
  local Hensel schemes.
Any local Hensel subscheme in $U^h\times Y$
  has an unique closed point and
  all this points lies in $x\times Y$. 
Since 
  local Hensel subschemes with different closed points
  don't intersect, 
$Hom$-groups between the objects with such supports are zero.

Hence 
  the category $Proj(U^h,Y)$ 
  splits into direct sum of subcategories
  $Proj_z(U^h,Y)$ consisting of 
    the modules that
    support are local Hensel scheme 
      with closed point $z\in x\times Y$.
This induce required splitting of Witt groups of this categories
  $WCor(U^h,Y) = \bigoplus\limits_z WCor_z(U^h,Y)$.

The isomorphism 
$WCor_z(U^h,V)\simeq WCor_z(U^h,Y)$
along etale neighbourhoods $(V,y)\colon (Y,y)$
holds on the level of subcategories too
and
follows from the existence of the lifts of 
morphisms from local Hensel schemas along etale neighbourhoods.

In fact 
  for any local Hensel subscheme 
    $Z\subset U^h\times Y$ finite over $U^h$ 
    and such that
    $z\in Z$
  there is a lift
  $l\colon Z\to U^h\times V$
  and any local Hensel subscheme 
    $Z\subset U^h\times V$ finite over $U^h$
    and such that $z\in Z$
  is isomorphic to the image $v(Z)$.  
This induce equivalences of subcategories 
  \begin{gather*}
  Proj_{Z}(U^h,Y)\simeq Proj_{l(Z)}(U^h,V),\;Z\subset U^h\times Y\\
  Proj_{v(Z)}(U^h,Y)\simeq Proj_{Z}(U^h,V),\;Z\subset U^h\times V
  \end{gather*}
consists of modules
  that supports containing in specified closed subschemes.
And taking the limit along 
  closed local Hensel subschemes $Z$ 
  finite over $U^h$ and 
  containing $y$
we get the equivalences
 $$Proj_z(U^h,Y)\simeq Proj_z(U^h,V).$$

{\em Lemma is proved.}

\begin{definition}
Since by
previous lemma the group $WCor_z(U^h,Y)$ 
doesn't change after changing of $Y$ by its 
etale neighbourhood at the point $y$ that is
  projection of the point $z$,
we can define the group 
$WCor_z(U^h,V^h)$ 
where $(V^h,y)$ is Henselisation of $Y$ at $y$.
And let's define 
$WCor(U^h,V^h) = \bigoplus\limits_{z\times x\times y}WCor_z(U^h,Y^h)$.

Let's note that 
it is possible to give strong definition
of the groups $WCor(U^h,V^h)$
as $Hom$-groups between the local Hensel schemes
(as Pro-objects of the category $WCor$) 
but here $WCor_z(U^h,V^h)$ and $WCor(U^h,V^h)$
are only notations.
\end{definition}

Preservation of sheaf quasi-isomorphisms
is equivalent to the preservation of
sheaf acyclicity.
i.e. that
$L(Wtr_{Pre})$ sends
Nisnevich acyclic complexes in the category $D^-(Pre)$
to a Nisnevich acyclic complexes in  $D^-(PreWtr)$. 
So we show that
$L(Wtr_{Pre})$
  sends
  full triangulated subcategory of Nisnevich sheaf acyclic complexes 
  in $D^-(Pre)$ 
  into
  full triangulated subcategory of Nisnevich sheaf acyclic complexes 
  in $D^-(PreWtr)$.

Now
  let's note that 
    left derived functors 
    on the category $\DPreWtr$
    can be computed 
    by using of resolvents 
      consisting of
      direct sums of representable presheaves 
      in $\mathbb Z(Sm)$,
      i.e. presheave of free abelian groups
           corresponding to 
           representable presheave of sets
             in $Sm$.

In fact
  any presheave is direct limit of representable ones,
  so for any complex in $\DPreWtr$ there is such resolvent.
  And since
  representable presheaves are projective objects in category $Pre$, 
     for any complex $\bF^\bullet\in\DPreWtr$
     consisting of sums of representable presheaves
     $L(Wtr)(\bF^\bullet) = Wtr_{Pre}(\bF^\bullet)$,
     i.e. to compute $L(Wtr)$ we can apply $Wtr_{Pre}$ to each member.

Thus 
  since subcategory of Nisnevich acyclic complexes
  is closed under quasi-isomorphisms
to prove the  theorem 
it is enough
to prove that for 
  any complex of representable presheaves $\bF^\bullet$
  if $\bF^\bullet$ is Nisnecich acyclic
  then $Wtr_{Pre}(\bF^\bullet)$ is Nisnevich acyclic too.
  
To finish the proof  let's note that 
Nisnevich acyclicity is equivalent to 
  acyclicity of the complexes of germs on Hensel local schemes
    $Wtr_{Pre}(\bF^\bullet)(U^h)$.
  and that  
  lemma \ref{splitWCor}
  allows to express
    $Wtr_{Pre}(\bF^\bullet)(U^h)$ in terms of $\bF(U^h)$
    for any presheave $\bF$ and local Hensel $U^h$.
In fact by the formula for left Kan extensions
  $$
  Wtr_{Pre}(\bF)(U) \simeq
  \varinjlim\limits_{{f\colon Y_1\to Y_2}}
    \bF(Y_2)\otimes WCor(U,Y_1).
  $$
By lemma \ref{splitWCor} for any local Hensel $U^h$ and any affine $Y$
$$WCor(U^h,Y) = \sum\limits_y WCor_y(U^h,Y).$$
Therefore by 
  including of 
  direct sums into direct limit
  and
  partly computing of direct limit along the morphisms
    $(V,y)\to (Y,y)$ that are 
    etale neighbourhoods of the closed point $y$
    $$
    \varinjlim\limits_{(V,y)}
      \bF(V) 
    \simeq
    \bF(Y^h_y)
    $$ 
we get that
  \begin{multline*}\label{WtrLocHens}
  Wtr_{Pre}(\bF)(U^h) \simeq
  \varinjlim\limits_{{f\colon Y_1\to Y_2}}
    \bF(Y_2)\otimes WCor(U^h,Y_1)
  \simeq\\\simeq
  \varinjlim\limits_{{f\colon Y_1\to Y_2} ,\, {z\in x\times Y_1} }
    \bF(Y_2)\otimes WCor_z(U^h,Y_1)
  \simeq
  \varinjlim\limits_{{f\colon V^h_1\to Y_2},\, {z\in x\times Y_2} }
    \bF(Y_2)\otimes WCor_z(U^h,V^h_1)
  \simeq\\\simeq
  \varinjlim\limits_{{f\colon V^h_1\to V^h_2},\, {z\in x\times y} }
    \bF(V^h_2)\otimes WCor_z(U^h,V^h_1)
  \simeq
  \varinjlim\limits_{{f\colon V^h_1\to V^h_2} }
    \bF(V^h_2)\otimes WCor(U^h,V^h_1)
    .
  \end{multline*}
where 
$V^h_1$ and $V^h_2$ denoted local Hensel schemes,
$x$ denotes closed point of $U^h$
and $y$ denotes the projection of $z$ to $Y_2$
or the closed point of $V^h$
.
Since 
  presheaves $\bF_i$ are direct sums of representable ones
  that are projective objects of category  $\DPre$,
  and since
  $\bF^\bullet(U^h)$ is acyclic complex,
complex $\bF^\bullet(U^h)$ is contractable. 
Then 
since tensor product and direct limit preserves direct sums
by equality  \eqref{WtrLocHens}
$Wtr_{Pre}(\bF^\bullet)(U^h)$ is acyclic.   

{\em Theorem \ref{WtrNTopEx} is proved.}
   
\begin{remark}
In the last prove
to reduce to the case of the
  complexes of direct sums of representable presheaves 
we use that 
  we work with the category $D^-$ and
  the existence of resolutions. 
But in fact 
  it is not essential
because 
  even in $D^b(PreWtr)$
  the subcategory  of Nisnevich acyclic complexes 
  is generated by the complexes
    $$0\to 
      \bZ(\tilde{U})\to \bZ(U)\oplus\bZ(\tilde{X})\to \bZ(X)
      \to 0$$
   corresponding to the elementary Nisnevich squares
   $$\xymatrix{
   \tilde{U}\ar[r]\ar[d] & \tilde{X}\ar[d] &\\
   U\ar[r] &X   
   &.}$$
\end{remark}   

\vspace{10pt} 
{\em Proof of the theorem \ref{ReprChWtr}.}

Let's now 
substitute in the isomorphism of 
  adjacency \eqref{ShNDWtrAdj}
  between $L(Wtr_{ShN})$ and $F_{Wtr}$
the Nisnevich sheaf associated to 
  the presheave $\bZ(X)$ 
    for any smooth affine $X$,
  that is regarded as complex concentrated  in degree zero, 
and
sheave with $Witt$-transfers $\bF$ as
  a complex concentrated in degree $-i$. 

We get natural isomorphism
$$Hom_{D^-(ShN)}(\bZ(X),\bF[i])\simeq 
  Hom_{\DShNWtr}(L(Wtr_{ShN})(\bZ(X)),\bF[i]).$$
Since $\bZ(X)$ is representable,
then
$L(Wtr_{Pre})(\bZ(X))\simeq \bZ_{Wtr}(X)$,
and by theorem \ref{WtrNTopEx} 
$$L(Wtr_{Sh})(\bZ(X))\simeq L(Wtr_{Pre})(\bZ(X)).$$

Thus the adjacency \eqref{ShNDWtrAdj}
  gives us that
$$Hom_{D^-(ShN)}(\bZ(X),\bF[i])\simeq Hom_{D^-(\ShNW)}(\bZ(X),\bF[i])$$
and in combination with isomorphisms of 
    $Hom$-groups in derived categories, 
    $Ext$-groups 
    sheaf cohomology groups
$$H^i_{Nis}(X,\bF)\simeq Ext^i_{ShNis}(\mathbb Z(X),\bF) \simeq 
  Hom_{D^-(ShN)}(\bZ(X),\bF[i]),$$
this gives us the required isomorphism from theorem \ref{ReprChWtr}. 

{\em Theorem \ref{ReprChWtr} is proved.}

\section{The construction of category of effective $Witt$-motives $DWM^-_{eff}(k)$}
\label{sec_DWMeff}

\newcommand{\iaffDPW}{D^-_{\aff-inv}(PreWtr)}
\newcommand{\caffDPW}{D^-_{\aff-contr}(PreWtr)}
\newcommand{\iaffDSNW}{D^-_{\aff-inv}(PreWtr)}
\newcommand{\caffDSNW}{D^-_{\aff-contr}(PreWtr)}

\newcommand{\Dinv}{}
\newcommand{\Dcontr}{}

\newcommand{\DWMeff}{ DWM^-_{eff}(k)}
\newcommand{\DWMeffl}{ DWM^{-,l}_{eff}(k)}
\newcommand{\DWMeffr}{ DWM^{-,r}_{eff}(k)}

The construction of the category $\DWMeff$
is based on the theorems \ref{DHI} and \ref{ReprChWtr}
and is similar to the construction of $DM^-_{eff}$.

The category $\DWMeff$ can be defined in two ways, 
in a few informal words
   as localisation $\DWMeffl$ of the category $\DShNWtr$ 
     by $\affl$-homotopy equivalences 
   and
   as full subcategory $DShNWtr$ of $\DShNWtr$
     consists of homotopy invariant objects. 
And for the first one
  there is a functor
  $p\colon \DShNWtr\to \DWMeffl$
and for the second one
  there is a functor
  $i\colon \DWMeffr\to DShNWtr$.

One of the basic supposed properties of the category of motives
is that
the functor  
  from the category of varieties 
  to the category of motives
  is an universal cohomology theory
    in that sense that
      any cohomology theory on category of varieties
      can be passed  throw it. 
So
more natural definition 
(in respect to this universal property) 
is the first one ($\DWMeffl$).
However, 
to compute represented in this category cohomology theories
in terms of $Hom$-groups of  $\DShNWtr$,
means to define right adjoin to the functor $p$.
And the second definition ($\DWMeffr$) just provides computation of this functor
as the functor of full embedding $i$.    

The equivalence of this two definition
can be proved by showing that
  structure of category with interval 
  on the category of varieties 
  (where interval is affine line)  
  induce
  semi-orthogonal decomposition of
  the category $\DShNWtr$.
That is decomposition to
    $\affl$-contractable and 
    $\affl$-homotopy invariant parts. 
And both definitions of $\DWMeff$ gives exactly
  homotopy invariant part of $\DShNWtr$. 
  
As was mentioned in the introduction
the category of motives
combines in some sense 
  the structure of category with interval
  and
  topological structure on $Sm_k$
    (in this case Nisnevich topology).
Ability of a good combining of this structures
is provided by their coherence
  proved in theorem \ref{DHI}. 
Therefore
  firstly we prove 
    that affine line as interval in $Sm$
    induces semi-orthogonal decomposition,
    of derived category of 
      presheaves with $Witt$-transfers $\DPreWtr$,  
     that doesn't deals with topology, 
  and then we push down 
    this decomposition to the 
    derived category of sheaves $\DShNWtr$.

As a result we get following commutative diagram 
  $$\xymatrix{	  
  &&& D^-(\ShNW_k)\ar@{-->}_{l_{\aff}}[rd]&\\
  Sm_k\ar[r]&  Witt_k\ar[r]^{com\quad} & 
      D^-(PreWtr_k)
      \ar_{l_{NIs}}[ru]
      \ar@{-->}_{l_{\aff}}[rd]& 
                               & DWM^-_{eff}(k)& \\
  &&& D^-_{\aff}(PreWtr_k)\ar@{-->}_{l_{Nis}}[ur]
  }$$
where 
  the category $D^-_{\aff}(PreWtr_k)$ 
  is homotopy invariant part of the category $\DPreWtr$, 
  and functors $l_{\aff}$ and $l_{Nis}$ 
  can be regarded as localisation functors
    by the morphisms corresponding to projections   
      $X\times\affl\to X$ 
        for all smooth affine $X$
    and by the morphisms
    $\tilde{X}\to Cone(U \coprod \tilde{X}\to X) $ 
        where $X$ , $U$ , $\tilde{X}$, $\tilde{U}$ 
        are vertices of elementary Nisnevich square
    respectively.
The second arrow, i.e. $com $ is the functor that sands 
  any variety $X$ 
  to 
  a complex concentrated in degree zero 
    and defined by the presheave $\Zwtr(X)$.

Now we proceed to prove of mentioned above
  semi-orthogonal decompositions of categories 
  $\DPreWtr$ and $\DShNWtr$. 

\newcommand{\bA}{\mathcal{A}}
\newcommand{\bB}{\mathcal{B}}
\newcommand{\bC}{\mathcal{C}}

Let's give following definition of
semi-orthogonal decompositions in triangulated categories. 
\begin{definition}\label{def-SemiDecTr}
A semi-orthogonal decomposition
  $<\mathcal A,\mathcal B>$
  of some triangulated category $\mathcal C$
is
  the pair of two 
    full triangulated subcategories 
    $\mathcal A$ and $\mathcal B$ 
  that are semi-orthogonal and generates the category $\mathcal C$,
  i.e. such that  
  $$Hom_{\mathcal C}(B^\bullet,A^\bullet)=0$$ for all
    $B^\bullet\in\mathcal B$ and $A^\bullet\in\mathcal A$,
  and
  for any $C^\bullet\in \mathcal C$
  there is a 
  the distinguished triangle 
  $$A^\bullet[-1]\to B^\bullet\to C^\bullet\to A^\bullet$$    
  with $B^\bullet\in\mathcal B$ and $A^\bullet\in\mathcal A$.
\end{definition}

\begin{remark}\label{rem_semi-ort-dec}
For any semi-orthogonal decomposition 
$\mathcal C = <\mathcal A,\mathcal B>$
  two pairs of adjoin fuctors 
  $$\begin{aligned}
  \mathcal B\; 
    {\substack{i_{\mathcal B}\\ \rightleftarrows\\ l_{\mathcal A} }} 
    \;&\mathcal C\; 
      {\substack{l_{\mathcal B}\\ \rightleftarrows\\ i_{\mathcal A} }} 
     \;\mathcal A
  \\   
  i_{\mathcal B}\dashv l_{\mathcal A},&\quad l_{\mathcal B}\dashv i_{\mathcal A}
  \end{aligned}$$
  where 
    $i_{\mathcal B}$ and $i_{\mathcal A}$ 
    are embedding functors 
      of subcategoried $\mathcal A$ and $\mathcal B$
  and
    $l_{\mathcal A}$ and $l_{\mathcal B}$
    are equivalent to localisation functors 
      at subcategories $\mathcal A$ and $\mathcal B$.
  And
  the compositions 
  $pr_{\mathcal A} = i_{\mathcal A} \colon l_{\mathcal A}$ and
  $pr_{\mathcal B} = i_{\mathcal B} \colon l_{\mathcal B}$
  are endofunctors-projectors on $\mathcal C$, 
  i.e. 
  there are natural transformations
  $\epsilon_{\mathcal A}\colon id_{\mathcal C}\to pr_{\mathcal A}$
  $\xi_{\mathcal B}\colon pr_{\mathcal B}\to id_{\mathcal C}$
  such that 
  $\epsilon(A^\bullet)$ is isomorphism 
    for any $A^\bullet\in\mathcal{A}$  
  and
  $\xi(B^\bullet)$ is isomorphism for any $B^\bullet\in\mathcal{B}$.
\end{remark}

Also we will use following standard fact
\begin{lemma}\label{DirFactor} 
Let $\bA$ be abelian category with (infinite) direct sums
and $\bB$ be its full abelian subcategory 
  closed under (infinite) direct sums. 

Then
  thick subcategory of $D^-(\bA)$
  generated by $\bB$
  (i.e. by complexes concntrated in degree zero 
    whose zero term is an object of $\bB$)
  equal to $\ker(D^-(\bA)\to D^-(\bA/\bB) )$
  and equal to subcategory of complexes that
    cohomology lies in $\bB$. 
    
And all complexes 
consisting of objects lying in $\bB$
(i.e. lying in $D^-(\bB)$)
lies in this category .    
\end{lemma}

Now let's give also one useful definition 

\begin{definition}\label{def_invcontrObj}
Let $\bC$ be additive category
and
$S\colon C\to C$ be endo-functor
with  
natural transformation
  $p\colon Id_\bC\to S$ and $s_0,s_1\colon S\to Id_\bC$
  such that $s_1\circ p = id_S = s_1\circ p $.
The object 
$C\in bC$ is called 
  \textbf{$S$-invariant  }if
  $p(C)$ is isomorphism,
and it is called 
  \textbf{$S$-contractable} if
  there is a morphism $h\colon C\to S(C)$ such that 
  $s_0(C)\circ h = id_C$ and $ s_0(C)\circ h = 0$.  
\end{definition}

\vspace{5pt}

The semi-orthogonal decompositions 
  on categories $\DPreWtr$ and $\DShNWtr$
are induced in some sense by the 'action'
of affine line on the category of varieties and 
  categories  $WCor$, $PreWtr$, $\DPreWtr$ and $\DShNWtr$. 
  
\begin{definition}\label{def_afficObj}
Let $Sm\to \bC$ be any additive category 
  under the category of smooth affine varieties
  with internal $\bHom$-functor. 
Then affine  line $\affl$,
projection homomorphism $p\colon \affl\to pt$,
and morphism of zero and unit sections $s_0,s_1\colon pt\to \affl$
considered as objects and morphisms in the category $\DPreWtr$
induce by applying of internal $\bHom$-functor 
an following endo-functor and natural transformations
\begin{gather*}
S=\bHom_\bC(\affl,-)\colon\DPreWtr,\\
p=\bHom_\bC(p,-)\colon Id_\DPreWtr\to S,\quad
s_0,s_1 = \bHom_\DPreWtr(p,-)\colon S\to Id_\DPreWtr
.\end{gather*}

And definition \ref{def_invcontrObj}
provides a notions of 
$\affl$-homotopy invariant and
$\affl$-contractable 
objects of $\bC$.
\end{definition}

Finally let's give following standard definition.
\begin{definition}
Let denote by $\Delta$
  the simplicical scheme 
  with 
\begin{gather*}
\Delta^n = Spec k[x_0,x_1,\dots,x_n]/(x_0+x_1+\dots x_n -1) \\
e_{n,i} \colon
  (x_0,x_1,\dots,x_n)\mapsto (x_0,\dots,x_i,0,x_{i+1},\dots, x_n)\\
d_{n,i}\colon
  (x_0,x_1,\dots,x_n)\mapsto (x_0,\dots,x_i+x_{i+1},\dots, x_n)
.
\end{gather*}
\end{definition}

\begin{theorem}\label{affsodDPW}
There is
semi-orthogonal decomposition of the category $\DPreWtr$
  $$\DPreWtr = <\caffDPW,\iaffDPW>$$
such that 
  $\caffDPW$ as full subcategory of  $\DPreWtr$ 
  consists of complexes quasi-isomorphic to 
    complexes of $\affl$-contractable presheaves 
  and  
  $\iaffDPW$ consists of complexes
    with homotopy invariant cohomologies.
\end{theorem}
{\em Proof of the theorem.}
Let's consider two full triangulated subcategories of $\DPreWtr$:
subcategory $\bA$ consisting of 
  homotopy invariant objects of $\DPreWtr$
  in sense of definition \ref{def_afficObj}
and 
subcategory $\bB$ that is thick subcategory generated by
  $\affl$-contractable objects.  
And let's show that this categories provides 
semi-orthogonal decomposition of $\DPreWtr$.

First let's show that 
$Hom_{\mathcal C}(B^\bullet,A^\bullet)=0$ for all
    $B^\bullet\in\mathcal B$ and $A^\bullet\in\mathcal A$.
In fact,
  $$s_0(A^\bullet)= s_1(A^\bullet)\colon S(A^\bullet)\to A^\bullet,$$
  since $A^\bullet$ is invariant object and
  $p(A^\bullet)$ induce isomorphism of $A^\bullet$ and $S(A^\bullet)$.
And since
  $B^\bullet$ is $\affl$-contractable object,
  there is an $\affl$-homotopy 
  $$
  h \in Hom_\DPreWtr(B^\bullet, \bHom(\affl,B^\bullet)\colon\quad
  s_0(B^\bullet)\circ h = id,\; s_1(B^\bullet)\circ h = 0
  .$$
So
  for any $f\in \bHom(B^\bullet,A^\bullet)$
  $$
  f =
  f\circ s_1(B^\bullet)\circ h =
  s_0(A^\bullet)\circ S(f)\circ h =
  s_1(A^\bullet)\circ S(f)\circ h = 
  f\circ s_0(B^\bullet)\circ h
  = 0
  .$$
$$\xymatrix{
\bHom(\affl,B^\bullet)
  \ar@<-1ex>[d]_{s_1}\ar@<1ex>[d]^{s_0} \ar[r]^{S(f)}&
\bHom(\affl,A^\bullet)
  \ar@<-1ex>[d]_{s_1}\ar@<1ex>[d]^{s_0} \\
B^\bullet
  \ar[u]|h \ar[r]^f &  
A^\bullet  
}$$  

Next let's consider the functor 
$$C^* = Hom_{\DPreWtr}(\Delta^\bullet, -)\colon \DPreWtr\to \DPreWtr.$$          
The canonical morphism of simplicial objects 
$\Delta^\bullet\to pt^\bullet\to pt$ 
  where $pt^\bullet$ denotes
    constant simplicial object
  and $pt$ simplicial object 
    concentrated at the degree zero
induces natural transformations 
$$C^*\xrightarrow{\epsilon} 
Hom_{\DPreWtr}(pt^\bullet, -)\to 
Hom_{\DPreWtr}(pt^\bullet, -)\simeq Id_{\DPreWtr}.$$
Let 
$\epsilon^\prime \colon C^*\to Id_{\DPreWtr}$ 
denotes the composition.
Then for any complex $C^\bullet$ there is a distinguished triangle
$$ Cone(\epsilon)[1]\to
     C \xrightarrow{\epsilon^\prime} 
       C^*(C) \to 
          Cone(\epsilon). $$          
Standard simplicial partition of 
  the cylinders 
  $\Delta^i\times\affl$
defines the $\affl$-homotopy between 
  zero and unit section
  $$s_0,s_1\colon C^*(C)\to \mathcal Hom(\affl,C^*(C))$$ 
that shows that $C^*(C)$ is homotopy invariant.
On other side 
\begin{multline*}
Cone(\epsilon) = \\
(
  \cdot\to
  \bHom_{PreWtr}(
     \bZ_{Wtr}(\Delta^i)/ \bZ_{Wtr}(pt),
     C
   )
  \to\cdots\to
  \bHom_{PreWtr}(
     \bZ_{Wtr}(\Delta^1)/ \bZ_{Wtr}(pt),
     C
   )
  \to 0 
)
.
\end{multline*} 
And its terms are contractable presheaves,
because
  linear homotopy of affine  simplexes 
  $$\begin{aligned}
  \affl\times\Delta^n&\rightarrow& \Delta^n
  \qquad\qquad\qquad\qquad\qquad\qquad
  \qquad\qquad \\
  (\lambda,x_0,x_1,\dots,x_n)&\mapsto&
    (x_0+\lambda\dot \sum_i x_i, 
      (1-\lambda)\dot x_1,\dots,(1-\lambda)\dot x_n
     )
  \end{aligned}$$
induce a homotopy of 
 $
 \bHom_{PreWtr}(
     \bZ_{Wtr}(\Delta^n)/ \bZ_{Wtr}(pt),
     C
   )
 . $  
Hence $Cone(\epsilon)\in \caffDPW$.
Any contractable presheave defines 
  the contractable object of $\DPreWtr$
so any complex consists 
  of contractable presheaves lies in $\caffDPW$.
So $\caffDPW\subset \mathcal B$.
Thus $Cone(\epsilon)\in \mathcal B$.

Thus we get semi-orthogonal decomposition 
  $$\bC = <\mathcal B,\mathcal A>.$$
In addition by the remark \ref{rem_semi-ort-dec} 
that fact that 
  the third term in the triangle 
  is the result of 
    the applying functor $C^*$,
implies that
  $C^*$ defines the functor form 
  $\DPreWtr$ to $\iaffDPW$ 
  that is
  left adjoin to the embedding of 
  $\iaffDPW$ into $\DPreWtr$.

To finish the proof 
it is enough to show that
$$\iaffDPW = \mathcal A, \caffDPW= \mathcal B.$$

To show the firs it is enough to note that   
the functor 
  $\bHom_{PreWtr}(\affl,-)$ is exact and
$$
h^i(bHom_{\DPreWtr}(\affl,C^\bullet))
  \simeq 
     \bHom_{PreWtr}(\affl,h^i(C^\bullet))
.$$
So the equalities 
$$
h^i(\bHom_{\DPreWtr}(\affl,C^\bullet))\simeq  h^i(C^\bullet)
$$
are equivalent to the equalities 
$$
\bHom_{PreWtr}(\affl,h^i(C^\bullet))\simeq  h^i(C^\bullet). $$ 
Since it is just proved that 
  $<\mathcal B,\mathcal A>$ is semi-orthogonal decomposition 
  of $\DPreWtr$ 
by remark \ref{rem_semi-ort-dec}
for any $B^\bullet\in \mathcal B$
the morphism $Cone(\epsilon(B^\bullet)[-1]\to B^\bullet$
is quasi-isomorphism.
But as was mentioned above 
  $Cone(\epsilon(B^\bullet)[-1]\in \caffDPW$.

{\em Theorem is proved.}

\renewcommand{\Dinv}{\iaffDPW}
\renewcommand{\Dcontr}{\caffDPW}
\begin{theorem}\label{genaffsodDPW}
The category $\Dinv$
is generated as triangulated subcategory
by homotopy invariant presheaves
  (considered as complexes concentrated in degree zero).
 
The category $\Dcontr$
is generated as triangulated subcategory
by presheaves $\bZ(\affl\times X)/\bZ(X)$
  (considered as complexes concentrated in degree zero) 
  for all smooth affine $X$.
\end{theorem}

{\em Proof of the theorem.}
The first statement is a particular case of the statement from
lemma \ref{DirFactor}. 

To prove the second
it is enough to find a 
  natural resolvent in $\DPreWtr$ 
  consisting of infinite direct sums of 
  presheaves $\bZ(\affl\times X)/\bZ(X)$,
because if such resolvent exists 
then any complex $\bB^\bullet$ consisting of contractable presheaves
is quasi-isomorphic to 
totalisation of bi-complex constituted 
by resolvents of $\bB^i$,
and this totalisation is 
a complex consisting of 
direct sums of terms of 
presheaves $\bZ(\affl\times X)/\bZ(X)$.
For any presheave $\bF$ there is 
a natural (in $\bF$) sequence 
\newcommand{\bP}{\mathcal P}
\begin{align*}\label{contrcover}
\xymatrix{
\cdots\ar[r]&
\bP_i  \ar[r]
  \ar[d]^{\epsilon^i} &
\bP_{i-1}  \ar[r]
  \ar[d]^{\epsilon^{i-1}} &
\cdots\ar[r] &
\bP_1  \ar[r]
  \ar[d]^{\epsilon^1} &
\bP_0  \ar[r]
  \ar[d]^{\epsilon} &
0    
\\
\cdots\ar[r]&
\bF_i\ar@{^(->}[ru] &
\bF_{i-1}\ar@{^(->}[ru] &
\cdots\ar[r]&
\bF_1\ar[ru] &
\bF\ar[ur]
}\\
\bP_i = 
  \sum\limits_{
    U,s\in \bF_i(\affl\times U)\colon j^U_0(s)=0
    }
  \bZ(\affl\times U)/\bZ(U),\\
j_0^X\colon 0\times X\hookrightarrow\affl\times X,\;
j_1^X\colon 1\times X\hookrightarrow\affl\times X,\\
{\epsilon^i}_{X,s}\colon
  \bZ(\affl\times X)/\bZ(X)\simeq \coker(j_0)\xrightarrow{s}\bF_i,\\
\bF_i=\ker(\epsilon^i),
\bF_0 = \bF \quad
.\end{align*}
(This sequence in fact defines 
 adjoin functor to the embedding functor of 
 subcategory in $\DPreWtr$ generated by presheaves
 $\bZ(\affl\times X)/\bZ(X)$). 
If $\bF\in \Dcontr$ then $C^*(\bF)$ is acyclic
and in particular $ h^0(C^*(\bF)) = 0$. 
But $h^0(C^*(\bF)) = \coker(\epsilon)$,
hence $\epsilon$ is surjective.
Next since 
  $\bZ(\affl\times X)/\bZ(X)\in \Dcontr$
  and
  $\epsilon$ is surjective,
$\bF_1 = Cone(\epsilon)[1]\in \Dcontr$.  
Then by induction we get that
all $\epsilon_i$ are surjective
(and all $\bF_i$ lies in $\Dcontr$).
Thus sequence \eqref{contrcover} gives resolvent of $\bF$.
{\em Theorem is proved.}

Let's proceed to the considering of the category $\DShNWtr$ 
and how semi-orthogonal decomposition consistent with Nisnevich topology.

\begin{theorem}\label{affsodDSNW}
There is a semi-ortogonal decomposition
  $$\DShNWtr = <\caffDSNW,\iaffDSNW>$$
such that
  the category $\caffDSNW$ 
  is thick subcategory generated 
  by Nisnevich sheaves 
    $\bZ_{Wtr,Nis}(\affl\times X)/\bZ_{Wtr,Nis}(X)$
  and 
  the category $\iaffDSNW$ 
  as full subcategory of $\DShNWtr$
  consists of 
    the complexes with homotopy invariant cohomologies.
\end{theorem}
{\em Proof of the theorem.}  
We can regard category $\DShNWtr$ as 
the localisation of $\DPreWtr$ 
at Nisnevich sheaf quasi-equivalences.
Let's denote this functor by 
  $$l_{Nis}\colon \DPreWtr \to \DShNWtr.$$
Let's 
  $\bA$ and $\bB$ be images of subcategories
  $\iaffDPW$ and $\caffDPW$ 
  $$
  \bA,\bB\subset \DShNWtr\colon\quad
  \bA = l_{Nis}(\iaffDPW),\;\bB = \caffDPW
  ,
  $$
  i.e.
  full subcategories that 
  consists of the complexes that are 
  sheaf quasi-isomorphic to the complexes 
  that lies in $\caffDPW$ and $\iaffDPW$
  respectively.
  
Let's show that 
\begin{equation}\label{eqimageAff_nis}
\bA = \iaffDSNW,\;\bB = \caffDSNW
.\end{equation}   
The functor $l_{Nis}$ sends
any presheave as a complex concentrated in degree zero
to its sheafication.
Then since localisation functor 
sends generators of thick subcategory to 
generators of its image,
theorem \ref{genaffsodDPW} implies that
subcategories $\bA$ and $\bB$ 
are thick subcategories generated 
  by sheafifications of homotopy invariant presheaves
  and by sheaves
    $\bZ_{Wtr,Nis}(\affl\times X)/\bZ_{Wtr,Nis}(X)$
    (that are Nisnevich sheafications of presheaves
     $\bZ_{Wtr}(\affl\times X)/\bZ_{Wtr}(X)$
     by definition).
So we get the equality for $\bB$.
To prove the equality for $\bA$
let's note that
  by lemma \ref{DirFactor} $\iaffDSNW$ is thick subcategory
  generated by homotopy invariant sheaves. 
By theorem 4 form \cite{AD_WtrSh}
sheafification of homotopy invariant presheave
  with $Witt$-transfers
is homotopy invariant.
Conversely any homotopy invariant sheaf   
is homotopy invariant presheave.
So the set of 
Nisnevich sheafifications of homotopy invariant presheaves
with $Witt$-transfers
is exactly the set of homotopy invariant sheaves
with $Witt$-transfers. 
Thus we get the equalities \eqref{eqimageAff_nis}.

Then let's show that
theorem \ref{DHI} and theorem \ref{ReprChWtr} 
  (about 
    strictly homotopy invariance of 
     Nisnewich sheafication of 
      homotopy invariant presheave 
       with  $Witt$-transfers
   and about
     isomorphism of 
       $Ext$-groups in $\ShNW$ and     
       Nisnevich cohomology groups of sheave 
         with  $Witt$-transfers) 
implies that the categories $\bA$ and $\bB$ are semi-orthogonal.
Really
to prove that 
  $$
  Hom_{\DShNWtr}(B^\bullet,A^\bullet)=0\colon\quad
  A^\bullet\in \bA, B^\bullet\in\bB
  $$
it is enough 
to check it 
  on generators of this subcategories,
  i.e. 
  for 
  $$
  A^\bullet = \bF[i],\quad
  B^\bullet = \bZ_{Wtr,Nis}(X\times\affl)/\bZ_{Wtr,Nis}(X)
  $$
    for any homotopy invariant presheave $\bF$ and
    smooth affine $X$.
But
\begin{multline*}
Hom_{\DShNWtr}
  (\bZ_{Wtr,Nis}(X\times\affl)/\bZ_{Wtr,Nis}(X),
   \bF[i])
\simeq \\
Ext_{ShNisWtr}^i
  (\bZ_{Wtr,Nis}(X\times\affl)/\bZ_{Wtr,Nis}(X),
   \bF_{Nis}) = 0 
\end{multline*}   
and the last group is zero because by 
  the remark to the theorems \ref{DHI} and \ref{ReprChWtr}
\begin{multline*}
Ext_{ShNisWtr}^i(\bZ_{Wtr,Nis}(X\times\affl),\bF_{Nis}) \simeq
H^i_{Nis}(X\times\affl,\bF_{Nis}) \simeq \\
H^i_{Nis}(X,\bF_{Nis}) \simeq 
Ext_{ShNis}^i(\bZ_{Wtr,Nis}(X),\bF_{Nis})
.
\end{multline*}

Now 
to prove that the pair $<\bB,bA>$ 
provides semi-orthogonal decomposition of $\DShNWtr$
it is enough to show that 
  for any object $C^\bullet\in \DShNWtr$
  there is distinguished triangle 
  $$
  A^\bullet[-1]\to B^\bullet\to C^\bullet \to A^\bullet\colon\quad
  B^\bullet\in\bB,\;A^\bullet\in\bA
  .$$
But since 
  $<\caffDPW,\iaffDPW>$ is semi-orthogonal decomposition 
    of $\DPreWtr$ 
for any complex $C$ there is 
distinguished triangle in $\DPreWtr$
  $A^\bullet[-1]\to B^\bullet\to C^\bullet \to A^\bullet$
  with $B^\bullet\in\caffDPW$ and $A^\bullet\in\iaffDPW$
and   
this triangle remains to be distinguished after localisation $l_Nis$.

{\em Theorem is proved.}  

\begin{definition}
The category of effective motives $\DWMeff$
is homotopy invariant part of $\DShNWtr$,
i.e. the category $\iaffDSNW$.
\end{definition}

\begin{remark}\label{adjAffDSh}
Due to semi-orthogonal 
  decomposition proved in theorem \ref{affsodDSNW}
there is adjacency 
\begin{equation}\label{adj_DShNW-iaff}\begin{aligned}
C^*\colon \DShNWtr &\rightleftarrows& \DWMeff\colon i_{\aff} 
&\; &\\
C^*&\dashv&i_{\aff}\qquad\qquad\qquad\qquad
&.\;& 
\end{aligned}\end{equation}
of the projection functor defined by $C^*$
and embedding. 

\end{remark}

\begin{remark}
Since 
  the projection functor $\DShNWtr\to \iaffDSNW$ 
  is equivalent 
  to the localization functor 
    at 
     the subcategory $\caffDSNW$,
and since 
  $\caffDSNW$ by the theorem \ref{affsodDSNW} 
  is generated by cones of projection
    $X\times\affl\to X$ for all $X$,
then
the category $\DWMeff$ 
is equivalent to 
the localisation of $\DShNWtr$ 
 at $X\times\affl\to X$.
\end{remark}

Let's define $Witt$-motives of smooth affine varieties.
\begin{definition}
The functor $$WM\colon Sm_k\to \DWMeff$$ 
defineing  
  $Witt$-motives of smooth affine varieties
is composition
$$DWM \stackrel{def}{=} 
  l_\aff\circ \tilde{-_{Nis}}\circ Wtr\circ \bZ(-).$$
\end{definition}

Now 
we get the required property 
by composing of the 
adjacency \eqref{adj_DShNW-iaff} 
with
adjacency \eqref{ShNDWtrAdj} from previous paragraph. 

\begin{theorem}
\label{RepH}
There is natural isomorphism
$$ Hom_{\DWMeff}(WM(X),\bF[i]) \simeq H^i_{Nis}(X,\bF) $$
for all smooth affine $X$ and homotopy invariant sheaf with $Witt$-transfers $\bF$.
\end{theorem}
{\em Proof of the theorem.}
By definition of the functor $WM$ 
  $$WM(X) = C^*(\tilde{Witt_{Nis}}(X)).$$
Due to adjacency of $C^*\dashv i_{\aff} $ from the remark \ref{adjAffDSh}
$$Hom_{\DWMeff}(C^*(\mathbb Z_{Wtr,Nis}(X)),\bF[i]) \simeq
  Hom_{\DShNWtr}(\mathbb Z_{Wtr,Nis}(X),\bF[i])$$    
Then by isomorphism from theorem \ref{ReprChWtr}
  $$Hom_{\DShNWtr}(\mathbb Z_{Wtr,Nis}(X),\bF[i]) \simeq 
    Ext_{NisSh}(\mathbb Z_{Wtr,Nis}(X),\bF[i]), $$
And finally
$$Ext_{NisSh}(\bZ(X),\bF[i])\simeq H^i_{Nis}(X,\bF).$$  
{\em Theorem is proved.}

\end{document}